\documentclass[aihp]{imsart}

\RequirePackage{amsthm,amsmath,amsfonts,amssymb}
\RequirePackage[numbers]{natbib}
\RequirePackage{graphicx}%% uncomment this for including figures

\startlocaldefs
\def\build#1_#2^#3{\mathrel{
\mathop{\kern 0pt#1}\limits_{#2}^{#3}}}
\def\llbracket{[\hspace{-.10em} [ }
\def\rrbracket{ ] \hspace{-.10em}]}

\theoremstyle{plain}
\newtheorem{theorem}{Theorem}
\newtheorem{proposition}[theorem]{Proposition}
\newtheorem{lemma}[theorem]{Lemma}

\theoremstyle{remark}

\def\w{\mathrm{w}}

\def\u{\mathcal{U}}
\def\t{\mathcal{T}}
\def\k{\mathcal{K}}

\def\W{\mathcal{W}}
\def\S{\mathcal{S}}
\def\T{\mathbb{T}}
\def\N{\mathbb{N}}
\def\M{\mathbb{M}}

\def\D{\mathbb{D}}
\def\P{\mathbb{P}}
\def\U{\mathbb{U}}
\def\E{\mathbb{E}}

\def\R{\mathbb{R}}
\def\z{\mathcal{Z}}
\def\y{\mathcal{Y}}
\def\n{\mathcal{N}}
\def\m{\mathcal{M}}
\def\cc{\mathcal{C}}
\def\ee{\mathcal{E}}
\def\ve{{\varepsilon}}
\def\la{\longrightarrow}
\def\da{\downarrow}
\def\ov{\overline}

\def\dd{\mathrm{d}}
\def\wh{\widehat}
\def\wt{\widetilde}
\def\tr{\mathrm{tr}}

\def\bbf{\mathbf{f}}

\def\Z{\mathbb{Z}}

\def\v{\mathcal{V}}
\def\bb{\mathbf{b}}
\def\xx{\mathbf{x}}

\def\dg{{\rm d_{gr}}}

\def\be{\mathbf{e}}
\def\bn{\mathbf{n}}
\def\bg{\mathbf{g}}

\endlocaldefs

\begin{document}

\begin{frontmatter}

\title{Spatial Markov property in Brownian disks}
%\title{A sample article title with some additional note\thanksref{T1}}
\runtitle{Spatial Markov property in Brownian disks}
%\thankstext{T1}{A sample of additional note to the title.}

\begin{aug}
\author[A]{\inits{J.}\fnms{Jean-Fran\c cois}~\snm{Le Gall}\ead[label=e1]{jean-francois.le-gall@universite-paris-saclay.fr}},
\author[B]{\inits{A.}\fnms{Armand}~\snm{Riera}\ead[label=e2]{riera@lpsm.paris}}
%%%%%%%%%%%%%%%%%%%%%%%%%%%%%%%%%%%%%%%%%%%%%%
%% Addresses                                %%
%%%%%%%%%%%%%%%%%%%%%%%%%%%%%%%%%%%%%%%%%%%%%%
\address[A]{Universit\'e Paris-Saclay, Math\'ematiques, B\^atiment 307, 91405 ORSAY, FRANCE\printead[presep={,\ }]{e1}}

\address[B]{Sorbonne Universit\'e, LPSM, 4 Place Jussieu, 75005 PARIS, FRANCE\printead[presep={,\ }]{e2}}
\end{aug}

\begin{abstract}
We derive a new representation of the Brownian disk in terms of a forest of labeled trees, where 
labels correspond to distances from a subset of the boundary. We then use this representation 
to obtain a spatial Markov property showing that the complement of a hull centered at
a boundary point of a Brownian disk is again a Brownian disk, with a random perimeter,
and is independent of the hull conditionally on its perimeter. Our proofs rely in part on a study of the peeling
process for triangulations with a boundary, which is of independent interest. The results 
of the present work will be applied to a continuous version of the peeling
process for the Brownian half-plane in a companion paper.
\end{abstract}

\begin{abstract}[language=french]
Nous donnons une nouvelle repr\'esentation du disque brownien en termes d'une for\^et d'arbres
\'etiquet\'es, o\`u les \'etiquettes correspondent aux distances depuis une partie de la
fronti\`ere. Nous utilisons cette repr\'esentation pour obtenir une propri\'et\'e de Markov
spatiale montrant que le compl\'ementaire d'une boule compl\'et\'ee centr\'ee en 
un point de la fronti\`ere est encore un disque brownien, avec un p\'erim\`etre
al\'eatoire, et est conditionnellement \`a ce
p\'erim\`etre ind\'ependant de la boule compl\'et\'ee. Les preuves reposent en partie sur une \'etude du processus d'\'epluchage
pour des triangulations avec fronti\`ere, qui est d'int\'er\^et ind\'ependant.  Les r\'esultats du pr\'esent 
travail seront appliqu\'es dans un article suivant \`a une version continue du processus d'\'epluchage pour le demi-plan brownien.
\end{abstract}



\end{frontmatter}

\section{Introduction}
Brownian disks are basic models of random geometry that arise as scaling limits
of random planar maps with a boundary, in the regime where the number 
of faces grows like the square of the boundary size. Brownian disks first appeared
in the work of Bettinelli \cite{Bet}, who obtained the existence of
subsequential limits of rescaled 
quadrangulations with a boundary in the Gromov-Hausdorff sense. The uniqueness
of the limit, which is called the Brownian disk, was then obtained in the
work \cite{BM} of Bettinelli and Miermont. In these scaling limits, \cite{Bet} and \cite{BM} mainly deal with
the case where both the boundary size and the volume are fixed, but it is
also of interest to study the so-called free Brownian disk for which the boundary 
size (also called the perimeter) is fixed but the volume is random, cf.~Section 1.5 in \cite{BM}. The free 
Brownian disk then appears as the limit of Boltzmann distributed random 
quadrangulations with a boundary, and in fact of much more general bipartite
planar maps \cite[Theorem 8]{BM}. In view of certain applications, it is desirable
to consider the case of random planar maps with a simple boundary. Convergence 
to the free Brownian disk in that case was obtained for quadrangulations by Gwynne and Miller \cite{GM1}
and for triangulations by Albenque, Holden and Sun \cite{AHS}. Both these papers prove 
convergence in a strong form of the Gromov-Hausdorff topology, which they
call the GHPU convergence, which includes the convergence of the so-called boundary curves.

In the present work, we are primarily interested in the free Brownian disk, and we
also consider the variant called the free pointed Brownian disk, where there is
a distinguished point in the interior of the disk --- if one ``forgets'' this distinguished point,
the distribution of the free pointed Brownian disk becomes a size-biased version of
the distribution of the free Brownian disk. The Bettinelli-Miermont construction
applied to the free pointed Brownian disk (see Section \ref{Bet-Mie} below)
relies on a random forest made of a collection of labeled $\R$-trees,
where labels correspond, up to a shift, to distances from the distinguished point
of the Brownian disk. Different constructions, still based on labeled $\R$-trees, have
been proposed in \cite{Disks,Repre} and shed light on various properties
of (free) Brownian disks. In the construction of \cite{Disks}, labels correspond to
distances from the boundary, and in \cite{Repre} they represent distances from
a distinguished point of the boundary. 
In the present work, we give yet another representation of the free Brownian disk
(Theorem \ref{identif-hull-comple}), where labels correspond to distances
from a part of the boundary.

Let us briefly describe our new representation of the free Brownian disk. Let us fix $\xi>0$, and let $C(\R_+,\R_+)$
denote the space of all continuous functions from $\R_+$ into $\R_+$. Then
let $\sum_{i\in I} \delta_{(t_i,e_i)}$ be a Poisson point measure on $[0,\xi]\times C(\R_+,\R_+)$
with intensity $2\,\dd t\,\bn(\dd e)$, where $\bn(\dd e)$ denotes the It\^o measure
of positive excursions of linear Browian motion. It is well known that each
excursion $e_i$ codes a compact $\R$-tree, which is denoted by $\t_{e_i}$. 
We then assign a real label $\ell_a$  to every point $a$ of $[0,\xi]$ and to every  $a\in\t_{e_i}$, $i\in I$,
in the following way. We consider a Brownian excursion $(\be_t)_{0\leq t\leq \xi}$ of duration $\xi$
and we declare that the label of every $a\in[0,\xi]$ is $\ell_a:=\sqrt{3}\,\be_a$. For every 
$i\in I$, we assign the label $\sqrt{3}\,\be_{t_i}$ to the root of $\t_{e_i}$, and then we require that
labels evolve like linear Brownian motion along the segments of $\t_{e_i}$ (independently
when $i$ varies). In other words, $(\ell_a)_{a\in\t_{e_i}}$ is distributed as Brownian 
motion indexed by $\t_{e_i}$ started from $\sqrt{3}\,\be_{t_i}$ at the root. 
The preceding objects (the normalized Brownian excursion $(\be_t)_{0\leq t\leq \xi}$ and
the labeled trees $\t_{e_i}$) are the basic ingredients of the Bettinelli-Miermont
construction of free Brownian disks \cite{Bet,BM} (see Section \ref{Bet-Mie} below), but
here we perform an additional step: we prune each tree $\t_{e_i}$ at levels where labels first hit the value $0$,
and write $\wt\t_{e_i}$ for the resulting pruned tree, which now carries nonnegative labels.
We then consider the union
$$\mathfrak{T}^\star:=[0,\xi] \cup\Bigg(\bigcup_{i\in I} \wt\t_{e_i}\Bigg),$$
where we identify the root of $\wt\t_{e_i}$ with the point $t_i$ of $[0,\xi]$. We then proceed in a way
very similar to known constructions of the Brownian sphere or the Brownian disk.
Precisely, if $a,b$ are two points of $\mathfrak{T}^\star$ with positive labels, we
let $D_\star^\circ(a,b)$ be the sum of the labels of $a$ and $b$ minus twice the 
minimal label ``between'' $a$ and $b$ (see Section \ref{conf-lim-space} for a more
precise definition) if this minimal label is positive, and if not we set $D_\star^\circ(a,b)=+\infty$.
We finally write $D_\star$ for the maximal pseudo-metric on $\mathfrak{T}^\star$
that is bounded above by $D^\circ_\star$. Theorem \ref{identif-hull-comple} then shows
that the quotient space\footnote{The notation $\mathfrak{T}^\star/\{D_\star=0\}$ refers to the quotient space of $\mathfrak{T}^\star$ for the equivalence relation defined by setting $a\simeq_\star
b$ if and only if $D_*(a,b)=0$.} 
 $\U:=\mathfrak{T}^\star/\{D_\star=0\}$
equipped with the metric induced by $D_\star$ is a free pointed
Brownian disk with (random) perimeter $\xi+\z$, where $\z$ is a random variable measuring, in some sense, the quantity of 
points with zero label in $\mathfrak{T}^\star$. Moreover, labels on $\U$ (which are inherited from the
labels on $\mathfrak{T}^\star$) correspond to distances from the subset of the boundary that
is the image of the set of points of $\mathfrak{T}^\star$ with zero label under the canonical
projection from $\mathfrak{T}^\star$ onto $\U$. The complementary part of the
boundary is the image of $[0,\xi]$ under the canonical projection.

An important motivation for Theorem \ref{identif-hull-comple} came from an application
to the complement of hulls centered at a boundary point of a Brownian disk. Let $\D'$ be a free Brownian disk 
with perimeter $\xi$ and boundary $\partial\D'$. One can define a ``standard boundary curve'' $(\Gamma(t))_{t\in[0,\xi]}$
that starts from a point uniformly distributed on $\partial \D'$ and runs along the boundary $\partial\D'$ 
at ``uniform speed'' (see Section \ref{Bet-Mie}). Then let $\alpha$ and $\beta$ be distinct real numbers in $[0,\xi]$
and consider the two points $a$ and $b$ of the boundary defined by $a=\Gamma(\alpha)$
and $b=\Gamma(\beta)$. 
Fix $r>0$ and write $B_r$ for the ball of radius $r$ centered at $a$ in $\D'$.
Conditionally on the event where the distance between $a$ and $b$
is greater than $r$, one may consider the connected component of $\D'\backslash B_r$ that contains $b$,
and we denote this component by $\wh B^\circ_r$.
By definition, the hull of radius $r$ centered at $a$, relative to $b$,
is $B^\bullet_r:=\D'\backslash\wh B^\circ_r$. 
Then, $\wh B^\circ_r$ (or rather its closure $\wh B^\bullet_r$) equipped with the appropriate intrinsic distance, is again a free Brownian disk, now with a random 
perimeter (Theorem \ref{spatial-Markov}). Moreover, conditionally on its boundary size,
this free Brownian disk is independent of the hull $B^\bullet_r$ also equipped 
with an intrinsic distance (Theorem \ref{indep-hull}). These results can be interpreted as a spatial Markov property of
the free Brownian disk. Imagine that one starts exploring the free Brownian disk from the point $a$ of the boundary. At the time
where one has discovered the ball $B_r(a)$ and all connected components of $\D'\backslash B_r(a)$ not containing $b$, 
what remains to be explored is again a free Brownian disk 
(with a random perimeter) which conditionally on its boundary size is independent
of what has already been discovered. This is also reminiscent of the peeling explorations
of random planar maps, which have found a number of striking applications
(see in particular \cite{AC,Cur,CLG}). 

The preceding results  take an even nicer form in the model 
called the Brownian half-plane \cite{CC,GM0,BMR}, which will be studied in the companion paper
\cite{ALG3}. The Brownian half-plane $\mathfrak{H}$
is a random non-compact metric space, which is homeomorphic to the usual half-plane
$\R\times\R_+$, so that it makes sense to define its boundary $\partial \mathfrak{H}$. The 
Brownian half-plane comes with a distinguished point $\xx$ on its boundary. For every
$r>0$, the hull of radius $r$ centered at $\xx$ is defined as the complement of the unbounded connected component
of $\mathfrak{H}\backslash B_r(\mathfrak{H})$, where $B_r(\mathfrak{H})$ denotes the closed ball of radius $r$
centered at $\xx$. 
Let $B^\bullet_r(\mathfrak{H})$ denote this hull. Then, the closure of $\mathfrak{H}\backslash B^\bullet_r(\mathfrak{H})$ equipped with the intrinsic metric
(and pointed at a boundary point which can be chosen in a deterministic way from the hull $B^\bullet_r(\mathfrak{H})$) is again 
a Brownian half-plane, which furthermore can be shown to be independent of the hull $B^\bullet_r(\mathfrak{H})$. This property is again
a continuous analog of the peeling process of infinite half-planar planar maps (see in particular \cite{AC}). The proof, whose details will be
given in \cite{ALG3}, is based on a passage to the limit from Theorems \ref{spatial-Markov} and \ref{indep-hull}. Similarly, a passage to the limit from
Theorem \ref{identif-hull-comple} yields a simple new representation of the Brownian half-plane. 

Let us discuss the relation between our main results and previous work. The paper \cite{Disks} shows that connected components of the complement 
of a ball centered at a distinguished point in the Brownian sphere are independent Brownian disks conditionally on their volumes and boundary sizes (see also \cite[Theorem 9]{Stars}
for a closely related result). The results of \cite{Disks,Stars} can be used to verify the equivalence of the definitions of Brownian disks given in 
Bettinelli and Miermont \cite{Bet,BM} and in Miller and Sheffield \cite{MS} (we note that \cite{MS} was motivated by strong connections with 
Liouville quantum gravity,  where Brownian disks correspond to the so-called quantum disks, see in particular \cite{MS2}). The proofs of \cite{Disks} rely on the excursion theory of \cite{ALG} and on a representation of Brownian disks 
in terms of labeled trees, where labels correspond to distances from the boundary. Here we are interested in the complement of the ball
(more precisely, the hull) centered at a boundary point of the Brownian disk,  and the proof of Theorem \ref{spatial-Markov} below depends on a very different representation
of the Brownian disk, which is provided by Theorem \ref{identif-hull-comple}. Much of what follows is in fact devoted to the proof of the latter result, which is rather involved and
requires a number of new ingredients. We view both \cite{Disks} and the present work as steps towards a
general form of the spatial Markov property in Brownian geometry.

Let us finally comment on the proofs of our results. 
The proof of Theorem \ref{identif-hull-comple} relies on
discrete approximations. The underlying idea is to start from a free pointed Brownian disk $\D$, and to consider hulls centered at
the distinguished point $\xx_*$ (which is now a point of the interior of $\D$ and not of $\partial\D$ as above) relative to the
boundary $\partial\D$. More precisely, we consider the hull $H$ with a radius $r_0$ which is the  distance between $\xx_*$ and the boundary $\partial\D$, and we write $U$ for the complement of $H$ in $\D$. In other words, $U$ is the connected component of the complement of the ball of radius $r_0$ centered at $\xx_*$ that contains all the boundary $\partial \D$ but the single point
$\xx_0$ realizing the distance between $\xx_*$ and $\partial\D$. 
The Bettinelli-Miermont construction of $\D$ (Section \ref{Bet-Mie}) allows one to get a representation of $U$
in terms of a Brownian excursion $\be$ and a collection $(\wt\t_{e_i})_{i\in I}$ of labeled trees having 
exactly the distribution described above. On the other hand,
one proves that the completion of $U$ for the appropriate intrinsic metric is also a free Brownian disk whose boundary can be
viewed as the union of $\partial \D$ and $\partial H$, provided that the unique point $\xx_0$ of $\partial \D\cap\partial H$ is split into two points.   
This identification of the law of (the completion of) $U$ is the difficult part of the proof of Theorem \ref{identif-hull-comple} in Section \ref{pass-lim},
and, for this, we first obtain an analogous discrete result: we observe that, for a Boltzmann distributed pointed triangulation
with a simple boundary, the analog of the set $U$, which is conveniently defined via a particular version of the 
peeling process, is again a Boltzmann distributed triangulation with a random boundary size, and we use
properties of the peeling process to investigate the asymptotics of this boundary size.
This part of the argument relies on a study of the peeling process for a 
Boltzmann distributed pointed triangulation with a boundary (Section \ref{discrete}), which is of independent interest. 

In order to derive the spatial Markov property of Theorems \ref{spatial-Markov} and \ref{indep-hull},
we rely on  Theorem \ref{identif-hull-comple} and we also use the representation of the Brownian disk in \cite{Repre}.
In this representation, the Brownian excursion $(\be_t)_{0\leq t\leq \xi}$ is replaced by a five-dimensional
Bessel excursion $(\bb_t)_{0\leq t\leq \xi}$, and the Poisson collection $(\t_{e_i})_{i\in I}$ of labeled trees
is conditioned to have only positive labels. Furthermore labels now correspond to distances
from a (uniformly distributed) point of the boundary. Thanks to this last property, the complement of the hull
of radius $r$ centered at the distinguished point of the boundary (and relative to another 
fixed point) can be coded by the ``subexcursion'' $\bb^{(r)}$ of $\bb$ above level $r$ that straddles a given time of $[0,\xi]$, 
and by the labeled subtrees $\t_{e_i}$ for indices $i$ such that $t_i$ belongs to the time interval associated with
$\bb^{(r)}$, provided these subtrees are pruned at levels where labels first hit $r$. Under an 
appropriate conditioning, the pair consisting of $\bb^{(r)}$ and the collection of
pruned labeled trees (where labels are shifted by $-r$) has the same distribution as the pair $(\be,(\wt\t_{e_i})_{i\in I})$ considered above, 
provided $\xi$ is replaced by the quantity $\xi'$ which is the duration of $\bb^{(r)}$. This allows one to 
apply Theorem \ref{identif-hull-comple}  in order to obtain that the hull complement is again a Brownian disk.

The paper is organized as follows. Section 2 gives several preliminaries. In particular, we recall the formalism
of curve-decorated measure metric spaces, and the associated Gromov-Hausdorff-Prokhorov-uniform distance $\dd_{GHPU}$, 
which has been introduced in \cite{GM0} and also used in \cite{AHS}. 
Moreover, we recall basic facts about snake trajectories and the Brownian snake excursion measure, 
which provide a convenient setting to deal with our labeled trees. In Section 3, we discuss the peeling process
of triangulations with a boundary, whose scaling limit is known to be the Brownian disk \cite{AHS}.
In Section~4, we introduce the space $\U$, and we explain how this space can be identified 
with (the completion of) the complement of a hull centered at the distinguished point in the free pointed Brownian disk. Section 5,
which is the most technical part of the paper, is devoted to the proof of Theorem \ref{identif-hull-comple}
identifying $\U$ as a Brownian disk. The general idea
is to pass to the limit from the analogous discrete result for triangulations, but unfortunately this
passage to the limit requires a number of technicalities. Finally, Section 6 presents the proof
of Theorems \ref{spatial-Markov} and \ref{indep-hull}. 

\newpage

\tableofcontents

\section{Preliminaries}

\label{sec:preli}

\subsection{Convergence of metric spaces}
\label{sec:convmetric}

In this work, we will consider different notions of convergence of a sequence of compact metric spaces, 
which we briefly present in this section. A bipointed compact metric space $(E,d,x,x')$ is just a compact metric space $(E,d)$ 
given with an ordered pair $(x,x')\in E\times E$ of distinguished points. We write $\M^{GH\bullet\bullet}$ for the set of all isometry classes 
of bipointed compact metric spaces (two pointed compact metric spaces $(E_1,d_1,x_1,x'_1)$ and $(E_2,d_2,x_2,x'_2)$
are isometry equivalent if there is an isometry $\Phi$ from $E_1$ onto $E_2$ such that $\Phi(x_1)=x_2$ and 
$\Phi(x'_1)=x'_2$). 
We can equip $\M^{GH\bullet\bullet}$ with the bipointed Gromov-Hausdorff distance $\dd_{GH\bullet\bullet}$, which is defined
by setting
$$
\dd_{GH\bullet\bullet}\big((E_1,d_1,x_1,x'_1),(E_2,d_2,x_2,x'_2)\big)
:= 
\inf\Big\{ d^E_{\mathrm{H}}(\Phi_1(E_1),\Phi_2(E_2)) \vee d(\Phi_1(x_1),\Phi_2(x_2))\vee d(\Phi_1(x'_1),\Phi_2(x'_2))\Big\},
$$
where the infimum is over all isometric embeddings $\Phi_1:E_1\la E$ and $\Phi_2:E_2\la E$ of 
$E_1$ and $E_2$ into the same compact metric space $(E,d)$, and $d^E_{\mathrm{H}}$
is the usual Hausdorff distance between compact subsets of $E$. Then, $(\M^{GH\bullet\bullet},\dd_{GH\bullet\bullet})$
is a Polish space. See in particular \cite{BBI} (proofs in \cite{BBI} are given in the non-pointed case, but are
immediately adapted). We can also define $\dd_{GH\bullet\bullet}$ in terms of correspondences. Recall 
that a correspondence between $E_1$ and $E_2$ is a subset $\cc$ of $E_1\times E_2$ such that
the restrictions to $\cc$ of both canonical projections $E_1\times E_2 \la E_1$ and $E_1\times E_2 \la E_2$ 
are surjective. The distortion of $\cc$ is then defined by
$$\mathrm{dis}(\cc):= \sup\{|d_1(y_1,z_1)-d_2(y_2,z_2)|:(y_1,y_2)\in\cc, (z_1,z_2)\in\cc\},$$
and the $\dd_{GH\bullet\bullet}$ distance can be expressed as
$$\dd_{GH\bullet\bullet}((E_1,d_1,x_1,x'_1),(E_2,d_2,x_2,x'_2))= \frac{1}{2}\,
\inf\{\mathrm{dis}(\cc)\}.$$
where the infimum is over all correspondences between $E_1$ and $E_2$ such that $(x_1,x_2)\in\cc$ and 
$(x'_1,x'_2)\in\cc$.

We will consider metric spaces equipped with additional structures. 
If $(E,d)$ is a compact metric space, we let $C_0(\R,E)$ be the space
of all continuous functions $\gamma:\R\la E$ such that, for every $\ve>0$,
there exists $T>0$ such that $d(\gamma(t),\gamma(T))<\ve $
and $d(\gamma(-t),\gamma(-T))<\ve $ for every $t\geq T$. By convention, if
$\gamma:[a,b]\la E$ is only (continuous and) defined on an interval $[a,b]$,
we view it as an element of $C_0(\R,E)$ by extending it so that it is constant on $(-\infty,a]$
and on $[b,\infty)$. 
A curve-decorated and pointed (compact) measure 
metric space is then a compact metric space $(E,d)$ equipped with a finite Borel measure $\mu$
(sometimes called the volume measure),
with a curve $\gamma\in C_0(\R,E)$, and with a distinguished point $x$. We write $\M^{GHPU\bullet}$ for the set of all
isometry classes of curve-decorated and pointed compact measure metric spaces (here $(E,d,\mu,\gamma,x)$
and $(E',d',\mu',\gamma',x')$ are isometry equivalent if there exists an isometry $\Phi$ from 
$E$ onto $E'$ such that $\Phi_*\mu=\mu'$, $\gamma'=\Phi\circ\gamma$, and $\Phi(x)=x'$).
One can equip $\M^{GHPU\bullet}$ with the so-called
Gromov-Hausdorff-Prokhorov-uniform distance $\dd_{GHPU\bullet}$, which is defined by
\begin{align*}
&\dd_{GHPU\bullet}((E_1,d_1,\mu_1,\gamma_1,x_1),(E_2,d_2,\mu_2,\gamma_2,x_2))\\
&:=\inf\Big\{d^E_{\mathrm{H}}(\Phi_1(E_1),\Phi_2(E_2)) \vee d^E_\mathrm{P}((\Phi_{1})_*\mu_1 ,(\Phi_2)_*\mu_2)
\vee\sup_{t\in\R}d(\Phi_1\circ\gamma_1(t),\Phi_2\circ\gamma_2(t))\vee d(\Phi_1(x_1),\Phi_2(x_2))\Big\},
\end{align*}
where the infimum is over all isometric embeddings $\Phi_1:E_1\la E$ and $\Phi_2:E_2\la E$ of 
$E_1$ and $E_2$ into the same compact metric space $(E,d)$, and $d^E_{\mathrm{P}}$
denotes the Prokhorov metric on the space of all finite measures on $E$. By a straightforward adaptation
of the arguments of \cite[Section 2.2]{GM0}, one verifies that $\dd_{GHPU\bullet}$ is a complete
separable metric on $\M^{GHPU\bullet}$.

Following \cite{GM0}, we will also use the space  $\M^{GHPU}$ of all isometry classes of (non-pointed) curve-decorated 
compact measure metric spaces, which is equipped with
distance $\dd_{GHPU}$ defined exactly as $\dd_{GHPU\bullet}$ in the last display by just omitting the last term $d(\Phi_1(x_1),\Phi_2(x_2))$.
Then $(\M^{GHPU},\dd_{GHPU})$ is again a Polish space \cite{GM0}.

\begin{proposition}
\label{GHPU-criterion}
Let $(E_n,d_n,\mu_n,\gamma_n,x_n)$, for $n\in\N$, and $(E_\infty,d_\infty,\mu_\infty,\gamma_\infty,x_\infty)$
be elements of $\M^{GHPU\bullet}$. Suppose that $(E_n,d_n,\mu_n,\gamma_n,x_n)$ converges to
$(E_\infty,d_\infty,\mu_\infty,\gamma_\infty,x_\infty)$ in $(\M^{GHPU\bullet},\dd_{GHPU\bullet})$, as $n\to\infty$.
Then, we can find a compact metric space $(E,d)$ and isometric embeddings $\Phi_n:E_n\la E$
and $\Phi_\infty:E_\infty\la E$ such that $\Phi_n(E_n)\la \Phi_\infty(E_\infty)$ for  the Hausdorff metric,
$(\Phi_{n})_*\mu_n\la (\Phi_{\infty})_*\mu_\infty$ for the Prokhorov metric, $\Phi_n\circ\gamma_n(t)\la \Phi_\infty\circ\gamma_\infty(t)$
uniformly in $t$, and $\Phi_n(x_n)\la \Phi_\infty(x_\infty)$, as $n\to\infty$.
\end{proposition}

This is the exact analog of \cite[Proposition 1.5]{GM0}, which deals with $\M^{GHPU}$ instead of $\M^{GHPU\bullet}$.
The proof is the same as in \cite{GM0}.  
In what follows, we will be interested in random metric spaces in $\M^{GHPU\bullet}$, and particularly in the special case where the
distinguished point is chosen ``uniformly'' according to the volume measure. The following lemma 
will be useful.

\begin{lemma}
\label{pointed-notpointed}
Let $(X^n,D^n,\Upsilon^n,\Gamma^n)$, for $n\in\N\cup\{\infty\}$, be random variables with values in $\M^{GHPU}$.
Assume that $(X^n,D^n,\Upsilon^n,\Gamma^n)$ converges to $(X^\infty,D^\infty,\Upsilon^\infty,\Gamma^\infty)$
in distribution when $n\to\infty$. Also assume that $0<\E[\Upsilon^n(X^n)]<\infty$
for every $n\in \N\cup\{\infty\}$, and that
\begin{equation}
\label{tec-point}
\E[\Upsilon^n(X^n)] \build{\la}_{n\to\infty}^{} \E[\Upsilon^\infty(X^\infty)].
\end{equation}
For every $n\in\N\cup\{\infty\}$, define a probability measure $\mathbf{\Theta}_n$ on $\M^{GHPU\bullet}$ by
setting, for every bounded continuous real function $F$ on $\M^{GHPU\bullet}$,
$$\int F\,\dd\mathbf{\Theta}_n = \frac{1}{\E[\Upsilon^n(X^n)] }\E\Big[\int \Upsilon^n(\dd x)\, F((X^n,D^n,\Upsilon^n,\Gamma^n,x))\Big].$$
Then $\mathbf{\Theta}_n$ converges weakly to $\mathbf{\Theta}_\infty$ as $n\to\infty$.
\end{lemma}

\proof By the Skorokhod representation theorem, we may assume that
$(X^n,D^n,\Upsilon^n,\Gamma^n)$ converges almost surely to $(X^\infty,D^\infty,\Upsilon^\infty,\Gamma^\infty)$. We then observe that,
if $F$ is bounded and continuous on $\M^{GHPU\bullet}$, the mapping
$$(E,d,\mu,\gamma)\mapsto \int \mu(\dd x)\,F((E,d,\mu,\gamma,x))$$
is continuous on $\M^{GHPU}$ (we leave the proof to the reader). It follows that we have a.s.
$$\int \Upsilon^n(\dd x)\, F((X^n,D^n,\Upsilon^n,\Gamma^n,x)) \build{\la}_{n\to\infty}^{} 
\int \Upsilon^\infty(\dd x)\, F((X^\infty,D^\infty,\Upsilon^\infty,\Gamma^\infty,x)).$$
Using dominated convergence and our assumption \eqref{tec-point}, we obtain that 
$\int F\,\dd\mathbf{\Theta}_n\la\int F\,\dd\mathbf{\Theta}_\infty$. \endproof

\subsection{Snake trajectories}
\label{sna-tra}

To construct the models of random geometry that we consider, we will use the formalism 
of snake trajectories. 
A (one-dimensional) finite path $\w$ is a continuous mapping $\w:[0,\zeta]\la\R$, where the
number $\zeta=\zeta_{(\w)}\geq 0$ is called the lifetime of $\w$. We let 
$\mathfrak{W}$ denote the space of all finite paths, which is a Polish space when equipped with the
distance
$$d_\mathfrak{W}(\w,\w'):=|\zeta_{(\w)}-\zeta_{(\w')}|+\sup_{t\geq 0}|\w(t\wedge
\zeta_{(\w)})-\w'(t\wedge\zeta_{(\w')})|.$$
The endpoint or tip of the path $\w$ is denoted by $\wh \w=\w(\zeta_{(\w)})$.
For $x\in\R$, we
set $\mathfrak{W}_x:=\{\w\in\mathfrak{W}:\w(0)=x\}$. The trivial element of $\mathfrak{W}_x$ 
with zero lifetime is identified with the point $x$ of $\R$. 

\smallskip
\noindent {\it Definition}. Let $x\in\R$. 
A snake trajectory with initial point $x$ is a continuous mapping $s\mapsto \omega_s$
from $\R_+$ into $\mathfrak{W}_x$ 
that satisfies the following two properties:
\begin{enumerate}
\item[\rm(i)] We have $\omega_0=x$ and the number $\sigma(\omega):=\sup\{s\geq 0: \omega_s\not =x\}$,
called the duration of the snake trajectory $\omega$,
is finite (by convention $\sigma(\omega)=0$ if $\omega_s=x$ for every $s\geq 0$). 
\item[\rm(ii)] {\rm (Snake property)} For every $0\leq s\leq s'$, we have
$\omega_s(t)=\omega_{s'}(t)$ for every $t\in[0,\displaystyle{\min_{s\leq r\leq s'}} \zeta_{(\omega_r)}]$.
\end{enumerate} 

\smallskip
We write $\S_x$ for the set of all snake trajectories with initial point $x$
and $\S=\bigcup_{x\in\R}\S_x$ for the set of all snake trajectories. If $\omega\in \S$, we often write $W_s(\omega):=\omega_s$ and $\zeta_s(\omega):=\zeta_{(\omega_s)}$
for every $s\geq 0$. The set $\S$ is a Polish space for the distance
$$d_{\S}(\omega,\omega'):= |\sigma(\omega)-\sigma(\omega')|+ \sup_{s\geq 0} \,d_\mathfrak{W}(W_s(\omega),W_{s}(\omega')).$$
We stress that a snake trajectory $\omega$ is completely determined 
by the knowledge of the lifetime function $s\mapsto \zeta_s(\omega)$ and of the tip function $s\mapsto \wh W_s(\omega)$: See \cite[Proposition 8]{ALG}.

Let $\omega\in \S$ be a snake trajectory and $\sigma=\sigma(\omega)$. The lifetime function $s\mapsto \zeta_s(\omega)$ codes a
compact $\R$-tree, which will be denoted 
by $\t_{(\omega)}$ and called the {\it genealogical tree} of the snake trajectory. This $\R$-tree is the quotient space $\t_{(\omega)}:=[0,\sigma]/\!\sim$ 
of the interval $[0,\sigma]$
for the equivalence relation
$$s\sim s'\ \hbox{if and only if }\ \zeta_s(\omega)=\zeta_{s'}(\omega)= \min_{s\wedge s'\leq r\leq s\vee s'} \zeta_r(\omega),$$
and $\t_{(\omega)}$ is equipped with the distance induced by
$$d_{(\omega)}(s,s'):= \zeta_s(\omega)+\zeta_{s'}(\omega)-2 \min_{s\wedge s'\leq r\leq s\vee s'} \zeta_r(\omega).$$
(notice that $d_{(\omega)}(s,s')=0$ if and only if $s\sim s'$).
We write $p_{(\omega)}:[0,\sigma]\la \t_{(\omega)}$
for the canonical projection, and the mapping $[0,\sigma]\ni t\mapsto p_{(\omega)}(t)$ can be viewed as a
cyclic exploration of $\t_{(\omega)}$. By convention, $\t_{(\omega)}$ is rooted at the point
$\rho_{(\omega)}:=p_{(\omega)}(0)$, and the volume measure on $\t_{(\omega)}$ is defined as the pushforward of
Lebesgue measure on $[0,\sigma]$ under $p_{(\omega)}$. If $u,v\in\t_{(\omega)}$, $\llbracket u,v\rrbracket$ denotes 
the geodesic segment between $u$ and $v$ in $\t_{(\omega)}$.
 The segment $\llbracket \rho_{(\omega)},u\rrbracket$
 is called the ancestral line of $u$.

By property (ii) in the definition of  a snake trajectory, the condition $p_{(\omega)}(s)=p_{(\omega)}(s')$ implies that 
$W_s(\omega)=W_{s'}(\omega)$. So the mapping $s\mapsto W_s(\omega)$ can be viewed as defined on the quotient space $\t_{(\omega)}$.
For $u\in\t_{(\omega)}$, we set $\ell_u(\omega):=\wh W_s(\omega)$ whenever  $s\in[0,\sigma]$ is such that $u=p_{(\omega)}(s)$  (by the previous observation, this does not
depend on the choice of $s$). We interpret $\ell_u(\omega)$ as a ``label'' assigned to the ``vertex'' $u$ of $\t_{(\omega)}$. 
Notice that the mapping $u\mapsto \ell_u(\omega)$ is continuous on $\t_{(\omega)}$, and that, for every $s\geq 0$, the path
$W_s(\omega)$ records the labels $\ell_u(\omega)$ along the ancestral line $\llbracket \rho_{(\omega)},p_{(\omega)}(s)\rrbracket$. 
We will use the notation $W_*(\omega):=\min\{\ell_u(\omega):u\in\t_{(\omega)}\}$.

We now introduce an important operation on snake trajectories in $\S$. 
Let $x,y\in \R$ with $y<x$. For every $\w\in\mathfrak{W}_x$, set
$$\tau_y(\w):=\inf\{t\in[0,\zeta_{(\w)}]: \w(t)=y\}$$
with the usual convention $\inf\varnothing =\infty$ (this convention will be in force throughout this work
unless otherwise indicated). Then, if 
$\omega\in \S_x$, we set, for every $s\geq 0$,
$$\eta_s(\omega):=\inf\Big\{t\geq 0:\int_0^t \mathrm{d}r\,\mathbf{1}_{\{\zeta_{(\omega_r)}\leq\tau_y(\omega_r)\}}>s\Big\}.$$
Note that the condition $\zeta_{(\omega_r)}\leq\tau_y(\omega_r)$ holds if and only if $\tau_y(\omega_r)=\infty$ or $\tau_y(\omega_r)=\zeta_{(\omega_r)}$.
Then, setting $\omega'_s=\omega_{\eta_s(\omega)}$ for every $s\geq 0$ defines an element $\omega'$ of $\S_x$,
which will be denoted by  $\tr_y(\omega)$ and called the truncation of $\omega$ at $y$
(see \cite[Proposition 10]{ALG}). The effect of the time 
change $\eta_s(\omega)$ is to ``eliminate'' those paths $\omega_s$ that hit $y$ and then survive for a positive
amount of time. The genealogical tree $\t_{(\tr_y(\omega))}$ is
canonically and isometrically identified to the closed set $$\big\{v\in\t_{(\omega)}:\ell_u(\omega)>y\hbox{ for every }u\in\llbracket \rho_{(\omega)},v\rrbracket\backslash\{v\}\big\},$$
and this identification preserves labels. In what follows, we will therefore view  $\t_{(\tr_y(\omega))}$
as a subset of $\t_{(\omega)}$. Informally, $\t_{(\tr_y(\omega))}$ is obtained from 
$\t_{(\omega)}$ by pruning branches at the level where labels first take the value $y$. 

We can then also define the excursions of $\omega$ away from a given level. Consider $\omega\in \S_x$
and  $y<x$. Let $(\alpha_j,\beta_j)$, $j\in J$, be the connected components of the open set
$$\{s\in[0,\sigma]:\tau_y(\omega_s)<\zeta_{(\omega_s)}\},$$
and notice that we have $\omega_{\alpha_j}=\omega_{\beta_j}$, for every $j\in J$, by the snake property.
For every $j\in J$, we define a snake trajectory $\omega^j\in\S_y$ by setting
$$\omega^j_{s}(t):=\omega_{(\alpha_j+s)\wedge\beta_j}(\zeta_{(\omega_{\alpha_j})}+t)\;,\hbox{ for }0\leq t\leq \zeta_{(\omega^j_s)}
:=\zeta_{(\omega_{(\alpha_j+s)\wedge\beta_j})}-\zeta_{(\omega_{\alpha_j})}\hbox{ and } s\geq 0.$$
We say that $\omega^j$, $j\in J$, are the excursions of $\omega$ away from $y$. We note that, for every $j\in J$,
the tree $\t_{(\omega^j)}$ is canonically identified to a subtree of $\t_{(\omega)}$ consisting of
descendants of $p_{(\omega)}(\alpha_j)=p_{(\omega)}(\beta_j)$. 

\subsection{The Brownian snake excursion 
measure on snake trajectories}
\label{sna-mea}

Let $x\in\R$. The Brownian snake excursion 
measure $\N_x$ is the $\sigma$-finite measure on $\S_x$ that satisfies the following two properties: Under $\N_x$,
\begin{enumerate}
\item[(i)] the distribution of the lifetime function $(\zeta_s)_{s\geq 0}$ is the It\^o 
measure of positive excursions of linear Brownian motion, normalized so that, for every $\ve>0$,
$$\N_x\Big(\sup_{s\geq 0} \zeta_s >\ve\Big)=\frac{1}{2\ve};$$
\item[(ii)] conditionally on $(\zeta_s)_{s\geq 0}$, the tip function $(\wh W_s)_{s\geq 0}$ is
a Gaussian process with mean $x$ and covariance function 
$$K(s,s'):= \min_{s\wedge s'\leq r\leq s\vee s'} \zeta_r.$$
\end{enumerate}
Informally, the lifetime process $(\zeta_s)_{s\geq 0}$ evolves under $\N_x$ like a Brownian excursion,
and conditionally on $(\zeta_s)_{s\geq 0}$, each path $W_s$ is a linear Brownian path started from $x$ with lifetime $\zeta_s$, which
is ``erased'' from its tip when $\zeta_s$ decreases and is ``extended'' when $\zeta_s$ increases.
The measure $\N_x$ can be interpreted as the excursion measure away from $x$ for the 
Markov process in $\mathfrak{W}_x$ called the Brownian snake.
We refer to 
\cite{Zurich} for a detailed study of the Brownian snake. 
For every $y<x$, we have
\begin{equation}
\label{hittingpro}
\N_x(W_*\leq y)={\displaystyle \frac{3}{2(x-y)^2}},
\end{equation}
where we recall the notation $W_*(\omega)$ for the minimal label on $\t_{(\omega)}$. See e.g. \cite[Section VI.1]{Zurich} for a proof. 

\medskip
\noindent{\it Exit measures.} Let $x,y\in\R$, with $y<x$. Under the measure $\N_x$, one can make sense of a quantity that measures ``how many'' paths 
$W_s$  hit $y$. One shows \cite[Proposition 34]{Disks} that the limit
\begin{equation}
\label{formu-exit}
L^y_t:=\lim_{\ve \da 0} \frac{1}{\ve^2} \int_0^t \dd s\,\mathbf{1}_{\{\tau_y(W_s)=\infty,\, \wh W_s<y+\ve\}}
\end{equation}
exists uniformly in $t\geq 0$, $\N_x$ a.e., and defines a continuous nondecreasing function, which is 
obviously constant on $[\sigma,\infty)$. 
The process $(L^y_t)_{t\geq 0}$ is called the exit local time from $(y,\infty)$, and the exit measure 
$\z_y$ is defined by $\z_y:=L^y_\infty=L^y_\sigma$. Then, $\N_x$ a.e., the topological support of the measure 
$\dd L^y_t$ is exactly the set $\{s\in[0,\sigma]:\tau_y(W_s)=\zeta_s\}$, and, in particular, $\z_y>0$ if and only if one of the paths $W_s$ hits $y$. The definition of $\z_y$
is a special case of the theory of exit measures (see \cite[Chapter V]{Zurich} for this general theory). We will use the formula for
the Laplace transform of $\z_y$: For $\lambda> 0$,
\begin{equation}
\label{Lap-exit}
\N_x\Big(1-\exp (-\lambda \z_y)\Big)= \Big((x-y)\sqrt{2/3} + \lambda^{-1/2}\Big)^{-2}.
\end{equation}
See formula (6) in \cite{CLG} for a brief justification. 

It is useful to observe that $\z_y$ can be defined in terms of the truncated snake $\tr_y(\omega)$. 
To this end, recall the time change $(\eta_s(\omega))_{s\geq 0}$ used to define 
$\tr_y(\omega)$ at the end of Section \ref{sna-tra}, and set $\wt L^y_t=L^y_{\eta_t}$
for every $t\geq 0$. Then $\wt L^y_\infty= L^y_\infty=\z_y$, whereas formula \eqref{formu-exit}
implies that
\begin{equation}
\label{formu-exit-bis}
\wt L^y_t=\lim_{\ve \da 0} \frac{1}{\ve^2} \int_0^t \dd s\,\mathbf{1}_{\{\wh W_s(\tr_y(\omega))<y+\ve\}}
\end{equation}
uniformly for $t\geq 0$, $\N_x$ a.e.

\medskip
\noindent{\it The special Markov property.} We use the notation introduced in  Section \ref{sna-tra}. More precisely, we write  $\omega^j$, $j\in J$, for  the excursions of $\omega$ below $y$ and $(\alpha_j,\beta_j)$, $j\in J$, for the associated time intervals. The special Markov property states that, conditionally on the truncation $\tr_y(\omega)$,  the point measure:
\begin{equation}\label{measure:special:mark}
\sum \limits_{j\in J} \delta_{(L^y_{\alpha_j}, \omega^{j})}
\end{equation}
is Poisson with intensity $\mathbf{1}_{[0,\z_y ]}(t)\, \dd t \,\mathbb{N}_y(\dd \omega)$. We refer to the Appendix of \cite{subor} for a proof. By combining the special 
Markov property with the fact that the ``law'' of $W_*$ under $\N_x$ has no atoms, one easily gets 
that, for every fixed $z\in(-\infty,x)$, the value $z$ is $\N_x$ a.e. not a local minimum 
of the function $s\mapsto \wh W_s$.

\subsection{A technical lemma}

In this section, we establish a lemma that will be useful in forthcoming proofs. This lemma is a direct consequence 
of arguments used in the proof of \cite[Proposition 31]{Disks}, which was the key
result needed for the extension of the distance to the boundary in the construction of 
the Brownian disk presented in \cite{Disks}. We use the notation $\N^{[0]}_r:=\N_r(\cdot\mid W_*\leq 0)$ for every $r>0$.
Under $\N^{[0]}_r(\dd\omega)$, we write $\wt\omega=\mathrm{tr}_0(\omega)$ to simplify notation.
Recall
that $\ell_a(\wt\omega)=\wh W_s(\wt\omega)$ if $a=p_{(\tilde\omega)}(s)$, and
note that $\ell_a(\wt\omega)\geq 0$ for every $a\in\t_{(\tilde\omega)}$, $\N^{[0]}_r(\dd\omega)$ a.s. 
The ``boundary'' of the tree $\t_{(\tilde\omega)}$ is then defined as the set $\partial \t_{(\tilde\omega)}:=\{a\in \t_{(\tilde\omega)}:\ell_a(\tilde \omega)=0\}$.
We set,
for every $s,t\in[0,\sigma(\wt\omega)]$,
$$\Delta^\circ_{(\tilde\omega)}(s,t):= \wh W_s(\wt \omega) + \wh W_t(\wt\omega) - 2 \min_{s\wedge t\leq r\leq s\vee t}\wh W_r(\wt \omega)$$
if the minimum in the last display is positive, and $\Delta^\circ_{(\tilde\omega)}(s,t):=\infty$ otherwise. We then set, for every $a,b\in \t_{(\tilde\omega)}\backslash \partial \t_{(\tilde\omega)}$,
$$\Delta^\circ_{(\tilde\omega)}(a,b):=\min\{\Delta^\circ_{(\tilde\omega)}(s,t):
s,t\in[0,\sigma(\wt\omega)],\,p_{(\tilde\omega)}(s)=a,p_{(\tilde\omega)}(t)=b\},$$
and $$\Delta_{(\tilde\omega)}(a,b):=\inf\Big\{ \sum_{i=1}^p \Delta^\circ_{(\tilde\omega)}(a_{i-1},a_i)\Big\}$$
where the infimum is over all choices of the integer $p\geq 1$ and of $a_0=a,a_1,\ldots,a_{p-1},a_p=b$ in $\t_{(\tilde\omega)}\backslash \partial \t_{(\tilde\omega)}$.
It is not hard to verify that the mapping $(a,b)\mapsto \Delta_{(\tilde\omega)}(a,b)$ takes finite values and is continuous on $(\t_{(\tilde\omega)}\backslash \partial \t_{(\tilde\omega)})
\times( \t_{(\tilde\omega)}\backslash \partial \t_{(\tilde\omega)})$. See the comments following Proposition 30 in \cite{Disks}.

\begin{lemma}
\label{tech-lem}
$\N^{[0]}_r$ a.s., the mapping $(a,b)\mapsto \Delta_{(\tilde\omega)}(a,b)$ has a unique continuous extension to $\t_{(\tilde\omega)}\times \t_{(\tilde\omega)}$. Moreover,
there exists a finite constant $C$, which does not depend on $r$,  such that
$$\N^{[0]}_r\Big( \sup_{a,b\in \t_{(\tilde\omega)}} \Delta_{(\tilde\omega)}(a,b)\Big)=C\,r.$$
\end{lemma}

\proof By scaling, it is enough to consider the case $r=1$. Let us start by proving the first assertion. Since $\Delta_{(\tilde\omega)}$ satisfies the triangle inequality, it is enough to
verify that, for any $a\in \partial\t_{(\tilde\omega)}$, if $(a_n)_{n\in\N}$
is a sequence in $\t_{(\tilde\omega)}\backslash\partial\t_{(\tilde\omega)}$ that 
converges to $a$, we have $\Delta_{(\tilde\omega)}(a_n,a_m)\la 0$ as $n,m\to\infty$. 
To get this, write $m_c$ for the minimal label along the ancestral line of $c$ in $\t_{(\tilde\omega)}$, 
for every $c\in\t_{(\tilde\omega)}$, and, for every $\delta >0$, let $\mathcal{C}^{\delta}_{j}$, $j\in\{1,\ldots,N_\delta\}$, be those connected components 
of the open set $\{c\in \t_{(\tilde\omega)}\backslash\partial\t_{(\tilde\omega)}:m_c<\delta\}$
whose closure intersects the ``boundary'' $\partial\t_{(\tilde\omega)}$. By formula (53) in \cite{Disks}, we have
$$\sup_{1\leq j\leq N_\delta} \Big(\sup_{b,b'\in \mathcal{C}^{\delta}_{j}} \Delta_{(\tilde\omega)}(b,b')\Big) \build{\la}_{\delta\to 0}^{} 0$$
(note that formula (53) in \cite{Disks} deals with a function $\Delta(x,y)$ which is defined in a slightly different way
than $\Delta_{(\tilde\omega)}(x,y)$, but the arguments apply as well to $\Delta_{(\tilde\omega)}(x,y)$). 
For every fixed $\delta>0$, there is a unique index $j$ such that $a$ belongs to the closure of $\mathcal{C}^{\delta}_{j}$, and, for $n$ large enough,
$a_n$ must belong to $\mathcal{C}^{\delta}_{j}$. The desired convergence of $\Delta_{(\tilde\omega)}(a_n,a_m)$ to $0$ then follows from the
last display. 

Let us turn to the second assertion. To simplify notation, we write $\ell_a$ instead of $\ell_a(\wt\omega)$
and $\rho$ instead of $\rho_{(\tilde\omega)}$. Our goal is to verify that
$$\N^{[0]}_1\Big( \sup_{a\in \t_{(\tilde\omega)}} \Delta_{(\tilde\omega)}(\rho,a)\Big)<\infty,$$
which will immediately give the second assertion (for $r=1$) since $\Delta_{(\tilde\omega)}$ satisfies the triangle inequality.
We need to recall some ingredients of the proof of Proposition 31 in \cite{Disks}. 
We first introduce the reduced tree of $\t_{(\tilde\omega)}$, which consists
 of all points $a$ of $\t_{(\tilde\omega)}\backslash \partial \t_{(\tilde\omega)}$ that have at least one descendant with label $0$ ($a$ belongs to the reduced tree 
 if there exists $b\in\t_{(\tilde\omega)}$ such that $\ell_b=0$ and $a\in\llbracket \rho,b\rrbracket$). Let $\t^\triangledown$
 stand for this subtree. Then the tree $\t^\triangledown$ is a binary $\R$-tree, which can be constructed by induction 
 as follows. One starts 
from a line segment connecting the root $\rho$ to a first branching point 
 $a_\varnothing$. To this branching point are attached two other line segments connecting
 $a_\varnothing$ to branching points $a_1$ and $a_2$, listed in the
order prescribed by the exploration $t\mapsto p_{(\tilde\omega)}(t)$ of $\t_{(\tilde\omega)}$. To $a_1$ (respectively to $a_2$) are then attached two line segments
 connecting $a_1$ (resp. $a_2$) to branching points $a_{(1,1)}$ and $a_{(1,2)}$ (resp.
 $a_{(2,1)}$ and $a_{(2,2)}$) and so on. The reason for introducing this reduced tree is the bound
 \begin{equation}
 \label{tec-bd10}
 \sup_{a\in \t_{(\tilde\omega)}} \Delta_{(\tilde\omega)}(\rho,a)\leq 2\sup_{(i_1,i_2,\ldots)\in \{1,2\}^\N}\Bigg(\sum_{n=0}^\infty \ell_{a_{(i_1,\ldots,i_n)}}\Bigg)
 + 4\, \sup_{a\in\t_{(\tilde\omega)}} \ell_a,
 \end{equation}
 which easily follows from the fact that $\Delta_{(\tilde\omega)}(a_{(i_1,\ldots,i_{n-1})}, a_{(i_1,\ldots,i_{n})})\leq \ell_{a_{(i_1,\ldots,i_{n-1})}}+\ell_{a_{(i_1,\ldots,i_n)}}$
 (see the end of the proof of \cite[Proposition 31]{Disks} for more details). The second term in the right-hand side of \eqref{tec-bd10} has finite expectation
 under $\N^{[0]}_1$ because, for every $x>1$,
 $$\N^{[0]}_1\Big( \sup_{a\in \t_{(\tilde\omega)}}\ell_a >x\Big)\leq \frac{2}{3}\,\N_1\Big( \sup_{a\in \t_{(\tilde\omega)}}\ell_a >x\Big) \leq (x-1)^{-2},$$
 using \eqref{hittingpro}. So it remains to verify that the first  term in the right-hand side of \eqref{tec-bd10} also has finite expectation
 under $\N^{[0]}_1$. To this end, we rely on the formula
 $$\N^{[0]}_1\Big((\ell_{a_{(i_1,\ldots,i_n)}})^{5/2}\Big)= \Big(\frac{24}{49}\Big)^{n+1},$$
 which is obtained in the proof of \cite[Proposition 31]{Disks} as a consequence of the recursive structure of the tree $\t^\triangledown$.
 We fix $\alpha\in(0,1)$ such that $2\alpha^{-5/2}<49/24$. Then, for every $x>1$,
 \begin{align*}
 \N^{[0]}_1\Bigg(\sup_{(i_1,i_2,\ldots)\in \{1,2\}^\N}\Bigg(\sum_{n=0}^\infty \ell_{a_{(i_1,\ldots,i_n)}}\Bigg) >x\Bigg)
 &\leq \sum_{n=0}^\infty  \N^{[0]}_1\Bigg(\Big(\sup_{(i_1,i_2,\ldots,i_n)\in \{1,2\}^n} \ell_{a_{(i_1,\ldots,i_n)}}\Big) >(1-\alpha)\alpha^nx\Bigg)\\
 &\leq \sum_{n=0}^\infty 2^n\,\sup_{(i_1,i_2,\ldots,i_n)\in \{1,2\}^n}\N^{[0]}_1\Big(\ell_{a_{(i_1,\ldots,i_n)}}>(1-\alpha)\alpha^nx\Big)\\
 &\leq \sum_{n=0}^\infty 2^n\times ((1-\alpha)\alpha^nx)^{-5/2}\times \Big(\frac{24}{49}\Big)^{n+1}\\
 &=c\,x^{-5/2}
 \end{align*}
with some constant $c<\infty$. This completes the proof. \endproof

\section{Peeling of a triangulation with a boundary}
\label{discrete}

Our goal
in this section is to discuss certain properties of Boltzmann distributed pointed 
triangulations with a simple boundary. More precisely, we are interested in the discrete hull 
 whose radius is the distance 
from the distinguished vertex to the boundary. Thanks to 
the results of \cite{AHS}, this study will allow us
to derive similar properties for the free pointed Brownian disk. In our investigation of 
Boltzmann distributed triangulations, it will
be convenient to use the peeling algorithm. 

\subsection{Peeling probabilities}
\label{peeling-proba}

For integers $L\geq 1$ and $k\geq 0$, we let
$\T^1(L,k)$ be the set of all  rooted planar triangulations of type I 
(i.e. loops and multiple edges are allowed) with a simple boundary
of length $L$ and $k$ inner vertices. By convention, triangulations with a simple boundary
are rooted on the boundary in such a way that the external
face (of degree $L$) lies to the left of the root edge (see e.g.~the introduction of \cite{AHS} 
for a more detailed presentation of triangulations with a boundary). Then (see e.g. Theorem~1.1 in \cite{BF}), 
$\T^1(1,0)=\varnothing$ and, for $(L,k)\not = (1,0)$,
\begin{equation}
\label{triangu1}
\# \T^1(L,k)= 4^{k-1}\,\frac{(2L+3k-5)!!}{k!\,(2L+k-1)!!}\,L{2L\choose L} \build{\sim}_{k\to\infty}^{} C^{(1)}(L) (12\sqrt{3})^k\,k^{-5/2},
\end{equation}
where
\begin{equation}
\label{triangu11}
C^{(1)}(L):=\frac{3^{L-2}}{4\sqrt{2\pi}}\,L{2L\choose L} \build{\sim}_{L\to\infty}^{} \frac{1}{36\pi\sqrt{2}}\sqrt{L}\,12^L.
\end{equation}
(When $L=2$ and $k=0$, formula \eqref{triangu1} is  valid with the convention $(-1)!!=1$, provided we consider the ``trivial triangulation'' as in \cite{CLG}.)
Assuming that $L\geq 2$, we have
\begin{equation}
\label{triangu12} 
Z(L):=\sum_{k=0}^\infty (12\sqrt{3})^{-k}\,\#\T^1(L,k)=
\frac{6^L(2L-5)!!}{8\sqrt{3}\,L!}.
\end{equation}
(see e.g.~\cite[Section 2.2]{AC}). We set $\T^1(L):=\bigcup_{k\geq 0}\T^1(L,k)$. A random triangulation $\tau$ in $\T^1(L)$ is said to be Boltzmann distributed if
$\P(\tau=\theta)= Z(L)^{-1}(12\sqrt{3})^{-k}$, for every $k\geq 0$ and $\theta\in\T^1(L,k)$. We will also consider rooted and pointed
triangulations with a boundary, which in addition to the root edge have a distinguished vertex, which
can be any inner vertex of the triangulation. We can then define
Boltzmann distributed rooted and pointed planar
triangulations in exactly the same way as we did in the non-pointed case (using the fact that
$\sum_{k=0}^\infty k\, (12\sqrt{3})^{-k}\# \T^1(L,k)<\infty$, by \eqref{triangu1}).

For integers $L\geq 1$, $p\geq 1$, and $k\geq 0$, let $\T^2(L,p,k)$ be the set of all planar triangulations with two simple boundaries of respective
lengths $L$ and $p$, and $k$ inner vertices, that are rooted on both boundaries (with the same 
convention for the orientation of the root edges). Notice that we distinguish the first and the
second boundary, and that
the size of the first one is  $L$. We refer to the introduction of \cite{BF} for a precise definition of
triangulations with several boundaries, and note in particular that the boundaries are assumed to be  vertex disjoint.
According to \cite{Kri} (see also \cite{BF}),
\begin{equation}
\label{triangu2}
\#\T^2(L,p,k)=\frac{4^k\,(2(L+p)+3k-2)!!}{k!\,(2(L+p)+k)!!} \,L {2L\choose L}\,p{2p\choose p}.
\end{equation}
Set
$$Z'(L,p):=\sum_{k=0}^\infty (12\sqrt{3})^{-k}\,\#\T^2(L,p,k).$$
Using calculations in Krikun \cite{Kri}, one checks that
\begin{equation}
\label{triangu21}
Z'(L,p)= \frac{1}{2}\,\frac{3^{L+p}}{L+p}\,\,L {2L\choose L}\,p{2p\choose p}=\frac{6^4\, \pi}{L+p} C^{(1)}(L)\times C^{(1)}(p).
\end{equation}
In the appendix below, we explain how formula \eqref{triangu21} can be deduced from \cite{Kri}.

 From now on, we always assume that $L\geq 2$. Let $\T^2(L,p)$ stand for the union of all $\T^2(L,p,k)$ for $k\geq 0$. Consider a random triangulation
$\tau$ of $\T^2(L,p)$ distributed according to Boltzmann weights. This means that, if $\theta$ is a
given triangulation of $\T^2(L,p,k)$ for some $k\geq 0$,
$$\P(\tau=\theta)= Z'(L,p)^{-1}\,(12\sqrt{3})^{-k}.$$

Consider a given edge of the second boundary of $\tau$. This edge, which will be called the revealed edge, is chosen in a deterministic manner
given the root of the second boundary. Let 
$\Delta$ be the triangle incident to this edge, which is called the revealed triangle.
Several configurations may occur (see Fig.~1 for an illustration).

\smallskip
1. The third vertex of $\Delta$ does not lie on any of the two boundaries. Then, if we ``remove'' $\Delta$ from $\tau$, we get a
triangulation of $\T^2(L,p+1,k)$ for some $k\geq 0$ --- the root edge on the second boundary can be chosen 
again in a deterministic manner from the position of the second root edge in $\tau$. We observe that configuration 1 occurs 
if and only if $\tau$ is obtained by filling the space between the first boundary and the (new) second boundary
by a triangulation of $\T^2(L,p+1,k)$ for some $k\geq 0$. For any fixed choice of the
latter triangulation, the probability of the corresponding event is
$$Z'(L,p)^{-1}\,(12\sqrt{3})^{-k-1}.$$
Finally, the probability of configuration 1 is
$$Z'(L,p)^{-1}\,\sum_{k=0}^\infty (12\sqrt{3})^{-k-1}\,\#\T^2(L,p+1,k)=\frac{1}{12\sqrt{3}}\,\frac{Z'(L,p+1)}{Z'(L,p)},$$
and this quantity is also equal to
\begin{equation}
\label{peel1}
\frac{1}{12\sqrt3} \, \frac{L+p}{L+p+1}\,\frac{C^{(1)}(p+1)}{C^{(1)}(p)}.
\end{equation}

\smallskip
2. The third vertex of $\Delta$ belongs to the second boundary, and the revealed triangle
disconnects the first boundary from $m$ edges of the second boundary, where $m\in\{0,1,\ldots,p-1\}$,
and these edges may lie either to the right or to the left of the revealed edge. Consider the case
where these edges lie to the right of the revealed edge (the other case is symmetric).

The complement of the
revealed triangle in the initial triangulation has two connected components (when $m=1$, one of them
may be the trivial triangulation). The one incident to the first boundary must be filled by 
a triangulation of $\T^2(L,p-m,k)$ for some $k\geq 0$, and the other one 
is filled by a triangulation of $\T^1(m+1,j)$ for some $j\geq 0$ ($j\geq 1$ if $m=0$). If these two triangulations
are fixed, the probability of the resulting event is
$$Z'(L,p)^{-1}\,(12\sqrt{3})^{-k-j}.$$
Hence, the probability of the configuration is
$$Z'(L,p)^{-1}\,\sum_{k,j=0}^\infty (12\sqrt{3})^{-k-j}\,\#\T^2(L,p-m,k)\,\#\T^1(m+1,j)=Z(m+1)\,\frac{Z'(L,p-m)}{Z'(L,p)}.$$
The last quantity is also equal to
\begin{equation}
\label{peel2}
\frac{L+p}{L+p-m}\,Z(m+1)\, \frac{C^{(1)}(p-m)}{C^{(1)}(p)}.
\end{equation}

3. The third vertex of the revealed triangle $\Delta$ belongs to the first boundary. To evaluate the probability of this event,
we first notice that there are $L$ possible choices for the third vertex. Then, given the revealed triangle,
the initial triangulation $\tau$ is determined from a triangulation of $\T^1(L+p+1,k)$ for some $k\geq 0$, and if this
triangulation is given, the probability is
$$Z'(L,p)^{-1}\,(12\sqrt{3})^{-k}.$$
The probability of configuration 3 is thus 
$$L\times Z'(L,p)^{-1}\,\sum_{k=0}^\infty (12\sqrt{3})^{-k}\,\#\T^1(L+p+1,k)= L\,\frac{Z(L+p+1)}{Z'(L,p)}.$$
When $L$ is large, we have
$$Z(L+p+1) \sim \frac{\sqrt{3}}{8\sqrt{\pi}}\,12^{L+p}\,(L+p)^{-5/2},$$
uniformly in $p$ (cf. Section 6.1 in \cite{CLG}), whereas 
$$Z'(L,p)\sim \frac{3^{L+p}\times\,4^{L}}{2(L+p)} \sqrt{\frac{L}{\pi }}\,p{2p\choose p}.$$
It follows that the probability of configuration 3 behaves, when $L$ and $p$ are large, like
$$\frac{\sqrt{3\pi}}{4}\,\sqrt{\frac{L}{p}}\,(L+p)^{-3/2}.$$
 
 \begin{figure}[!h]
 \label{Fig1}
 \begin{center}
    \includegraphics[height=5.2cm,width=5.2cm]{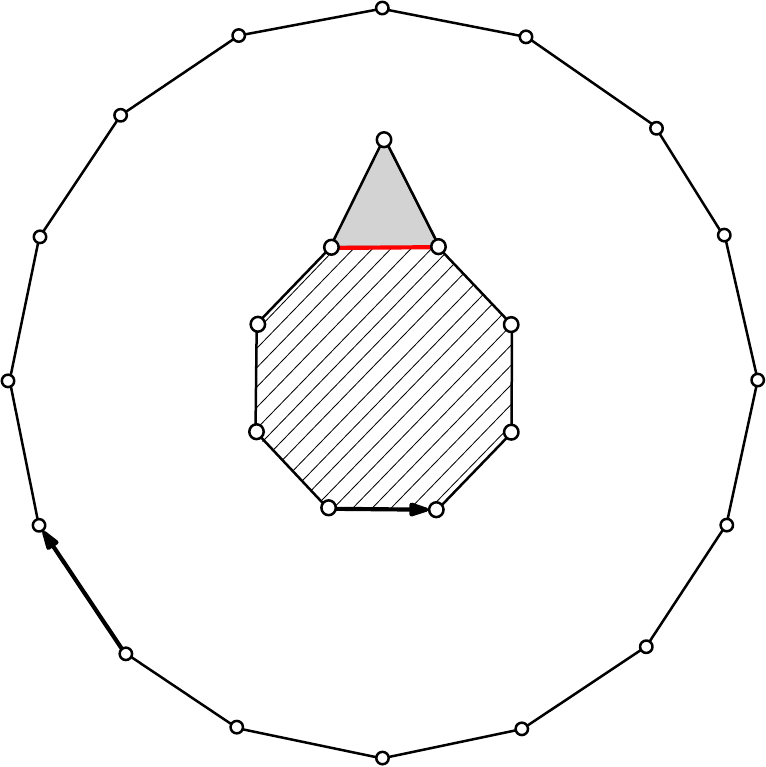} \hspace{0.5cm}\includegraphics[height=5.2cm,width=5.2cm]{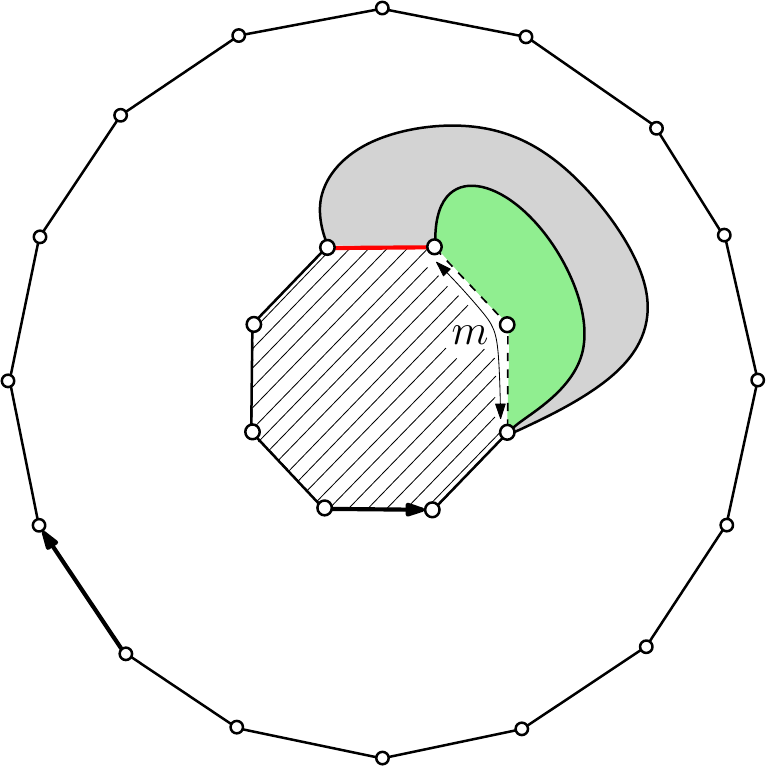} \hspace{0.5cm}\includegraphics[height=5.2cm,width=5.2cm]{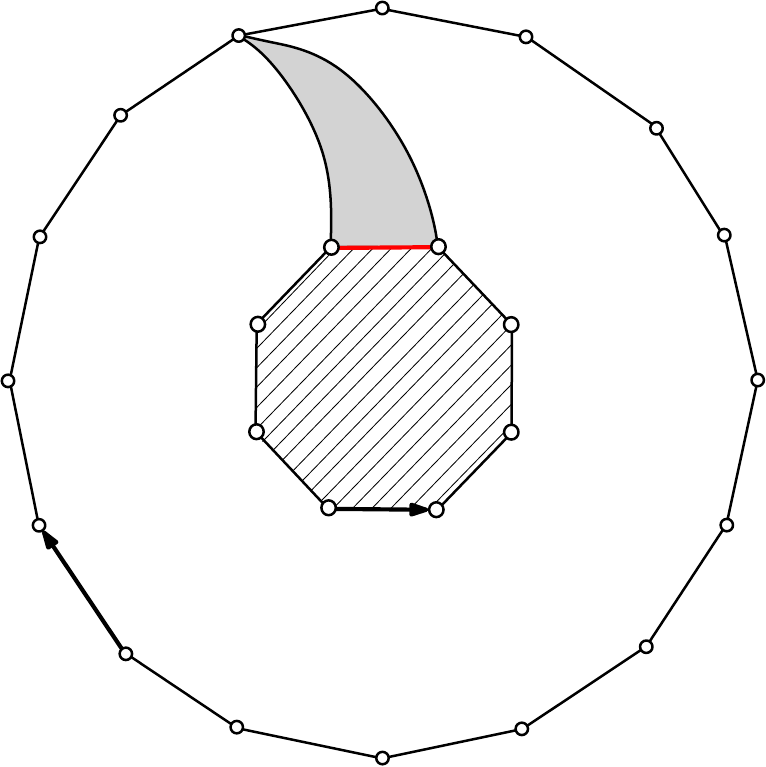}
 \end{center}
    \caption{Illustration from left to right of configurations 1,2 and 3. }
  \end{figure}

\subsection{The peeling process}
\label{peeling-pro}

In this section, the integer $L\geq 2$ is fixed, and we also fix an integer $p_0\geq 1$. As previously, we
consider a Boltzmann triangulation $\tau^L$ in $\T^2(L,p_0)$. We define 
a peeling algorithm giving rise to a sequence $(\tau^L_n)_{0\leq n<\zeta_L}$ of 
triangulations with two boundaries, where $\zeta_L\geq 1$ is a random integer.
Precisely, we take $\tau^L_0=\tau^L$, and then we proceed inductively as follows.
\begin{itemize}
\item[] At step $n$, assuming that $n<\zeta_L$, we choose an edge of the second boundary of $\tau^L_n$ and reveal the
triangle incident to this edge in the way explained in the previous section. If configuration 1 occurs, we let $\tau^L_{n+1}$ be obtained by
removing the revealed triangle in $\tau^L_n$. If configuration 2 occurs,
$\tau^L_{n+1}$ is obtained by removing both the revealed triangle and those triangles 
that are disconnected by the revealed triangle from the first boundary. Finally, if configuration 3 occurs, we
take $\zeta_L:=n+1$. 
\end{itemize}
The preceding description is a little informal since we need to specify how the revealed edge is
chosen at each step. To this end, define, for every $n<\zeta_L$, the planar map $\check\tau^L_n$
obtained by considering the second boundary of $\tau^L$ and all triangles of $\tau^L$ that
do not appear in $\tau^L_n$. We view $\check\tau^L_n$ as a rooted planar map whose root
edge is the root of the second boundary of $\tau^L$ and call
$\check\tau^L_n$ the revealed region at step $n$. The revealed region $\check\tau^L_n$ is 
given with a ``boundary'' that consists of all its edges that are incident to $\tau^L_n$. We also define the revealed region 
at step $\zeta_L$ by adding to $\check\tau^L_{\zeta_L-1}$ the revealed triangle at step $\zeta_L$
(which necessarily has a vertex on the first boundary of $\tau^L$), and defining the
boundary in the obvious manner. Then, at each step $n<\zeta_L$, the choice of 
the revealed edge is made on the boundary of $\check \tau^L_n$ as a deterministic function of $\check\tau^L_n$. Furthermore, 
$\tau^L_n$ is viewed as a triangulation with two boundaries, and the first root edge is the 
same as in $\tau^L$, whereas the second root edge is chosen on the boundary of $\check \tau^L_n$ as a deterministic function of $\check\tau^L_n$.

For $n<\zeta_L$, write $P^L_n$ for the size of the second boundary of $\tau^L_n$
and set $P^L_n=\dagger$ for $n\geq \zeta_L$, where $\dagger$ serves as a cemetery point. The discussion  of the
previous section shows that, conditionally on the event $\{n<\zeta_L,P^L_n=k\}$,
$\tau^L_n$ is distributed as a Boltzmann triangulation of $\T^2(L,k)$. 
Furthermore, $(P^L_n)_{n\geq 0}$ is a Markov chain with
values in $\N\cup\{\dagger\}$ with transition probabilities
$$\P(P^L_{n+1}=p+1\mid P^L_n=p)= \frac{1}{12\sqrt3} \frac{Z'(L,p+1)}{Z'(L,p)}=:q_L(p,p+1)$$
and, for every $m\in\{0,1,\ldots,p-1\}$,
$$\P(P^L_{n+1}=p-m\mid P^L_n=p)=2\,\frac{Z'(L,p-m)}{Z'(L,p)}\,Z(m+1)=:q_L(p,p-m),$$
and finally 
$$\P(P^L_{n+1}=\dagger\mid P^L_n=p)=1-q_L(p,p+1)- \sum_{m=0}^{p-1}q_L(p,p-m)=:q_L(p,\dagger).$$
Now recall \eqref{peel1} and \eqref{peel2}, and use the notation 
$q_\infty(p,j)$ for the transition probabilities of the peeling process 
of the UIPT of type I, which is discussed in \cite[Section 6.1]{CLG}. According to \cite{CLG}, the nonzero values
of $q_\infty(p,j)$ are determined as follows.
For every integer $p\geq 1$, we have
$$q_\infty(p,p+1)=\frac{1}{12\sqrt{3}}\,\frac{C^{(1)}(p+1)}{C^{(1)}(p)},$$
and, for every $m\in\{0,1,\ldots,p-1\}$,
$$q_\infty(p,p-m)=2\,Z(m+1)\,\frac{C^{(1)}(p-m)}{C^{(1)}(p)}.$$
Comparing the last two displays with \eqref{peel1} and \eqref{peel2}, we see that we have
$$q_L(p,p+j)=\frac{L+p}{L+p+j}\,q_\infty(p,p+j)$$
for every $j\in\{1,0,-1,-2,\ldots,-p+1\}$. In other words, $(P^L_n)_{n\geq 0}$ is a $h$-transform of the 
peeling process 
of the UIPT, for the function $h=h_L$ defined by
$$h_L(j):=\frac{L}{L+j},$$
for $j=1,2,\ldots$. 

\subsection{Asymptotics for the peeling process} 

We now want to derive asymptotics when $L\to\infty$ (the integer $p_0$ remains fixed). We set $Z_L=P^L_{\zeta_L-1} + 1$ (thus $L+Z_L$ is interpreted as 
the boundary size of the triangulation that needs to be ``pasted'' to the revealed region at time $\zeta_L$ to
recover $\tau^L$). We also write $(P^\infty_n)_{n\geq 0}$ for the Markov chain with transition probabilities $q_\infty(p,j)$ started at $p_0$, which is known to be transient \cite{CLG}.

\begin{proposition}
\label{asymptotics-Z}
We have 
$$\frac{Z_L}{L} \build{\la}_{L\to\infty}^{\rm(d)} \Lambda,$$
where $\Lambda$ has density $\frac{3}{2}\,(1+x)^{-5/2}$ on $\R_+$.
\end{proposition}

\proof
For every integer $j\geq 1$, we get, using the $h$-transform relation 
between the Markov chains $(P^L_n)_{n\geq 0}$ and $(P^\infty_n)_{n\geq 0}$,
\begin{align}
\label{asymp-Z}
\P(Z_L=j+1)&=\sum_{n=0}^\infty \P(P^L_n=j,\zeta_L=n+1)\nonumber\\
&=\sum_{n=0}^\infty  \P(P^L_n=j)\,q_L(j,\dagger)\nonumber\\
&=\frac{1}{h_L(p_0)}\sum_{n=0}^\infty  \P(P^\infty_n=j)\,h_L(j)\,q_L(j,\dagger)\nonumber\\
&=\frac{1}{h_L(p_0)}\, U(p_0,j)\,h_L(j)\,q_L(j,\dagger),
\end{align}
where we have written $U(k,\ell)$ for the potential kernel
of the Markov chain $(P^\infty_n)_{n\geq 0}$. We can explicitly compute $U(p_0,j)$ when $j\geq p_0$. To this end, set for every 
integer $k\in\{1,0,-1,-2,\ldots\}$,
$$q_k:=\lim_{p\to\infty} q_\infty(p,p+k)=\left\{\begin{array}{ll}
\displaystyle{\frac{1}{\sqrt{3}}}&\hbox{if } k=1,\\
\noalign{\medskip}
\displaystyle{2\,Z(k+1)\,12^{-k}}\quad&\hbox{if }k\leq 0.
\end{array}
\right.
$$
From \cite{CLG}, $(q_k)_{k\leq 1}$ defines a probability measure with mean zero on $\Z$.
Let $(S_n)_{n\geq 0}$ denote the (recurrent) random walk with jump distribution $(q_k)_{k\leq 1}$, and
set $T^S_0=\inf\{n\geq 0:S_n\leq 0\}$. Consider the killed random walk $(S^\bullet_n)_{n\geq 0}$
defined by $S^\bullet_n=S_n$ if $n<T^S_0$ and $S^\bullet_n=\dagger$ if $n\geq T^S_0$. 
According to \cite{CLG}, $(P^\infty_n)_{n\geq 0}$ 
is the $h$-transform of the Markov chain $(S^\bullet_n)_{n\geq 0}$, for the function
$$h_\bullet(p):=12^{-p} C^{(1)}(p)$$
for every $p\geq 1$.
For the random walk $S$ started from $0$, the expected number of visits 
of $j\geq 1$ before the first return to $0$ is equal to $1$, and is also equal 
to $1/\sqrt{3}$ times the expected number of visits of $j$ for the Markov chain
$S^\bullet$ started from $1$. So, if $U^\bullet$ denotes the potential kernel
of $S^\bullet$, we have $U^\bullet(1,j)=\sqrt{3}$ for every $j\geq 1$. The 
$h$-transform relation 
between the Markov chains $(P^\infty_n)_{n\geq 0}$ and $(S^\bullet_n)_{n\geq 0}$
then gives
$$U(1,j)=\frac{h_\bullet(j)}{h_\bullet(1)}\,\sqrt{3},$$
and since $(P^\infty_n)_{n\geq 0}$ is transient and its positive jumps are of size $1$, it
is immediate that $U(p,j)=U(1,j)$ whenever $p\leq j$. 
From \eqref{triangu11}, we have $h_\bullet(1)=1/(72\sqrt{2\pi})$ and
$$h_\bullet(p)\build{\sim}_{p\to\infty}^{} \frac{1}{36\pi\sqrt{2}}\,\sqrt{p}.$$
It follows that,
for $j\geq p_0$,
$$U(p_0,j)= 72\sqrt{6\pi}\,h_\bullet(j)\build{\sim}_{j\to\infty}^{} \frac{2\sqrt{3}}{\sqrt{\pi}}\,\sqrt{j}.$$
Now recall formula \eqref{asymp-Z}. From the end of Section \ref{peeling-proba}, we
know that $q_L(j,\dagger)$ behaves like
$$\frac{\sqrt{3\pi}}{4}\,\sqrt{\frac{L}{j}}\,(L+j)^{-3/2}$$
when both $L$ and $j$ are large. We thus get
$$\P(Z_L=j+1)\build{\sim}_{L,\,j\to\infty}^{} \frac{2\sqrt{3}}{\sqrt{\pi}}\,\sqrt{j}
\times \frac{L}{L+j}\times\frac{\sqrt{3\pi}}{4}\,\sqrt{\frac{L}{j}}\,(L+j)^{-3/2}=\frac{3}{2}\,L^{3/2}\,(L+j)^{-5/2}.$$
The result of the proposition follows. \endproof

\subsection{Convergence of rescaled triangulations}
\label{conve-tria}
For every integer $L\geq 1$, let $\t'_L$ be a Boltzmann distributed rooted triangulation with a simple boundary of size $L$.
We let $\dg$ stand for the graph distance on the vertex set $V(\t'_L)$.
We write $\partial\t'_L$ for the set of all boundary vertices and
we denote the set of all inner vertices by $V_i(\t'_L):=V(\t'_L)\backslash \partial\t'_L$. We also let 
$\nu'_L$ be the counting measure on $V_i(\t'_L)$ scaled by the factor $\frac{3}{4}L^{-2}$. 
We finally consider the ``boundary path'' $\Theta'_L=
(\Theta'_L(k))_{0\leq k\leq L}$, which is obtained by letting
$\Theta'_L(0)=\Theta'_L(L)$ be the root vertex of $\t'_L$ and then letting
$\Theta'_L(1),\Theta'_L(2),\ldots,\Theta'_L(L-1)$ be the points of
$\partial \t'_L$ enumerated in clockwise order from $\Theta'_L(0)$. We also set $\wh\Theta'_L(t)=\Theta'_L(\lfloor Lt\rfloor)$
for $t\in[0,1]$. According to  Theorem 1.1 of \cite{AHS} we have
\begin{equation}
\label{convGHPU-unpointed}
(V(\t'_L),\sqrt{3/2}\,L^{-1/2}\dg,\nu'_L,\wh\Theta'_L)\build{\la}_{L\to\infty}^{(d)} (\D',D',\mathbf{V}',\Gamma').
\end{equation}
where the convergence holds in distribution in $(\M^{GHPU},\dd_{GHPU})$, and $(\D',D',\mathbf{V}',\Gamma')$ 
is a curve-decorated free Brownian disk of perimeter $1$ 
(the precise definition of this limiting space will be given below in Section \ref{sec:limspace}). 
We note that $\E[\mathbf{V}'(\D')]=1$ (the density of $\mathbf{V}'(\D')$ is the function $r\mapsto (2\pi r^5)^{-1/2}\exp(-1/(2r))$,
cf. \cite[Section 1.5]{BM}). 

The 
convergence \eqref{convGHPU-unpointed} seems to be stated incorrectly since the
paths $\wh\Theta'_L$ are {\it not continuous} 
and thus the random space $(V(\t'_L),\sqrt{3/2}\,L^{-1/2}\dg,\nu'_L,\wh\Theta'_L)$ does not belong to 
$\M^{GHPU}$. There is however a straightforward way of 
overcoming this difficulty, by replacing $V(\t'_L)$ with the union of all 
its edges, each edge being represented by a copy of the interval $[0,1]$, so that the
boundary path can be made continuous and its range will be the union
of the boundary edges --- see \cite{AHS} for more details.

We now want to argue that a result similar to \eqref{convGHPU-unpointed} holds for 
rooted {\it and pointed} triangulations. So, for every integer $L\geq 1$, let $\t_L$ be a Boltzmann distributed rooted
and pointed triangulation with a simple boundary of size $L$. We define $V_i(\t_L)$, $\nu_L$ and $\Theta_L$
in the same way as $V_i(\t'_L)$, $\nu'_L$ and $\Theta'_L$ were defined above, and we also write $v_*^{(L)}$
for the distinguished vertex of $\t_L$. Then, we have
\begin{equation}
\label{convGHPU-pointed}
\big(V(\t_L),\sqrt{3/2}\,L^{-1/2}\dg,\nu_L,\wh\Theta_L,v_*^{(L)}\big)\build{\la}_{L\to\infty}^{(d)} \big(\D,D,\mathbf{V},\Gamma,\xx_*\big),
\end{equation}
where the convergence holds in distribution in $(\M^{GHPU\bullet},\dd_{GHPU\bullet})$, and the limit $(\D,D,\mathbf{V},\Gamma,\xx_*)$ is now a curve-decorated free pointed Brownian disk of perimeter $1$ (see Section \ref{sec:limspace} below).

Let us explain why \eqref{convGHPU-pointed}
follows from \eqref{convGHPU-unpointed}. To simplify notation, write $\mathbf{X}^L$, respectively $\mathbf{X}'^L$, for the space $(V(\t_L),\sqrt{3/2}\,L^{-1/2}\dg,\nu_L,\wh\Theta_L,v_*^{(L)})$, resp. for $(V(\t'_L),\sqrt{3/2}\,L^{-1/2}\dg,\nu'_L,\wh\Theta'_L)$. Also write 
$\langle\nu'_L,1\rangle$ for the total mass of $\nu'_L$. Then, for every bounded continuous function $F$
on $(\M^{GHPU\bullet},\dd_{GHPU\bullet})$,
$$\E[F(\mathbf{X}^L)]=\frac{\E\Big[\sum_{x\in V_i(\t'_L)} F((\mathbf{X}'^L,x))\Big]}{\E[\#V_i(\t'_L)]}
=\frac{\E\Big[\int \nu'_L(dx)F((\mathbf{X}'^L,x))\Big]}{\E[\langle\nu'_L,1\rangle]},$$
where $(\mathbf{X}'^L,x)$ obviously denotes the pointed space derived from $\mathbf{X}'_L$ by distinguishing the point $x$. 
We then claim that
\begin{equation}
\label{tec-point2}
\E[\langle\nu'_L,1\rangle] \build{\la}_{L\to\infty}^{} 1=\E[\mathbf{V}'(\D')].
\end{equation}
Assuming that \eqref{tec-point2} holds, we can apply Lemma \ref{pointed-notpointed}, which implies that
$\mathbf{X}^L$ converges in distribution (in the space $(\M^{GHPU\bullet},\dd_{GHPU\bullet})$)
to the random space $\mathbf{X}^\infty$ whose law is characterized by
$$\E[F(\mathbf{X}^\infty)]= \E\Big[\int \mathbf{V}'(\dd x)\,F((\D',D',\mathbf{V}',\Gamma',x))\Big].
$$
The last display exactly means that $\mathbf{X}^\infty$ is a (curve-decorated) free pointed Brownian disk of perimeter $1$ --- see e.g. the discussion
in \cite[Section 6]{Repre}. It only remains to justify our claim \eqref{tec-point2}. We already know
(by \eqref{convGHPU-unpointed}) that $\langle\nu'_L,1\rangle$ converges in distribution to $\mathbf{V}'(\D')$,
and therefore it suffices to verify that $\E[\langle\nu'_L,1\rangle\,\mathbf{1}_{\{\langle\nu'_L,1\rangle\geq a\}}]$
tends to $0$ as $a\to +\infty$, uniformly in $L$. This can be checked 
from the explicit formulas \eqref{triangu1},\eqref{triangu11},\eqref{triangu12} and we omit the details. 

\section{The limiting space}
\label{sec:limspace}

\subsection{The Bettinelli-Miermont construction}
\label{Bet-Mie}

In this section, we recall the Bettinelli-Miermont construction of the free Brownian disk \cite{Bet,BM},
which will play an important role in our proofs. We follow the presentation 
of Section 6 in \cite{spine}, which is slightly different from \cite{Bet,BM}.  

We fix $\xi>0$, which will correspond to the boundary size
of the Brownian disk. We consider a Brownian excursion $(\be_t)_{0\leq t\leq \xi}$ of duration $\xi$,
and, conditionally on $(\be_t)_{0\leq t\leq \xi}$, a Poisson point measure $\n=\sum_{i\in I}\delta_{(t_i,\omega_i)}$ on $[0,\xi]\times \S$
with intensity
$$2\,\dd t\,\N_{\sqrt{3}\,\be_t}(\dd \omega).$$
We let $\mathfrak{T}$ be
the compact metric space obtained from
the disjoint union
\begin{equation}
\label{tree-disk}
[0,\xi] \cup \Big(\bigcup_{i\in I} \t_{(\omega_i)}\Big)
\end{equation}
by identifying the root $\rho_{(\omega_i)}$ of $\t_{(\omega_i)}$
with the point $t_i$ of $[0,\xi]$, for every $i\in I$. The metric $\dd_\mathfrak{T}$ on $\mathfrak{T}$ is defined as follows. 
First, the restriction of $\dd_\mathfrak{T}$ to each tree $\t_{(\omega_i)}$ is 
the metric $d_{(\omega_i)}$. Then, if $u,v\in[0,\xi]$, we take $\dd_\mathfrak{T}(u,v):=|v-u|$. If $u\in [0,\xi]$, and $v\in\t_{(\omega_i)}$ for some $i\in I$, we take $\dd_\mathfrak{T}(u,v):=|u-t_i|+ d_{(\omega_i)}(\rho_{(\omega_i)},v)$. Finally
  if $u\in \t_{(\omega_i)}$ and $v\in \t_{(\omega_j)}$, with $j\not =i$, we let
$$\dd_\mathfrak{T}(u,v):= d_{(\omega_i)}(u,\rho_{(\omega_i)})+|t_i-t_j|+ d_{(\omega_j)}(\rho_{(\omega_j)},v).$$
The volume measure on $\mathfrak{T}$ is the sum of the volume measures on the trees $\t_{(\omega_i)}$, $i\in I$.

If $\Sigma:=\sum_{i\in I}\sigma(\omega_i)$ is the total mass of the volume measure, we define a clockwise exploration $(\ee_t)_{0\leq t\leq\Sigma}$
of $\mathfrak{T}$, informally by concatenating the mappings $p_{(\omega_i)}:[0,\sigma(\omega_i)]\la \t_{(\omega_i)}$ in the 
order prescribed by the $t_i$'s. To give a more precise definition, set
$$\beta_s:= \sum_{i\in I} \mathbf{1}_{\{t_i\leq s\}}\,\sigma(\omega_i)\,,\ \beta_{s-}:= \sum_{i\in I} \mathbf{1}_{\{t_i< s\}}\,\sigma(\omega_i)\,,$$
for every $s\in[0,\xi]$. 
Then, for every $t\in[0,\Sigma]$, we define $\ee_t\in\mathfrak{T}$ as follows. We observe that there is a unique $s\in[0,\xi]$ such that $\beta_{s-}\leq t\leq \beta_{s}$, and:
\begin{description}
\item[$\bullet$] Either there is a (unique) $i\in I$ such that $s=t_i$, and we set
$\ee_t:= p_{(\omega_i)}(t-\beta_{t_i-})$.
\item[$\bullet$] Or there is no such $i$ and we set $\ee_t:=s$.
\end{description}
Note that $\ee_0=0$ and $\ee_\Sigma=\xi$.

The clockwise exploration allows us to define ``intervals''
in $\mathfrak{T}$. Let us make the convention that, if 
$s,t\in[0,\Sigma]$ and $s>t$, the (real) interval $[s,t]$ is defined by $[s,t]:=[s,\Sigma]\cup [0,t]$ (of course, if 
$s\leq t$, $[s,t]$ is the usual interval).  
Then, for every $u,v\in\mathfrak{T}$, such that $u\not =v$, there is a smallest interval $[s,t]$, with $s,t\in[0,\Sigma]$, such that
$\ee_s=u$ and $\ee_t=v$, and we define  
$$[|u,v|]:=\{\ee_r:r\in[s,t]\}.$$
Observe that in general  $[|u,v|]\not =[|v,u|]$. We also take $[|u,u|]=\{u\}$. Note that we use the notation 
$[|u,v|]$ rather than $[u,v]$ to avoid confusion with intervals of the real line.

We then assign labels $(\ell_a)_{a\in\mathfrak{T}}$ to the points of $\mathfrak{T}$.
 If $a=s\in[0,\xi]$, we
take $\ell_a:=\sqrt{3}\,\be_s$, and if $a\in \t_{(\omega_i)}$ for some $i\in I$, we 
simply let $\ell_a$ be the label of $a$ in $ \t_{(\omega_i)}$. The function $a\mapsto\ell_a$
is continuous on $\mathfrak{T}$. The following simple fact will be important for us: For every $\ve>0$, 
formula \eqref{hittingpro} and the property $\int_0^\ve \dd t/(\be_t)^2=\infty$ imply that some
of the trees $\t_{(\omega_i)}$ such that $t_i<\ve$ carry negative labels.

For every $a,b\in\mathfrak{T}$, we set
$$D^\circ(a,b):=\ell_a+\ell_b - 2\max\Big(\min_{c\in[|a,b|]}\ell_c,\min_{c\in[|b,a|]}\ell_c\Big).$$
Notice that $D^\circ(0,\xi)=0$ (because $\ell_0=\ell_\xi=0$ and the ``interval'' $[|\xi,0|]$ is the pair $\{0,\xi\}$).
We define a pseudo-metric on $\mathfrak{T}$ by setting, for every $a,b\in\mathfrak{T}$, 
$$D(a,b):=\inf_{a_0=a,a_1,\ldots,a_{p-1},a_p=b}\,\sum_{i=1}^p D^\circ(a_{i-1},a_i)$$
where the infimum is over all choices of the integer $p\geq 1$
and of the points $a_1,\ldots,a_{p-1}$ in $\mathfrak{T}$. 
One can prove \cite[Theorem 13]{Bet} that $D(a,b)=0$ if and only if $D^\circ(a,b)=0$ (the ``if'' part is trivial).
We set $\D:=\mathfrak{T}/\{D=0\}$, where the notation $\mathfrak{T}/\{D=0\}$ refers to the quotient space
of $\mathfrak{T}$ for the equivalence relation defined by setting $a\simeq
 b$ if and only if $D(a,b)=0$, and this quotient space is equipped with the metric induced by $D$ --- similar notation will be used several times in what follows. It is immediate that
$D(a,b)\geq|\ell_a-\ell_b|$, and therefore $D(a,b)=0$ implies that $\ell_a=\ell_b$,
so that we can make sense of labels on $\D$, for which we keep the
same notation $\ell_x$. We write $\Pi$ for the canonical projection 
from $\mathfrak{T}$ onto $\D$. 
The volume measure $\mathbf{V}$ on $\D$ is the pushforward of the volume measure on $\mathfrak{T}$ under $\Pi$. 
The metric space $(\D,D)$ is a.s.~homeomorphic to the closed unit disk of the plane \cite{Bet}, and, in any such 
homeomorphism, the unit circle corresponds to the ``boundary'' $\partial\D:=\Pi([0,\xi])$ (which is therefore 
the set of all points of $\D$ that have no neighborhood
homeomorphic to the open unit disk).

There is a unique point $a_*\in\mathfrak{T}$ such that $\ell_{a_*}=\min\{\ell_a:a\in\mathfrak{T}\}<0$
and we set $\xx_*=\Pi(a_*)$. For every $x\in \D$, we have $D(\xx_*,x)=\ell_x-\ell_{\xx_*}$. In particular, the
distance from $\xx_*$ to $\partial \D$ is $-\ell_{\xx_*}$, and $\Pi(0)=\Pi(\xi)$ is  the unique point of
$\partial \D$ at minimal distance from $x_*$.

The {free pointed Brownian disk with perimeter} $\xi$ may then be defined as 
the random measure  metric space $(\D,D,\mathbf{V})$
with the distinguished point $\xx_*$ but, for our purposes, it will be convenient
to view $\D$ as a curve-decorated space.
We first observe that the mapping
$[0,\xi]\ni t\mapsto \Pi(t)$ is a simple loop 
(recall that $\Pi(0)=\Pi(\xi)$) whose range is $\partial \D$. More precisely, the loop $[0,\xi]\ni t\mapsto \Pi(t)$ is a {\it standard boundary curve} in the following sense.
We first recall from
 \cite[Theorem 9]{Repre} that the measures $\ve^{-2}\mathbf{1}_{\{D(x,\partial \D)\leq \ve\}}
\mathbf{V}(\dd x) $ converge weakly (a.s.) to a measure on the boundary $\partial \D$, which we denote by $\mu_{\partial \D}$
and whose total mass is the perimeter $\xi$ of $\D$ (the measure $\mu_{\partial \D}$ is known as the uniform measure on the boundary of $\D$).
We then say that $f:[0,\xi]\to \partial \D$ is a standard boundary curve of $\D$ if $f$ is a simple loop 
whose range is $\partial\D$, and if the pushforward of Lebesgue measure on $[0,\xi]$ under $f$ is $\mu_{\partial \D}$.
The latter property is equivalent to
\begin{equation}\label{def:bou:curve}
t= \lim \limits_{\ve\to 0}\ve^{-2}\int \mathbf{V}(\dd x) \mathbf{1}_{\{D(x,f([0,t]))\leq \ve\}}, \quad \text{ for every } t\in[0,\xi].
\end{equation}
It then follows from
 \cite[Theorem 9]{Repre} that the loop $[0,\xi]\ni t\mapsto \Pi(t)$ is a standard boundary curve, and the same holds for the
time-reversed loop $\check\Pi(t):=\Pi(\xi-t)$. Moreover, for every $x\in\partial\D$, there are exactly two standard boundary curves
with starting point $x$, which are obtained by changing the origin of the loops $\Pi(t)$ and $\check\Pi(t)$. 

We observe that the starting point $\Pi(0)$ of the loop $t\mapsto \Pi(t)$ is not a ``typical'' point of $\partial \D$, 
since it is the point of $\partial\D$ at minimal distance from $\xx_*$.
So we will consider  the loop $t\mapsto\Pi(t)$  ``re-rooted at a uniform boundary point''. To this end,
let $\mathfrak{U}$ 
be uniformly distributed over $[0,\xi]$ and independent of the random quantities involved in the definition of
$\D$. We set $\Gamma(t)=\Pi(\mathfrak{U}+t)$ for $t\in [0,\xi-\mathfrak{U}]$, and $\Gamma(t)=\Pi(\mathfrak{U}+t-\xi)$
for $t\in(\xi-\mathfrak{U},\xi]$. Then $\Gamma$ is again a standard boundary curve (now rooted at a uniform
boundary point). 

The {\it curve-decorated free pointed Brownian disk with perimeter $\xi$} that appears (for $\xi=1$) in formula \eqref{convGHPU-pointed} is the random space
$(\D,D,\mathbf{V},\Gamma,\xx_*)$, which is a random variable taking values in $\M^{GHPU\bullet}$. The
curve-decorated free Brownian disk
with perimeter $\xi$ (appearing in \eqref{convGHPU-unpointed}) is then the random space $(\D',D',\mathbf{V}',\Gamma')$ in $\M^{GHPU}$ whose distribution
is characterized by
\begin{equation}
\label{pointed/nonpointed}
\E[F((\D',D',\mathbf{V}',\Gamma'))]= \xi^2\, \E\Big[ \frac{F((\D,D,\mathbf{V},\Gamma))}{\mathbf{V}(\D)}\Big].
\end{equation}
See the discussion at the beginning of \cite[Section 6]{Repre}. It will sometimes be convenient to  ``forget'' the curve $\Gamma'$
and to keep track only of its initial point: The pointed measure metric space $(\D',D',\mathbf{V}',\Gamma'(0))$ is the
free Brownian disk of perimeter $\xi$ {\it pointed at a uniform boundary point}, which is discussed in
\cite[Section 6]{Repre} (informally, given the unpointed space $(\D^\prime, D^\prime, \mathbf{V}^\prime)$, the 
distinguished point $\Gamma'(0)$ is chosen according to the uniform measure on $\partial\D'$).

Let us finally discuss simple geodesics in $\D$. Let $x=\Pi(a)$ be a point of $\D$, and $r\in[0,\Sigma]$
such that $a=\ee_r$. For every $t\leq \ell_x-\ell_*=D(\xx_*,x)$, set
$$\varphi_r(t):=\left\{\begin{array}{ll}
\inf\{s\in[r,\Sigma]:\ell_{\ee_s}=\ell_x-t\}\quad&\hbox{if }\{s\in[r,\Sigma]:\ell_{\ee_s}=\ell_x-t\}\not=\varnothing,\\
\inf\{s\in[0,r]:\ell_{\ee_s}=\ell_x-t\}&\hbox{otherwise.}
\end{array}
\right.
$$
Then $\big(\Pi(\ee_{\varphi_r(t)})\big)_{0\leq t\leq D(\xx_*,x)}$ is a geodesic from $x$ to $\xx_*$, which is
called a simple geodesic. It is easy to verify that, if $x=\Pi(a)$ and $y=\Pi(b)$ are two points of $\D$,
there are two simple geodesics starting from $x$ and from $y$ respectively that coalesce at a 
point whose label is 
$$\max\Big(\min_{c\in[|a,b|]}\ell_c,\min_{c\in[|b,a|]}\ell_c\Big).$$
Consequently, the quantity $D^\circ(a,b)$ is the length of a path from $\Pi(a)$ to $\Pi(b)$ obtained by concatenating
two simple geodesics up to the point where they merge.

\subsection{Construction of the limiting space}
\label{conf-lim-space}

We will now slightly modify the preceding construction of the Brownian disk to get another 
random metric space, which later will be identified to a particular subset of $\D$ equipped with
an intrinsic metric. We start from the same Poisson point measure $\n=\sum_{i\in I}\delta_{(t_i,\omega_i)}$ as in
the previous section, but, for every $i\in I$, we now consider the truncation  $\wt \omega_i:=\mathrm{tr}_0(\omega_i)$ of $\omega_i$ at level $0$, and we write
$\mathrm{tr}_0(\n)=\sum_{i\in I}\delta_{(t_i,\tilde\omega_i)}$.
Recall that the genealogical tree $\t_{(\tilde\omega_i)}$ is identified to the subtree obtained from $\t_{(\omega_i)}$
by pruning the branches when labels first hit $0$. We then consider the geodesic space
$\mathfrak{T}^\star$ which is obtained from the disjoint union
$$[0,\xi] \cup \Bigg(\bigcup_{i\in I} \t_{(\tilde\omega_i)}\Bigg)$$
by identifying the root of $\t_{(\tilde\omega_i)}$ with the point $t_i$ of $[0,\xi]$, for every $i\in I$. Then
$\mathfrak{T}^\star$ is identified to a closed subset of $\mathfrak{T}$, so that labels $\ell_a$
are defined for every $a\in\mathfrak{T}^\star$ and we can
equip $\mathfrak{T}^\star$ with the restriction of the distance on $\mathfrak{T}$. We also define the volume measure on $\mathfrak{T}^\star$ as the sum of the volume measures on the trees $\t_{(\tilde\omega_i)}$, $i\in I$.

We note that
labels $\ell_a$ are
nonnegative for every $a\in\mathfrak{T}^\star$, because we have replaced $\omega_i$ by $\mathrm{tr}_0(\omega_i)$. We 
define the boundary of $\mathfrak{T}^\star$ by $\partial\mathfrak{T}^\star:=\{a\in\mathfrak{T}^\star: \ell_a=0\}$ and note that $0,\xi\in\partial\mathfrak{T}^\star$.

We can introduce a clockwise exploration process $(\ee^\star_s)_{0\leq s\leq \Sigma^\star}$ of $\mathfrak{T}^\star$ in exactly the same way
as we did for $\mathfrak{T}$, and use this exploration process to define the ``interval'' $[|a,b|]_\star$ for every $a,b\in\mathfrak{T}^\star$
(we have in fact $[|a,b|]_\star=[|a,b|]\cap\mathfrak{T}^\star$). We then set, for every $a,b\in\mathfrak{T}^\star\backslash\partial\mathfrak{T}^\star$,
\begin{equation}
\label{Dbar1}
D_\star^\circ(a,b):=\ell_a+\ell_b - 2\max\Big(\min_{c\in[|a,b|]_\star}\ell_c,\min_{c\in[|b,a|]_\star}\ell_c\Big)
\end{equation}
if the maximum in the right-hand side is positive, and $D_\star^\circ(a,b):=\infty$ otherwise. 
Finally, in exactly the same way as we defined $D$ from $D^\circ$, we set,
for every $a,b\in\mathfrak{T}^\star\backslash\partial\mathfrak{T}^\star$, 
\begin{equation}
\label{Dbar2}
D_\star(a,b):=\inf_{a_0=a,a_1,\ldots,a_{p-1},a_p=b}\sum_{i=1}^p D_\star^\circ(a_{i-1},a_i)
\end{equation}
where the infimum is over all choices of the integer $p\geq 1$
and of the points $a_1,\ldots,a_{p-1}$ in $\mathfrak{T}^\star\backslash\partial\mathfrak{T}^\star$. It is not hard to
verify that $D_\star(a,b)<\infty$ (see Proposition 30 (i) in \cite{Disks} for a very similar
argument). The mapping $(a,b)\mapsto D_\star(a,b)$
is continuous on $(\mathfrak{T}^\star\backslash\partial\mathfrak{T}^\star)\times (\mathfrak{T}^\star\backslash\partial\mathfrak{T}^\star)$, and we have again
$D_\star(a,b)\geq |\ell_a-\ell_b|$. We also notice that $D_\star^\circ(a,b)=D^\circ(a,b)$
whenever $D_\star^\circ(a,b)<\infty$ (if labels are positive on $[|a,b|]_\star$, this implies that $[|a,b|]_\star=[|a,b|]$). Consequently,
$D_\star(a,b)\geq D(a,b)$ for every $a,b\in \mathfrak{T}^\star\backslash\partial\mathfrak{T}^\star$.

\begin{proposition}
\label{def-lim-space}
Almost surely, the function $(a,b)\mapsto D_\star(a,b)$ has a continuous extension to $\mathfrak{T}^\star\times \mathfrak{T}^\star$, 
which satisfies the triangle inequality and the bound $D_\star(a,b)\geq |\ell_a-\ell_b|$ for every $a,b\in \mathfrak{T}^\star$. Furthermore,
the property $D_\star(a,b)=0$ holds if and only if either $a$ and $b$ both belong to
$\mathfrak{T}^\star\backslash\partial\mathfrak{T}^\star$ and $D^\circ_\star(a,b)=0$, or $a$ and $b$
both belong to $\partial\mathfrak{T}^\star$ and we have $\{a,b\}=\{\ee^\star_s,\ee^\star_t\}$,
with $0\leq s\leq t\leq \Sigma^\star$, and $\ell_{\ee^\star_r}>0$ for every $r\in(s,t)$.
\end{proposition}

\proof Let us start by proving the first assertion. Since $D_\star$ satisfies the triangle inequality, it is enough to
verify that, for any $a\in \partial\mathfrak{T}^\star$, if $(a_n)_{n\in\N}$
is a sequence in $\mathfrak{T}^\star\backslash\partial\mathfrak{T}^\star$ that 
converges to $a$, we have $D_\star(a_n,a_m)\la 0$ as $n,m\to\infty$. 
If $a\in \partial\mathfrak{T}^\star\backslash\{0,\xi\}$, this is a straightforward consequence of Lemma \ref{tech-lem}.
Indeed, let $i\in I$ such that $a\in\t_{(\tilde\omega_i)}$. Then, we have also $a_n\in\t_{(\tilde\omega_i)}$
when $n$ is large. However, when $a_n$ and $a_m$ both belong to $\t_{(\tilde\omega_i)}$, we immediately get from our definitions that 
$D_\star(a_n,a_m)\leq \Delta_{(\tilde\omega_i)}(a_n,a_m)$
with the notation introduced before Lemma \ref{tech-lem}. The convergence of  $D_\star(a_n,a_m)$ to $0$ now follows
as a consequence of the first assertion of Lemma \ref{tech-lem}.

Let us turn to the case where $a=0$ or $a=\xi$. For definiteness, take $a=0$ (the 
other case is similar). We can list all indices $i\in I$ such that $t_i\in[0,\xi/2]$ and $W_*(\omega_i)\leq 0$
in a sequence $i_1,i_2,\ldots$ such that $\xi/2>t_{i_1}>t_{i_2}>\cdots$. We set $h_j=\sqrt{3}\,\be_{t_{i_j}}$
for every $j\geq 1$. Then, conditionally on $(\be_t)_{0\leq t\leq \xi}$ and on the sequence $(t_{i_1},t_{i_2},\ldots)$,
the snake trajectories
$\omega_{i_1},\omega_{i_2},\ldots$ are independent and the conditional distribution 
of $\omega_{i_j}$ is $\N^{[0]}_{h_j}$. As already mentioned, we have for every
$a,b\in \t_{(\tilde\omega_{i_j})}\backslash \partial \t_{(\tilde\omega_{i_j})}$,
\begin{equation}
\label{bd-tec-lem}
D_\star(a,b) \leq \Delta_{(\tilde\omega_{i_j})}(a,b),
\end{equation}
For every $j\geq 1$,
set
$$H_j:=\sup\big\{\Delta_{(\tilde\omega_{i_j})}(a,b):a,b\in \t_{(\tilde\omega_{i_j})}\backslash \partial \t_{(\tilde\omega_{i_j})}\big\}.$$
Thanks to Lemma \ref{tech-lem} and a scaling argument, we have $H_j=h_jH'_j$
where the random variables $H'_j$
are independent of $( \be,t_{i_1},t_{i_2},\ldots)$ and have the same distribution with finite expectation. Next observe that
$$\E\Big[ \sum_{j=1}^\infty h_j\,\Big|\, \be\Big]= 2\int_0^{\xi/2} \dd t\,\sqrt{3}\,\be_t\,\N_{\sqrt{3}\be_t}(W_*\leq 0)
=3\sqrt{3} \int_0^{\xi/2} \frac{\dd t}{\be _t}<\infty\quad\hbox{a.s.}$$
and thus
$$\sum_{j=1}^\infty h_j<\infty\quad\hbox{a.s.}$$
It also follows that
$$\sum_{j=1}^\infty H_j<\infty\quad\hbox{a.s.}$$
since 
$$\E\Big[ \sum_{j=1}^\infty H_j\,\Big|\, \be,t_{i_1},t_{i_2},\ldots\Big]=\sum_{j=1}^\infty h_j\,\E[H'_j]= C\,\sum_{j=1}^\infty h_j.$$
where the constant $C$ is as in Lemma \ref{tech-lem}.
For every $j\geq 1$, write $\rho_j=t_{i_j}$  for the root of $\t_{(\tilde\omega_{i_j})}$. Fix an integer $k\geq 1$, and 
let $a\in [|0,\rho_{k}|]_\star$ such that $\ell_a>0$. Then there is an index $j>k$ such that either $a\in \t_{(\tilde\omega_{i_j})}$ or
$a\in [|\rho_j,\rho_{j-1}|]_\star$. We observe that in both cases we have
\begin{equation}
\label{bd-tec-lem2}
D_\star(a,\rho_j)\leq H_j + \ell_a.
\end{equation}
If $a\in \t_{(\tilde\omega_{i_j})}$, this is immediate from the bound \eqref{bd-tec-lem} and the definition of $H_j$. 
If $a\in [|\rho_j,\rho_{j-1}|]_\star$, we choose $\ve>0$ smaller than the minimal label 
in $[|\rho_j,\rho_{j-1}|]_\star$ and we write $a_{(\ve)}$ for the last vertex 
(in the clockwise exploration of $\mathfrak{T}^\star$) of $\t_{\tilde\omega_{i_j}}$ with label $\ve$. Then
we have $D_\star(a,\rho_j)\leq D_\star(a_{(\ve)},\rho_j)+ D_\star(a,a_{(\ve)})$, and on one
hand $D_\star(a_{(\ve)},\rho_j)\leq H_j$, on the other hand, $D_\star(a,a_{(\ve)})= \ell_a -\ve$ (because
labels ``between'' $a_{(\ve)}$ and $a$ remain greater than $\ve$). 

In the case $a=\rho_{j-1}$, \eqref{bd-tec-lem2} gives $D_\star(\rho_{j-1},\rho_j)\leq H_j + h_{j-1}$. 
Using the triangle inequality, it follows that, for every $a,b\in  [|0,\rho_{k}|]_\star$ with $\ell_a\wedge \ell_b>0$,
$$D_\star(a,b)\leq \sum_{j=k}^\infty(H_j + h_j) + 2 \sup\{\ell_c:c\in [|0,\rho_{k}|]_\star\}.$$
The right-hand side can be made arbitrarily small by choosing $k$ large. This proves that $D_\star(a,b)$
tends to $0$ when $a,b\to 0$ in  $\mathfrak{T}^\star\backslash\partial\mathfrak{T}^\star$. This completes the
proof of the first assertion of the proposition. 

We note that the continuous extension of $(a,b)\mapsto D_\star(a,b)$
also satisfies the triangle inequality and the bound $D_\star(a,b)\geq |\ell_a-\ell_b|$.

Let us turn to the second assertion. Trivially the property $D^\circ_\star(a,b)=0$
for $a,b\in\mathfrak{T}^\star\backslash\partial\mathfrak{T}^\star$ implies $D_\star(a,b)=0$. 
Then suppose that $a,b\in \partial\mathfrak{T}^\star$, and $a=\ee^\star_s$, $b=\ee^\star_t$,
with $0\leq s< t\leq \Sigma^\star$, and $\ell_{\ee^\star_r}>0$ for every $r\in(s,t)$. Take $u=(s+t)/2$
and for every $\ve\in(0,\ell_{\ee^\star_u})$, define 
$$s_\ve:=\sup\{r\in[s,u]:\ell_{\ee^\star_r}=\ve\}\,,\quad t_\ve:=\inf\{r\in[u,t]:\ell_{\ee^\star_r}=\ve\}.$$
Then, $D_\star(\ee^\star_{s_\ve},\ee^\star_{t_\ve})=0$ and by letting $\ve\to 0$ we get 
$D_\star(a,b)=0$.

Conversely, if $D_\star(a,b)=0$, with $\ell_a=\ell_b>0$, this implies a fortiori that $D(a,b)=0$
and, from results recalled in Section \ref{Bet-Mie}, this can only occur if
$D^\circ(a,b)=0$. In particular labels are greater than or equal to $\ell_a$ on $[|a,b|]$ (or on $[|b,a|]$) and it readily follows
that we have also $D^\circ_\star(a,b)=0$. 

Finally, suppose that $D_\star(a,b)=0$ for distinct points $a,b\in\partial\mathfrak{T}^\star$. 
Let us write $a=\ee^\star_s$, $b=\ee^\star_t$,
with $0\leq s< t\leq \Sigma^\star$, and exclude the case $\{s,t\}=\{0,\Sigma^\star\}$. Note that $[|a,b|]_\star=\{\ee^\star_r:s\leq r\leq t\}$.
Argue by contradiction and suppose that $\ell_c$ vanishes for some
$c\in[|a,b|]_\star\backslash\{a,b\}$. Then  there must also exist $c'\in[|a,b|]$ such that $\ell_{c'}<0$
(because $0$ cannot be a local minimum for the function $s\mapsto \wh W_s(\omega_i)$, for any $i\in I$, see
the end of Section \ref{sna-mea}).
Since it is also clear that the minimum of labels on $[|b,a|]$ is negative, it
follows that $D^\circ(a,b)>0$, which in turn implies that $D(a,b)>0$, and a fortiori $D_\star(a,b)>0$ (the bound 
$D_\star(a,b)\geq D(a,b)$ remains valid for any $a,b\in \mathfrak{T}^\star$ by continuity). This contradicts 
our initial assumption $D_\star(a,b)=0$.

The preceding argument does not work for
$a=0$ and $b=\xi$ because $D^\circ(0,\xi)=0$, but we can argue as follows.
Suppose that $D_\star(0,\xi)=0$. For every integer $n>2/\xi$, we choose $r_n\in(0,1/n)$ small enough so that
$D_\star(0,r_n)<1/n$ and $D_\star(\xi-r_n,\xi)<1/n$. Then $D_\star(r_n,\xi-r_n)<2/n$
by the triangle inequality. It follows that we can find $a^{(n)}_0=r_n,a^{(n)}_1,\ldots,a^{(n)}_{p_n-1},a^{(n)}_{p_n}=\xi-r_n$ in $\mathfrak{T}^\star
\backslash \partial\mathfrak{T}^\star$
such that $\sum_{j=1}^{p_n} D_\star^\circ(a^{(n)}_{j-1},a^{(n)}_j)<2/n$. Next fix any $i_0\in I$ such that labels on $\t_{(\tilde\omega_{i_0})}$ 
vanish. For $n$ large enough so that $r_n<t_{i_0}<\xi-r_n$, at least one of the 
points $a^{(n)}_j$, $1\leq j\leq p_{n-1}$, say $a^{(n)}_{j_n}$, must belong to $\t_{(\tilde\omega_{i_0})}$ (otherwise we would have 
$D_\star^\circ(a^{(n)}_{j-1},a^{(n)}_j)=\infty$ for some $j$). By extracting a convergent sequence from
the sequence $(a^{(n)}_{j_n})$, we get a point $a^{(\infty)}$ of $\t_{(\tilde\omega_{i_0})}$
such that $D_\star(0,a^{(\infty)})=0$. By the cases treated previously, this is impossible, and
this contradiction completes the proof. \endproof

\bigskip
 We then consider the quotient space $\U:=\mathfrak{T}^\star/\{D_\star=0\}$, and the canonical projection $\Pi_\star:\mathfrak{T}^\star\la \U$. In contrast with the 
construction of $\D$, we observe that $\Pi_\star(0)\not=\Pi_\star(\xi)$. 
The function  $(a,b)\mapsto D_\star(a,b)$ induces 
a metric on $\U$, which we still denote by $D_\star$, and the metric space 
$(\U,D_\star)$ is equipped with the pushforward of the 
volume measure on $\mathfrak{T}^\star$ under $\Pi_\star$. This measure will be 
denoted by $\mathbf{V}_\star$. We also write $\partial_0\U=\Pi_\star([0,\xi])$ and
$\partial_1\U=\Pi_\star(\partial\mathfrak{T}^\star)$. 
We note that labels $\ell_x$ make sense for $x\in\U$ (because 
$D_\star(a,b)=0$ implies $\ell_a=\ell_b$). Furthermore, we have $D_\star(x,\partial_1\U)=\ell_x$
for every $x\in\U$ (here and below, we use the notation $D_\star(x,A):=\inf\{D_\star(x,y):y\in A\}$). Indeed, the bound $D_\star(x,\partial_1\U)\geq\ell_x$ is immediate
since $D_\star(x,y)\geq |\ell_x-\ell_y|$ for every $y\in \U$. Conversely, if $x=\Pi_\star(\ee^\star_r)$
and $s=\inf\{t\geq r:\ell_{\ee^\star_t}=0\}$, it is easy to verify that $D_\star(x,\ee^\star_s)=\ell_x$. 

Our goal is to verify that the 
random measure metric space $(\U,D_\star,\mathbf{V}_\star)$ equipped with an
appropriate standard boundary curve
is a curve-decorated free Brownian disk with a random perimeter. The boundary of this Brownian disk
will be $\partial_0\U \cup \partial_1\U$. 
To this end, we will first identify $(\U,D_\star,\mathbf{V}_\star)$ with another space
constructed directly from the Brownian disk $\D$. 

\subsection{Identification of $\U$}
\label{sec:identi}

We consider the free pointed Brownian disk $\D$ with perimeter $\xi$
constructed in Section \ref{Bet-Mie}. Recall the notation $\Pi$ for the canonical projection from $\mathfrak{T}$ onto $\D$, and $\xx_*$ for the distinguished point of $\D$.
To simplify notation, we will also write $\xx_0=\Pi(0)=\Pi(\xi)$ for the unique point 
of $\partial \D$ at minimal distance from $\xx_*$, and $r_0=D(\xx_*,\xx_0)=-\ell_{\xx_*}$. 
Let $\mathcal{B}_\D(\xx_*,r_0)$ stand for the closed ball of radius $r_0$ centered at $\xx_*$ in $\D$. 
The hull $H$ of radius $r_0$
centered at $\xx_*$ is the closed subset of $\D$ defined by saying that $\D\backslash H$ is the connected component
of $\D\backslash\mathcal{B}_\D(\xx_*,r_0)$ that
contains $\partial \D\backslash\{\xx_0\}$. We also let $U$ be the closure of $\D\backslash H$,
and write $\mathrm{Int}(U)=\D\backslash H$ for the topological interior of $U$ and $\partial U=U\backslash \mathrm{Int}(U)=\partial H$ for its topological boundary. Notice that
$\mathrm{Int}(U)$ contains $\partial\D\backslash\{\xx_0\}$. We also observe
that $D(\xx_*,x)=r_0$ for every $x\in\partial H$. Since  $\D$ is a geodesic space, it follows that, for every $x\in\mathrm{Int}(U)$,
$D(x,\partial H)=D(x,\xx_*)-r_0=\ell_x$.

Let $d_\infty$ stand for the
intrinsic distance on $\mathrm{Int}(U)$. For every $x,y\in \mathrm{Int}(U)$, $d_\infty(x,y)$ is the infimum 
of lengths of paths connecting $x$ to $y$ that stay in $\mathrm{Int}(U)$ ---
here lengths are of course evaluated with respect to the metric $D$ on $\D$. The fact that
$\D$ is a geodesic space easily implies that $d_\infty(x,y)<\infty$ for every $x,y\in \mathrm{Int}(U)$.
We let $(U',d_\infty)$ stand for the completion of the metric space 
$(\mathrm{Int}(U),d_\infty)$. 
We write $\mathbf{V}_U$ for the 
restriction of the volume measure $\mathbf{V}$ to $\mathrm{Int}(U)$, which 
may also be viewed as a measure on $U'$ since $\mathrm{Int}(U)$ is an open subset of $U'$. 
Finally, we recall the abuse of notation that consists in viewing $\mathfrak{T}^\star$
as a subset of $\mathfrak{T}$.

\begin{proposition}
\label{ident-lim-space}
The measure metric spaces $(U',d_\infty,\mathbf{V}_U)$ and $(\U,D_\star,\mathbf{V}_\star)$
are almost surely equal. More precisely, the following property holds a.s.: there is a one-to-one mapping $\Psi$ from
$\mathrm{Int}(U)$ into $\U$ such that $\Psi$ maps $\Pi(a)$ to $\Pi_\star(a)$, for every
$a\in \mathfrak{T}^\star\backslash \partial\mathfrak{T}^\star$, and the mapping $\Psi$
extends to a measure-preserving isometry from $(U',d_\infty,\mathbf{V}_U)$ onto $(\U,D_\star,\mathbf{V}_\star)$.
\end{proposition}

\proof Let us first
identify the hull $H$ in terms of our construction of $\D$. To this end,
for every $x\in \D$, set $m_x=\ell_x$ if $x\in[0,\xi]$, and, if $x=\Pi(a)$ where $a\in\t_{(\omega_i)}$ 
for some $i\in I$, let $m_x$ be the minimal label along the 
ancestor line of $a$ in $\t_{(\omega_i)}$. Notice that this definition does not
depend on the choice of $a$ such that $x=\Pi(a)$. Then an easy extension
of the so-called cactus bound (cf. Proposition 3.1 in \cite{Geodesics}) shows
that any continuous curve from $x$ to $\partial \D$ has to come within distance $r_0+m_x$
from $\xx_*$. On the other hand, if $x=\Pi(a)$ with $a\in\t_{(\omega_i)}$, the image of the ancestral line of $a$ under $\Pi$ provides 
a path from $x$ to $\partial\D\backslash\{\xx_0\}$ whose minimal distance from $\xx_*$ is $r_0+m_x$. It follows that
 $x\in \D\backslash H$ if and only if $m_x> 0$, and consequently
$$\mathrm{Int}(U)=\D\backslash H=\Pi(\mathfrak{T}^\star\backslash \partial\mathfrak{T}^\star),\quad
H=\Pi(\mathfrak{T}\backslash \mathfrak{T}^\star) \cup \Pi(\partial \mathfrak{T}^\star),\quad
\partial H=\partial U=\Pi(\partial\mathfrak{T}_\star)\,.
$$

The next step is to prove that, if $x=\Pi(a)\in \mathrm{Int}(U)$
and $y=\Pi(b)\in \mathrm{Int}(U)$, where $a,b\in \mathfrak{T}^\star\backslash \partial\mathfrak{T}^\star$,
we have 
\begin{equation}
\label{ident-dist}
d_\infty(x,y)=D_\star(a,b).
\end{equation}
 The upper bound $d_\infty(x,y)\leq D_\star(a,b)$
is easy because each quantity of the form
$$\sum_{i=1}^p D^\circ_\star(a_{i-1},a_i),$$
where $a_0=a$, $a_p=b$ and $a_1,\ldots,a_{p-1}\in\mathfrak{T}^\star\backslash \partial\mathfrak{T}^\star$  are such that $D^\circ_\star(a_{i-1},a_i)<\infty$ for $1\leq i\leq p$,
is the length of a curve from $x$ to $y$ that stays in $\mathrm{Int}(U)$ (use the simple geodesics defined at the end
of Section \ref{Bet-Mie} and note that the condition
$D^\circ_\star(a_{i-1},a_i)<\infty$ precisely prevents the curve
from hitting the set of points with zero label). 

Let us now justify the reverse bound $d_\infty(x,y)\geq D_\star(a,b)$. It is enough to verify that,
if $\gamma=(\gamma(t))_{0\leq t\leq 1}$ is a curve connecting $x$ to $y$ and
staying in $\mathrm{Int}(U)$, the length of $\gamma$ is bounded below by $D_\star(a,b)$.
Since $\gamma$ stays in $\mathrm{Int}(U)$, 
a compactness argument shows that $\mathrm{dist}_D(\gamma,\partial U)>0$
(here we write $\mathrm{dist}_D(\gamma,\partial U)$ for the minimal $D$-distance between 
a point of the curve $\gamma$ and $\partial U$). So we can fix an integer $n\geq 1$
such that, for every $i\in\{1,\ldots,n\}$,
we have
\begin{equation}
\label{ass-reverse}
D\Big(\gamma(\frac{i-1}{n}),\gamma(\frac{i}{n})\Big) <\frac{1}{2} \mathrm{dist}_D(\gamma,\partial U).
\end{equation}
Since the length of $\gamma$ is bounded below by 
\begin{equation}
\label{keysum}
\sum_{i=1}^n D\Big(\gamma(\frac{i-1}{n}),\gamma(\frac{i}{n})\Big)
\end{equation}
it is enough to verify that the latter sum is bounded below by $D_\star(a,b)$. Fix $\ve\in(0,\frac{1}{4}\mathrm{dist}_D(\gamma,\partial U))$. By the definition of the 
metric $D$, for every $i\in\{1,\ldots,n\}$, we can find $a^i_0,a^i_1,\ldots, a^i_{p_i}\in \mathfrak{T}$
with $\Pi(a^i_0)=\gamma(\frac{i-1}{n})$ and $\Pi(a^i_{p_i})=\gamma(\frac{i}{n})$ such that
\begin{equation}
\label{ass-reverse2}
D\Big(\gamma(\frac{i-1}{n}),\gamma(\frac{i}{n})\Big) 
\geq \sum_{j=1}^{p_i} D^\circ(a^i_{j-1},a^i_j) - \frac{\ve}{n}.
\end{equation}
This implies in particular that, for every $k\in\{1,\ldots,p_i\}$,
$$D(a^i_0,a^i_k)\leq \sum_{j=1}^{k} D^\circ(a^i_{j-1},a^i_j)\leq D\Big(\gamma(\frac{i-1}{n}),\gamma(\frac{i}{n})\Big) + \frac{\ve}{n} <\mathrm{dist}_D(\gamma,\partial U).$$
Since $\Pi(a^i_0)$ belongs to the range of $\gamma$, it follows that $\Pi(a^i_k)\in \mathrm{Int}(U)$. Thus all points $a^i_k$ must belong to
$\mathfrak{T}^\star\backslash \partial \mathfrak{T}^\star$ (recall that $\Pi(\mathfrak{T}\backslash\mathfrak{T}^\star) \cup \Pi(\partial \mathfrak{T}^\star)=H$). We now claim that 
\begin{equation}
\label{claim-ident}
D^\circ(a^i_{j-1},a^i_j)=D^\circ_\star(a^i_{j-1},a^i_j),
\end{equation}
for every $i\in\{1,\ldots,n\}$ and $j\in\{1,\ldots,p_i\}$. If the claim is proved, we obtain that
$$\sum_{i=1}^n D\Big(\gamma(\frac{i-1}{n}),\gamma(\frac{i}{n})\Big)
\geq \sum_{i=1}^n \sum_{j=1}^{p_i} D^\circ_\star(a^i_{j-1},a^i_j) - \ve \geq D_\star(a,b) -\ve,$$
by the very definition of $D_\star(a,b)$. Since $\ve$ was arbitrary, this shows that 
the sum \eqref{keysum} is bounded below by $D_\star(a,b)$, as desired.

Let us prove our claim \eqref{claim-ident}. From our definitions, it is 
enough to verify that 
\begin{equation}
\label{claim-ident2}
\max\Big(\min_{c\in[|a^i_{j-1},a^i_j|]}\ell_c,\min_{c\in[|a^i_{j},a^i_{j-1}|]}\ell_c\Big)>0
\end{equation}
(note that, if for instance $\min_{c\in[|a^i_{j-1},a^i_j|]}\ell_c>0$, this implies that 
$[|a^i_{j-1},a^i_j|]=[|a^i_{j-1},a^i_j|]_\star$).
Fix $i\in\{1,\ldots,n\}$ and $j\in\{1,\ldots,p_i\}$. Using \eqref{ass-reverse}, \eqref{ass-reverse2} and
 the 
bound $D^\circ(u,v)\geq |\ell_u-\ell_v|$, we have
\begin{align*}
\mathrm{dist}_D(\gamma,\partial U)&\geq \sum_{k=1}^{p_i} D^\circ(a^i_{k-1},a^i_k)\\
&\geq |\ell_{\gamma(\frac{i-1}{n})}-\ell_{a^i_{j-1}}| 
+\ell_{a^i_{j-1}} +\ell_{a^i_{j}}-2\max\Big(\min_{c\in[|a^i_{j-1},a^i_j|]}\ell_c,\min_{c\in[|a^i_{j},a^i_{j-1}|]}\ell_c\Big) +|\ell_{\gamma(\frac{i}{n})}
-\ell_{a^i_{j}}|\\
&\geq \ell_{\gamma(\frac{i-1}{n})} + \ell_{\gamma(\frac{i}{n})} -2\max\Big(\min_{c\in[|a^i_{j-1},a^i_j|]}\ell_c,\min_{c\in[|a^i_{j},a^i_{j-1}|]}\ell_c\Big)\\
&\geq 2\Bigg(\mathrm{dist}_D(\gamma,\partial U) - \max\Big(\min_{c\in[|a^i_{j-1},a^i_j|]}\ell_c,\min_{c\in[|a^i_{j},a^i_{j-1}|]}\ell_c\Big)\Bigg),
\end{align*}
thanks to the equality $D(x,\partial U)=\ell_x$ for $x\in\mathrm{Int}(U)$.
Clearly this implies that \eqref{claim-ident2} holds, completing the proof of \eqref{ident-dist}.

Recall that $\mathrm{Int}(U)=\Pi(\mathfrak{T}^\star\backslash \partial\mathfrak{T}^\star)$,
and also recall the notation $\partial_1\U=\Pi_\star(\partial\mathfrak{T}^\star)$. By
Proposition \ref{def-lim-space}, two points 
$a$ and $b$ of $\mathfrak{T}^\star\backslash \partial\mathfrak{T}^\star$ are identified in the quotient space 
defining $\U$ if and only if they are identified in the quotient space $\D$. 
This shows that $\U\backslash\partial_1\U=\Pi_\star(\mathfrak{T}^\star\backslash \partial\mathfrak{T}^\star)$ is
identified as a set to $\mathrm{Int}(U)$, and \eqref{ident-dist} entails that
this identification is an isometry when $\mathrm{Int}(U)$ is equipped with $d_\infty$ and
$\U\backslash\partial_1\U$ is equipped with $D_\star$. Since $\U\backslash\partial_1\U$ is dense in the
compact set $\U$, it follows that the completion $U'$ can
also be identified isometrically to $\U$. Furthermore, it is clear that
the volume measure $\mathbf{V}_U$ corresponds to $\mathbf{V}_\star$ in this identification.
This completes the proof of the proposition. \endproof

In the identification of $U'$ with $\U$, the ``boundary'' $\partial U':=U'\backslash \mathrm{Int}(U)$
is mapped bijectively onto $\partial_1\U$. Note that every point of $\partial U'=\partial_1\U$ 
has to correspond to a point of $\partial U$. This correspondence is one-to-one except that the two points
$\Pi_\star(0)$ and $\Pi_\star(\xi)$ of $\partial_1\U$ correspond to the same point $\xx_0$ of $\partial U$.

We now record two technical properties that will be useful in the proofs
of the next section. 

\begin{proposition}
\label{techni}
The following properties hold a.s.

\noindent{\rm(i)} For every $\delta>0$, let $U_{(\delta)}$ be the set of all 
$x\in \D$ such that there is a continuous path from $x$ to $\partial\D$ 
that stays at distance at least $r_0-\delta$ from $\xx_*$. Then, for every $\ve>0$, there exists 
$\delta>0$ such that $U_{(\delta)}\subset \{x\in\D:D(x,U)<\ve\}$. 

\noindent{\rm(ii)} For every $\ve>0$, there exists $\delta>0$ such that, if
$x\in \D$ and $D(x,H)\geq\ve$, there is a continuous path from $x$ to
$\partial \D$ that stays at distance at least $r_0+\delta$ from $\xx_*$.
\end{proposition}

\proof (i) Recall the notation $m_x$ introduced in the proof
of Proposition \ref{ident-lim-space} for $x\in \D$. As it was explained in this proof,
any continuous path from $x$ to $\partial \D$ has to come within distance $r_0+m_x$
from $\xx_*$.
So the proof boils down to verifying that, given $\ve>0$, we can find
$\delta>0$ such that 
$$\{x\in \D:m_x\geq -\delta\}\subset \{x\in\D:D(x,U)<\ve\}.$$
Note that the mapping $x\mapsto m_x$ is continuous. Moreover, the property $m_x=0$
may hold only if $x\in\partial U$. To see this,  first note that, a.s. for every $x\in\D$ of the form $x=\Pi(a)$ with $a\in \t_{(\omega_i)}$, the property $m_x=0$ implies that $\ell_x=0$
and that labels on $\llbracket \rho_{(\omega_i)},a\rrbracket\backslash\{a\}$ are positive. Indeed, if this were not the case, this would mean that
the Brownian path describing the labels along the ancestral line $\llbracket \rho_{(\omega_i)},a\rrbracket$ has a local minimum
equal to $0$ at a point distinct from $a$, which does not occur for any $a\in\mathfrak{T}$, a.s. It follows that $a\in\partial \mathfrak{T}_\star$ and $x\in\partial U$. 

The desired result then follows from a
compactness argument, since the intersection of compact sets
$$\bigcap_{\delta>0} \Big(\{x\in \D:m_x\geq- \delta\}\backslash  \{x\in\D:D(x,U)<\ve\}\Big)$$
is empty by the preceding considerations and the fact that $m_x>0$ implies $x\in \mathrm{Int}(U)$. 

\smallskip
\noindent(ii) For every $x\in \D\backslash H$, let $\varphi(x)$ be the supremum of 
all $\delta>0$ such that there is a continuous path from $x$ to $\partial\D$
that stays at distance at least $r_0+\delta$ from $\xx_*$. Then $\varphi(x)>0$
for every $x\in \D\backslash H$ and the mapping $x\mapsto \varphi(x)$ is continuous.
It follows that $\inf\{\varphi(x):x\in \D, D(x,H)\geq \ve\}>0$, giving the desired result. \endproof

\section{Passage to the limit}
\label{pass-lim}

\subsection{Preliminaries}
\label{preli}

For every integer $L\geq1$, let $\t_L$ be a Boltzmann distributed rooted and pointed type I triangulation with 
a simple boundary of size $n_L\geq 1$. We assume throughout this section that
$n_L/L \la \xi>0$ as $L\to\infty$. We use a similar notation as in Section \ref{conve-tria}. In particular, 
$\nu_L$ is the counting measure on the set 
$V_i(\t_L)$ scaled by the multiplicative factor $\frac{3}{4}L^{-2}$, $\Theta_L=(\Theta_L(0),\Theta_L(1),\ldots,\Theta_L(n_L))$
is the boundary path of $\t_L$, and $\wh\Theta_L(t)=\Theta_L(\lfloor Lt\rfloor)$ for $0\leq t\leq n_L/L$, and $v_*^{(L)}$ is the distinguished vertex of $\t_L$. From (a trivial extension of) \eqref{convGHPU-pointed}, we have
\begin{equation}
\label{convGHPU}
(V(\t_L),\sqrt{3/2}\,L^{-1/2}\dg,\nu_L,\wh\Theta_L,v_*^{(L)})\build{\la}_{L\to\infty}^{(d)} (\D,D,\mathbf{V},\Gamma,\xx_*),
\end{equation}
where the convergence holds in distribution in $(\M^{GHPU\bullet},\dd_{GHPU\bullet})$, and the
limit $(\D,D,\mathbf{V},\Gamma,\xx_*)$ is the curve-decorated free pointed Brownian disk with perimeter $\xi$ 
as defined at the end of Section \ref{Bet-Mie}. Recall that the range of $\Gamma$ is the boundary $\partial\D$.

Assuming that $\dg(v_*^{(L)},\partial\t_L)\geq 2$, we can make sense of the hull of radius $1$ centered at $v_*^{(L)}$, which we denote by $H_1(v_*^{(L)})$: We first 
define the ball $B_1(v_*^{(L)})$ as the union of all faces of $\t_L$ that are incident to $v_*^{(L)}$, and the 
hull $H_1(v_*^{(L)})$ is obtained by adding to the ball $B_1(v_*^{(L)})$ all faces of $\t_L$ that
are disconnected from $\partial\t_L$ by the ball. This hull 
can be viewed as a triangulation with a boundary, such that every boundary vertex is 
at distance $1$ from $v_*^{(L)}$.

Note that the probability of the event $\{\dg(v_*^{(L)},\partial\t_L)\geq 2\}$ tends to $1$
as $L\to\infty$. From now on, we will argue conditionally on this event. Then the complement of the hull $H_1(v_*^{(L)})$ in $\t_L$  is a triangulation $\t'_L$ with two boundaries. The first boundary 
is $\partial \t_L$ and the second one is the boundary of the hull $H_1(v_*^{(L)})$ --- we may choose
the root on the second boundary uniformly at random, independently of $\t_L$.
Furthermore, conditionally on the event that the boundary size of the hull $H_1(v_*^{(L)})$ is equal to $p\geq 1$,
the triangulation $\t'_L$ is Boltzmann distributed on $\T^{2}(n_L,p)$, so that we can apply the
results of Section \ref{discrete}. In view of this application, we also observe that the 
boundary size of $H_1(v_*^{(L)})$ (that is, the size of the second boundary of 
$\t'_L$) remains bounded in probability when $L\to\infty$. This follows from an
easy counting argument since this boundary size is bounded above by the degree of $v_*^{(L)}$
in $\t_L$. 

We can then run the peeling algorithm of $\t'_L$ as described in Section \ref{discrete}.
With a slight abuse of terminology, we will consider that the revealed region at step $n\leq\zeta_L$, as defined in
Section \ref{peeling-pro}, also includes the hull $H_1(v_*^{(L)})$, and thus can be viewed 
as a triangulation with a (simple) boundary of size $P^L_n$.
We consider
the ``peeling by layers'' algorithm, which involves particular choices of the revealed edge at
each step (see e.g.~\cite{CLG} for more details). The only property of this algorithm that we will use
is the fact that, for every step $n=0,1,\ldots,\zeta_L$, there exists an integer $r_n\geq 1$ such 
that every vertex of the boundary of the revealed version at step $n$ is at distance $r_n$
or $r_n+1$ from $v_*^{(L)}$. 
We let $r_0^{(L)}$ be $\sqrt{3/2}\,L^{-1/2}$ times the maximal graph distance between  $v_*^{(L)}$ and a point of the boundary of the 
revealed region at time $\zeta_L$. So a path from any point of the revealed region at time $\zeta_L$ to the boundary $\partial\t_L$ must come within distance  
$\sqrt{2/3}\,L^{1/2}\,r_0^{(L)}$ from the distinguished point $v_*^{(L)}$. On the other hand, by the properties of the 
peeling by layers algorithm, 
we have $\dg(v_*^{(L)},\partial\t_L)=\sqrt{2/3}\,L^{1/2}\,r_0^{(L)}$.

\begin{figure}[!h]
 \begin{center}
    \includegraphics[height=6cm,width=6cm]{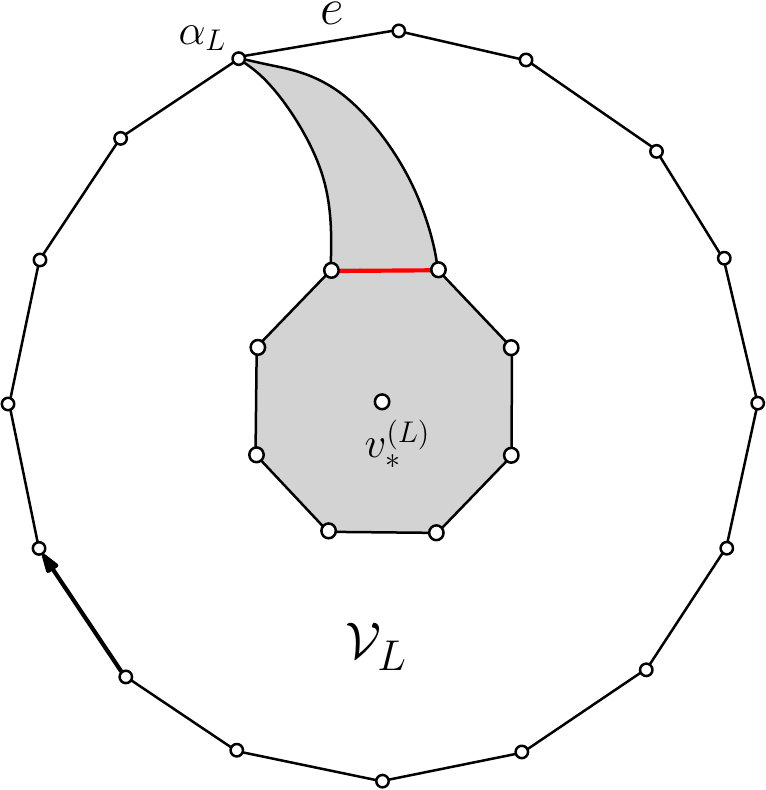} \hspace{1cm}\includegraphics[height=6cm,width=6cm]{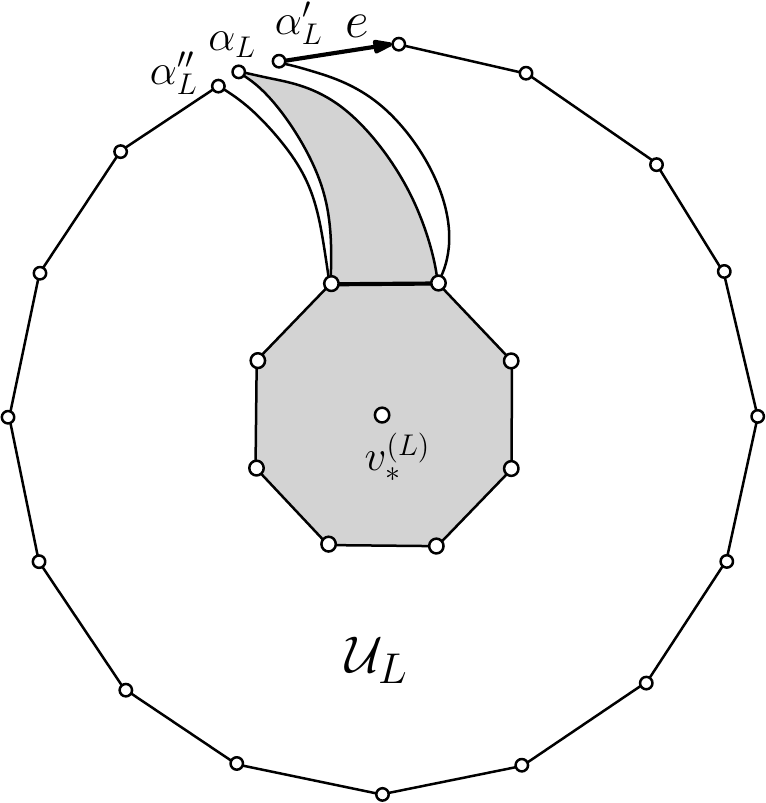} 
 \end{center}
    \caption{\label{fig:unglued} Illustration of the triangulation $\mathcal{T}_L$ before and after  ungluing the vertex $\alpha_L$. We represent the revealed and unrevealed regions  
    at time $\zeta_L$ in gray and in white respectively.}
  \end{figure}

The collection of all faces of $\t_L$ that do not belong to the revealed region at time $\zeta_L$
will be called the unrevealed 
region. It will be convenient to write $\v_L$ for the set of all vertices of the unrevealed region, and $\W_L$
for the set of all vertices of the revealed region (at time $\zeta_L$). The unrevealed region is rooted at the root of $\t_L$ and then corresponds to a triangulation with a simple boundary 
where two boundary vertices, say $\alpha'_L$ and $\alpha''_L$, have been glued to give the unique vertex $\alpha_L$ of $\partial\t_L$ 
that belongs to the revealed region at step $\zeta_L$ (see Figure \ref{fig:unglued}).  If these two vertices are ``unglued'', we get 
a triangulation $\u_L$ with a simple boundary $\partial\u_L$ of size $n_L+Z_L$, where we know 
from Proposition \ref{asymptotics-Z} that $Z_L/n_L$ converges in distribution to a random variable denoted by $\Lambda$. Moreover, if
$e$ stands for the unique oriented edge of $\partial\t_L$ whose tail
is $\alpha_L$, and such that the external face of $\t_L$ lies to the left of $e$, we may root $\u_L$ at the
edge corresponding to $e$, whose origin is $\alpha'_L$  (see again Figure \ref{fig:unglued}). From the discussion 
in Section \ref{peeling-proba}, one easily checks that, conditionally on $\{Z_L=k\}$, the triangulation $\u_L$
is Boltzmann distributed in $\T^1(n_L+k)$. 

Since, conditionally on $Z_L$, the map $\u_L$ is a Boltzmann distributed rooted triangulation with perimeter $n_L+Z_L$,
and $L^{-1}Z_L$ converges in distribution to $\y:=\xi\,\Lambda $, we can also 
apply the convergence \eqref{convGHPU-unpointed} to the 
triangulations $\u_L$. Write $(\Theta'_L(j))_{0\leq j\leq n_L+Z_L}$ for the boundary path of $\u_L$, and $\wh \Theta'_L(t)=\Theta'_L(\lfloor Lt\rfloor)$ 
for $0\leq t\leq (n_L+Z_L)/L$. Also let $d_L$ stand for the graph distance on $V(\u_L)$ scaled by the factor
$\sqrt{3/2}L^{-1/2}$, and let $\theta_L$ be the counting measure on $V_i(\u_L)$ scaled by $\frac{3}{4}L^{-2}$. Then 
\begin{equation}
\label{convGHPU2}
(V(\u_L),d_L,\theta_L,\wh \Theta'_L) \build{\la}_{L\to\infty}^{(d)} \Big(\wt \D, \wt D, \wt{\mathbf{V}},\wt\Gamma\Big),\end{equation}
where $(\wt \D, \wt D, \wt{\mathbf{V}},\wt\Gamma)$ is a curve-decorated free (non-pointed) Brownian disk with  perimeter $\xi+\y$,
and the convergence in distribution holds in $(\M^{GHPU},\dd_{GHPU})$.

\subsection{Paths avoiding part of the boundary}

The boundary $\partial\u_L$ of $\u_L$ coincides with the set $\{\Theta'_L(j):0\leq j\leq n_L+Z_L\}$. We will consider the subset 
$\partial_1\u_L$ of this boundary corresponding
to the boundary of the revealed region at time $\zeta_L$. Precisely, we set $\partial_1\u_L:=\{\Theta'_L(j):n_L\leq j\leq n_L+Z_L\}$
and $\partial_0\u_L:=\{\Theta'_L(j):0\leq j\leq n_L\}$. We also define $\partial_1\wt \D:=\{\wt\Gamma(t):\xi\leq t\leq \xi+\mathcal{Y}\}$. 

By the Skorokhod representation theorem, we may and will assume in this section  that the convergence \eqref{convGHPU2} holds a.s.,
and moreover $Z_L/L\la \mathcal{Y}$ a.s. 

\begin{lemma}
\label{avoiding1}
Let $\delta>0$ and $\eta>0$. Then a.s.
there exists $\ve_0>0$
such that, for every large enough $L$, for every $x,y\in V(\u_L)$
with $d_L(x,\partial_1\u_L)\geq \delta$ and $d_L(y,\partial_1\u_L)\geq \delta$,
there is a path from $x$ to $y$ that stays at $d_L$-distance at least $\ve_0$
from $\partial_1 \u_L$, and whose $d_L$-length is bounded above
by $d_L(x,y)+\eta$. 
\end{lemma}

\proof Let us fix $\omega$ in the underlying probability space such that the (a.s.)
convergence \eqref{convGHPU2} holds and $Z_L/L\la \mathcal{Y}$. 
By \cite[Proposition 1.5]{GM0} (or Proposition \ref{GHPU-criterion}), we may assume that all metric spaces $(V(\u_L),d_L)$ and $(\wt \D,\wt D)$ are embedded isometrically simultaneously 
in the same compact metric space $(E,\Delta)$, in such a way that we have
$V(\u_L)\la \wt\D$ for the Hausdorff distance and $(\wh\Theta'_L(t\wedge (n_L+Z_L)/L))_{t\geq 0}$
converges uniformly to $(\wt\Gamma(t\wedge (\xi+\mathcal{Y})))_{t\geq 0}$. In particular, $\partial_1\u_L\la \partial_1\wt\D$ in the sense of the
Hausdorff distance between compact subsets of $(E,\Delta)$.

We may assume that $\eta<\delta$. We argue by contradiction.
If the statement of the lemma does not hold, we can find a sequence $\ve_k\downarrow 0$
and a sequence $L_k\uparrow\infty$ such that, for every $k$, there exist
$x_k,y_k\in V(\u_{L_k})$
with $d_{L_k}(x_k,\partial_1\u_{L_k})\geq \delta$ and $d_{L_k}(y_k,\partial_1\u_{L_k})\geq \delta$,
such that any path from $x_k$ to $y_k$  that stays at distance at least $\ve_k$
from $\partial_1 \u_{L_k}$ has $d_{L_k}$-length at least $d_{L_k}(x_k,y_k)+\eta$. 

By compactness, we may assume that $x_k\la x_\infty$ and $y_k\la x_\infty$, where
$x_\infty,y_\infty\in \wt\D$, and $\wt D(x_\infty,\partial_1\wt\D)\geq \delta$, $\wt D(y_\infty,\partial_1\wt\D)\geq \delta$. We take $k$ large enough so
that $\Delta(x_k,x_\infty)<\eta/16$ and $\Delta(y_k,y_\infty)<\eta/16$ and then we have also
$d_{L_k}(x_k,y_k)=\Delta(x_k,y_k)\geq \Delta(x_\infty,y_\infty)-\eta/8=\wt D(x_\infty,y_\infty)-\eta/8$. 
Next, by Lemma 19 in \cite{BM}, we can find a point
$x'_\infty\in\wt\D$ with $\wt D(x_\infty,x'_\infty)<\eta/16$ and a point $y'_\infty\in\wt\D$ with
$\wt D(y_\infty,y'_\infty)<\eta/16$, such that there is a $\wt D$-geodesic from $x'_\infty$
to $y'_\infty$ that does not hit $\partial\wt\D$, and thus stays at distance at least $\alpha>0$
from $\partial\wt\D$, for some $\alpha>0$. We can assume that $\alpha<\eta/8$, and we write $(\gamma(t))_{0\leq t\leq \tilde D(x'_\infty,y'_\infty)}$
for the latter geodesic. To simplify notation, we set $d_*=\wt D(x'_\infty,y'_\infty)$ and we note that
$d_{L_k}(x_k,y_k)\geq \wt D(x_\infty,y_\infty)-\eta/8\geq d_*-\eta/2$.

Then pick an integer $N\geq 4$ large enough so that $\frac{d_*}{N}<\frac{\alpha}{4}$. 
Taking $k$ even larger if necessary, we can choose, for every $0\leq i\leq N$, a point
$x^{(k)}_i$ in $V(\u_{L_k})$ such that $\Delta(\gamma(\frac{id_*}{N}),x^{(k)}_i)<\alpha/(2N)$. 
Hence, for every $0\leq i\leq N-1$, we have
$$\Delta\big(x^{(k)}_i,x^{(k)}_{i+1}\big)\leq \Delta\big(x^{(k)}_i,\gamma(\frac{i d_*}{N})\big)+ \Delta\big(\gamma(\frac{i d_*}{N}),\gamma(\frac{(i+1)d_*}{N})\big)+\Delta\big(\gamma(\frac{(i+1)d_*}{N}),
x^{(k)}_{i+1}\big)<\frac{d_*}{N} + \frac{\alpha}{N}\leq \frac{\alpha}{2}.$$
Moreover,
$$\Delta\big(x_k,x^{(k)}_0\big)\leq \Delta\big(x_k,x_\infty\big)+\Delta\big(x_\infty,x'_\infty\big)+\Delta\big(x'_\infty, x^{(k)}_0\big)<\frac{\eta}{16}+\frac{\eta}{16} + \frac{\alpha}{2N}\leq \frac{3\eta}{16},
$$
and similarly $\Delta(y_k,x^{(k)}_N)<\frac{3\eta}{16}$. If we now concatenate a geodesic from $x_k$ to $x^{(k)}_0$
in $\u_{L_k}$
with geodesics from $x^{(k)}_i$ to $x^{(k)}_{i+1}$, for $0\leq i\leq N-1$, and finally with a geodesic from
$x^{(k)}_N$ to $y_k$, we get a path from $x_k$ to $y_k$ with length smaller that $d_*+\alpha+\frac{3\eta}{8}<d_{L_k}(x_k,y_k)+\eta$.
Furthermore, recalling that the $\Delta$-Hausdorff distance between $\partial_1\u_L$ and $\partial_1\wt\D$ tends to $0$ when $L\to\infty$, and using
the bounds
$\wt D(x_\infty,\partial_1\wt\D)\geq \delta>\eta$, $\wt D(y_\infty,\partial_1\wt\D)\geq \delta>\eta$, 
$\Delta(x_k,x_0^{(k)})<\frac{3\eta}{16}$, $\Delta(y_k,x_N^{(k)})<\frac{3\eta}{16}$,
and the fact that the geodesic $\gamma$
stays at distance at least $\alpha$
from $\partial\wt\D$, one easily verifies that the path from $x_k$ to $y_k$ that we constructed stays at distance
at least $\alpha/4$ from $\partial_1\u_{L_k}$ when 
$k$ is large. Taking $k$ such that $\ve_k<\alpha/4$, we get a contradiction with our initial assumption. 
This contradiction completes the proof. \endproof

For $\delta>0$, we set
$$\v_{L,\delta}:=\{x\in V(\u_L):d_L(x,\partial_1\u_L)\geq \delta\}.
$$
From the a.s. convergence \eqref{convGHPU2}, the latter set is not empty
for $\delta$ small.
We will view $\v_{L,\delta}$ as a bipointed metric space with distinguished points
\begin{align*}
z_{L,\delta}&:=\Theta'_L(\inf\{j\in\{0,1,\ldots,n_L\}:d_L(\Theta'_L(j),\partial_1\u_L)\geq \delta\}),\;\\
z'_{L,\delta}&:=\Theta'_L(\sup\{j\in\{0,1,\ldots,n_L\}:d_L(\Theta'_L(j),\partial_1\u_L)\geq \delta\}).
\end{align*}
Again this definition makes sense if $\delta$ is small enough, which we assume in what follows.
We also view $(V(\u_L),d_L)$ as a bipointed metric space whose distinguished 
points are $\Theta'_L(0)$ and $\Theta'_L(n_L)$. Finally, we
recall the notation $\dd_{GH\bullet\bullet}$ for the bipointed Gromov-Hausdorff distance
(Section \ref{sec:convmetric}). 

\begin{lemma}
\label{main2}
Almost surely, we have
$$\lim_{\delta\to 0}\Bigg(\limsup_{L\to\infty} \dd_{GH\bullet\bullet}\Big(\big(\v_{L,\delta},d_L,z_{L,\delta},
z'_{L,\delta}\big),\big(V(\u_L),d_L,\Theta'_L(0),\Theta'_L(n_L)\big)\Big)
\Bigg)= 0$$
where the limit holds as $\delta\to 0$, except possibly along a countable set of values of $\delta$.
\end{lemma}

\proof We may assume that all metric spaces $(V(\u_L),d_L)$ and $(\wt \D,\wt D)$ are embedded isometrically 
in $(E,\Delta)$ in the way explained at the beginning of the proof of
Lemma \ref{avoiding1}. Set
$$\wt \D_\delta:=\big\{x\in \wt \D: \wt D(x,\partial_1\wt\D)\geq \delta\big\}.$$
Then $\wt \D_\delta\not=\varnothing$ when $\delta>0$ is small enough, which we assume from now on. We view $(\wt\D_\delta,\wt D)$ as a bipointed metric space whose distinguished
points are 
$$\wt z_\delta:=\wt\Gamma\big(\inf\big\{t\in[0,\xi]: \wt D(\wt\Gamma(t),\partial_1\wt \D)\geq \delta\big\}\big),\quad
\wt z'_\delta:=\wt\Gamma\big(\sup\big\{t\in[0,\xi]: \wt D(\wt\Gamma(t),\partial_1\wt \D)\geq \delta\big\}\big)$$ 
(again, this makes
sense for $\delta>0$ small, which we assume in the following). 
We claim that
\begin{equation}
\label{main21}
\v_{L,\delta}\build{\la}_{L\to\infty}^{} \wt\D_\delta,
\end{equation}
for the $\Delta$-Hausdorff distance, provided that $\delta>0$ is not a local maximum of 
the function $x\mapsto \wt D(x,\partial_1\wt\D)$ on $\wt \D$ (there are only countably many 
such local maxima). 

Indeed, suppose that, for some $\alpha>0$, there exists a sequence $L_k\uparrow\infty$
and, for every $k$, a point $x_k\in\wt\D_\delta$ such that $\Delta(x_k,\v_{L_k,\delta})>\alpha$.
By compactness, we may assume that $x_k\la x_\infty\in \wt\D_\delta$. From the condition
we imposed on $\delta$, we can find 
$x'_\infty\in\wt\D$ such that $\wt D(x_\infty,x'_\infty)<\alpha/3$ and $\wt D(x'_\infty,\partial_1\wt \D)>\delta$. 
Since  $V(\u_{L_k})$ converges to $\wt\D$ in the sense of the
$\Delta$-Hausdorff distance, we can find a point
$x'_k\in V(\u_{L_k})$, for every $k$, in such a way that
$\Delta(x'_k,x'_\infty)\la 0$ as $k\to\infty$. This last property and the convergence of
$\partial_1\u_{L_k}$ to $\partial_1\wt\D$ ensure that,
for $k$ large enough, we have $\Delta(x'_k,\partial_1\u_{L_k})>\delta$
 and thus $x'_k\in \v_{L_k,\delta}$. Finally, writing $d_{L_k}(x_k,x'_k)\leq \Delta(x_k,x_\infty)+
 \Delta(x_\infty,x'_\infty)+ \Delta(x'_\infty,x'_k)$, we see that we have also $d_{L_k}(x_k,x'_k)<\alpha$
 when $k$ is large, which contradicts $\Delta(x_k,\v_{L_k,\delta})>\alpha$. This contradiction
 shows that, for every $\ve>0$, the set $\wt\D_\delta$ is contained in 
 $\{x\in E:\Delta(x,\v_{L,\delta})<\ve\}$ when $L$ is large enough. 
 
 A similar (easier) argument shows that, for every $\ve>0$, the set $\v_{L,\delta}$
 is contained in $\{x\in E:\Delta(x,\wt\D_\delta)<\ve\}$ when $L$ is large. This completes the
 proof of our claim \eqref{main21}. 
 
 We then observe that we have also, on the event where $\wt z_\delta$ and $ \wt z'_\delta$ are defined,
 \begin{equation}
\label{main22}
z_{L,\delta} \build{\la}_{L\to\infty}^{}  \wt z_\delta\,,\quad
z'_{L,\delta} \build{\la}_{L\to\infty}^{}  \wt z'_\delta\,,
\end{equation}
provided that $\delta$ is not a local maximum
of the function $[0,\xi+\y]\ni t\mapsto \wt D(\wt\Gamma(t),\partial_1\wt \D)$. In fact, under the last assumption,
using the uniform convergence of $\wh\Theta'_L$ to $\wt\Gamma$ and the 
fact that $\partial_1\u_L$ converges to $\partial_1\wt \D$ for the $\Delta$-Hausdorff distance,
it is easy to verify that $\inf\{t\in[0,n_L/L]: d_L(\wh\Theta'_L(t),\partial_1\u_L)\geq \delta\}$
converges to $\inf\{t\in[0,\xi]: \wt D(\wt\Gamma(t),\partial_1\wt \D)\geq \delta\}$ as $L\to\infty$,
which gives the first convergence in \eqref{main22}. A
similar argument applies to the second convergence.

From \eqref{main21} and \eqref{main22}, we conclude that $\dd_{GH\bullet\bullet}((\v_{L,\delta},d_L,z_{L,\delta},z'_{L,\delta}),(\wt\D_\delta,
\wt D,\wt z_\delta,\wt z'_\delta))$ converges to $0$ as $L\to\infty$, except possibly
for $\delta$ belonging to a countable set. Finally writing\newpage
\begin{align*}
 &\dd_{GH\bullet\bullet}((\v_{L,\delta},d_L,z_{L,\delta},z'_{L,\delta}),(V(\u_L),d_L,\Theta'_L(0),\Theta'_L(n_L)))\\
 &\qquad\leq \dd_{GH\bullet\bullet}((\v_{L,\delta},d_L,z_{L,\delta},z'_{L,\delta}),(\wt\D_\delta,
\wt D,\wt z_\delta,\wt z'_\delta))\\
&\qquad+ \dd_{GH\bullet\bullet}((\wt\D_\delta,
\wt D,\wt z_\delta,\wt z'_\delta),(\wt\D,\wt D,\wt\Gamma(0),\wt\Gamma(\xi)))\\
&\qquad+\dd_{GH\bullet\bullet}((\wt\D,\wt D,\wt\Gamma(0),\wt\Gamma(\xi)), (V(\u_L),d_L,\Theta'_L(0),\Theta'_L(n_L))),
\end{align*}
we obtain (for all but countably many values of $\delta$),
$$\limsup_{L\to\infty} \dd_{GH\bullet\bullet}((\v_{L,\delta},d_L,z_{L,\delta},z'_{L,\delta}),(V(\u_L),d_L,\Theta'_L(0),\Theta'_L(n_L)))
\leq \dd_{GH\bullet\bullet}((\wt\D_\delta,
\wt D,\wt z_\delta,\wt z'_\delta),(\wt\D,\wt D,\wt\Gamma(0),\wt\Gamma(\xi))).$$
This completes the proof since the right-hand side tends to $0$ as $\delta\to 0$. \endproof

\subsection{Comparing the length of discrete and continuous paths}
\label{compa-length}

We again use Skorokhod's representation theorem, which allows us to assume that \eqref{convGHPU} holds
almost surely. 
We then fix $\omega$ (in the underlying probability space)
such that \eqref{convGHPU} holds. By Proposition \ref{GHPU-criterion}, we may
embed all measure metric spaces $(V(\t_L),\sqrt{3/2}\,L^{-1/2}\dg)$ and $(\D,D)$ isometrically in the same compact 
metric space $(E,\Delta)$, in such a way that
$$V(\t_L)\build{\la}_{L\to\infty}^{}\D$$
for the Hausdorff distance $\Delta_{\rm Haus}$, the distinguished point $v_*^{(L)}$ of $V(\t_L)$
converges to the distinguished point $\xx_*$ of $\D$, the boundary path $\wh \Theta_L$
converges uniformly to $\Gamma$ (and in particular $\partial\t_L$
converges to $\partial\D$ for the Hausdorff distance), and we also have the weak convergence
$$\nu_L\build{\la}_{L\to\infty}^{} \mathbf{V}.$$
Notice that we have then $r_0^{(L)}\la r_0$, where we recall that $r_0=D(x_*,\partial \D)$ and $r_0^{(L)}$ is the rescaled graph
distance between $v_*^{(L)}$ and $\partial \t_L$.

We state a simple lemma that follows from the preceding convergences. The notation $\mathrm{Length}(\gamma)$
stands for the length of a continuous path $\gamma$ in $(\D,D)$. 

\begin{lemma}
\label{easy-lemma}
Let $\omega$ be fixed as explained above, and let $A>0$. For every $\ve>0$, we can find $L_0$
such that the following properties hold for every $L\geq L_0$.

\smallskip
\noindent{\rm(i)} Let $(\gamma(t))_{0\leq t\leq 1}$ be a continuous path in $\D$. Then we can find a discrete path
$\gamma_L=(u_0,u_1,\ldots,u_p)$ in $\t_L$ such that $\Delta(\gamma(0),u_0)<\ve$, $\Delta(\gamma(1),u_p)<\ve$
and $\Delta(u_i,\{\gamma(t):0\leq t\leq 1\})<\ve$ for every $i\in\{1,\ldots,p-1\}$.
If $\mathrm{Length}(\gamma)\leq A$, we can construct $(u_0,u_1,\ldots,u_p)$ so that
the rescaled $\dg$-length of $\gamma_L$ is bounded above by $\mathrm{Length}(\gamma)+\ve$. In
all cases, if
$\gamma(1)\in\partial \D$, we can take $u_p\in\partial\t_L$.

\smallskip
\noindent{\rm(ii)} If $(u_0,u_1,\ldots,u_p)$ is a discrete path in $\t_L$, we can find a continuous 
path $\gamma=(\gamma(t))_{0\leq t\leq 1}$ in $\D$ such that $\Delta(\gamma(0),u_0)<\ve$, $\Delta(\gamma(1),u_p)<\ve$
and, for every $t\in [0,1]$, $\Delta(\gamma(t),\{u_0,u_1,\ldots,u_p\})<\ve$. If the rescaled $\dg$-length of 
$(u_0,\ldots,u_p)$ is smaller than $A$, we can construct $\gamma$ so that
the length of $\gamma$ is bounded above by the rescaled $\dg$-length of $(u_0,\ldots,u_p)$ plus $\ve$. In all cases, if $u_p\in\partial \t_L$, we can
take $\gamma(1)\in\partial \D$.
\end{lemma}

\proof (i) Let $N\geq 2$ be an integer such that $N(\ve/4)\geq A$. We choose $L_0$ large enough so that $\Delta_{\rm Haus}(V(\t_L), \D)<\frac{\ve}{4N}$
and $\Delta_{\rm Haus}(\partial\t_L, \partial\D)<\frac{\ve}{4N}$ for every $L\geq L_0$.
Then, let $0=t_0<t_1<\cdots<t_k=1$
be such that $D(\gamma(t_{j-1}),\gamma(t_j))\leq\ve/4$ for every $1\leq j\leq k$. For every 
$0\leq j\leq k$, we can find $v_j\in V(\t_L)$ with $\Delta(v_j,\gamma(t_j))<\ve /(4N)$
(and in the case $\gamma(1)\in \partial \D$ we can take $v_k\in\partial\t_L$).  We note that, for $1\leq j\leq N$,
$$\Delta\big(v_{j-1},v_j\big)\leq \Delta\big(v_{j-1},\gamma(t_{j-1})\big)+\Delta\big(\gamma(t_{j-1}),\gamma(t_{j})\big)
+ \Delta\big(\gamma(t_{j}),v_j\big)\leq \Delta\big(\gamma(t_{j-1}),\gamma(t_{j})\big)+\frac{\ve}{2N}\leq\frac{\ve}{2}.$$
We construct
the path $(u_0,\ldots,u_p)$ as the concatenation of (graph distance) geodesics between
$v_{j-1}$ and $v_j$, for $1\leq j\leq k$. Then, for $i\in\{0,1,\ldots,p\}$, the vertex $u_i$
belongs to the geodesic between
$v_{j-1}$ and $v_j$, for some $j$, and $\Delta(u_i,\gamma(t_j))\leq \Delta(u_i,v_j)+\Delta(v_j,\gamma(t_j))
\leq \Delta(v_{j-1},v_j) +\ve/4<\ve$. 

If $\mathrm{Length}(\gamma)\leq A$, we observe that we can then take $k=N$ in the preceding considerations
(thanks to the condition $N(\ve/4)\geq A$). It follows that the rescaled $\dg$-length of the path 
$(u_0,\ldots,u_p)$ is equal to
$$\sum_{j=1}^N \Delta(v_{j-1},v_j)\leq \sum_{j=1}^N  \Delta(\gamma(t_{j-1}),\gamma(t_{j})) + N\,\frac{\ve}{2N} \leq \mathrm{Length}(\gamma)+\frac{\ve}{2}.$$

The proof of (ii) is similar and omitted. \endproof

\subsection{The key technical result}
\label{key-Tec}

Recall the convergences in distribution \eqref{convGHPU} and \eqref{convGHPU2}. By a tightness
argument, we may find a sequence $L_n\uparrow \infty$ along which these
two convergences hold jointly. From now on, we restrict our attention to
values of $L$ belonging to this sequence. Then, by an application of the 
Skorokhod representation theorem, we may assume that both convergences
 \eqref{convGHPU} and \eqref{convGHPU2}
hold a.s. 

\begin{proposition}
\label{key-techni}
Under the preceding assumptions, the bipointed metric space $(\wt \D,\wt D,\wt\Gamma(0),\wt\Gamma(\xi))$ 
coincides a.s. with the space $(\U,D_\star,\Pi_\star(0),\Pi_\star(\xi))$.
\end{proposition}

The proof of Proposition \ref{key-techni} occupies the remaining part of this section, and will
rely on several lemmas.
From now on until the end of this section, we fix $\omega$ such that both convergences
 \eqref{convGHPU} and \eqref{convGHPU2}
hold. We may and will also assume that $\omega$ has been chosen so that the properties of Proposition \ref{techni} hold. 

 We embed all measure metric spaces $(V(\t_L),\sqrt{3/2}\,L^{-1/2}\dg)$ and $(\D,D)$ isometrically in the same compact 
metric space $(E,\Delta)$ in the way explained at the beginning of Section \ref{compa-length}. 
In the following, the word
``distance'' refers to the distance $\Delta$ on the space $E$ in which the measure metric spaces
$(V(\t_L),\sqrt{3/2}\,L^{-1/2}\dg)$ and $(\D,D)$ are embedded isometrically.

Recall that $\W_L$ denotes the set of all vertices of the revealed region at time $\zeta_L$, and that
$\v_L$ is the set of all vertices of the unrevealed region at the same time. We note that $V(\t_L)= \v_L\cup\W_L$,
and that $\v_L\cap\W_L=\partial \W_L$ is the boundary of the revealed region. We also observe that
$\partial \W_L \cap \partial\t_L$ consists of the single point $\alpha_L$, and every point of $\v_L$
except $\alpha_L$ corresponds to a single point of $V(\u_L)$. We already noticed that any path from $\W_L$
to $\partial \t_L$ must visit a vertex whose distance from $v_*^{(L)}$ is bounded above
by $r_0^{(L)}$.

\begin{lemma}
\label{conv-set-lem}
For $\omega$ fixed as above, we have both
\begin{equation}
\label{conv-set}
\v_L\build{\la}_{L\to\infty}^{}U
\end{equation}
and
\begin{equation}
\label{conv-set2}
\W_L\build{\la}_{L\to\infty}^{} H
\end{equation}
for the Hausdorff distance in $(E,\Delta)$.
\end{lemma}

\proof For any (nonempty) subset $A$ of $E$, we write $A_\ve:=\{x\in E: \Delta(x,A)<\ve\}$. Let us first discuss the convergence
\eqref{conv-set}. Fix $\ve>0$. We start by verifying that 
$U\subset (\v_L)_\ve$ for all large enough $L$. To this end, first suppose that $x\in \mathrm{Int}(U)$. Then
$\Delta(x,H)>0$ and we can find a continuous path from $x$
to $\partial \D$ that does not intersect $H$,
and therefore this path stays at 
distance at least $r_0+\delta$ from the distinguished point $\xx_*$, for some
$\delta>0$ depending on $x$. For $L$ large enough, Lemma \ref{easy-lemma} (i) then allows us to find 
$v_{(L)}\in V(\t_L)$ such that $\Delta(x,v_{(L)})<\ve$, and a path from $v_{(L)}$ to $\partial\t_L$ in $\t_L$ that stays at distance at least $r_0+\delta/2>r_0^{(L)}$
from the distinguished point $v_*^{(L)}$.
Then $v_{(L)}$ must belong to $\v_L$, and we conclude that $x\in (\v_L)_\ve$
for $L$ large enough. This also holds when $x\in \partial U$, because
we can then approximate $x$ by a point of $\mathrm{Int}(U)$ close to $x$ and use the 
preceding conclusion with $\ve$ replaced by $\ve/2$. We can now use a compactness
argument to verify that $U\subset (\v_L)_\ve$ for all sufficiently large $L$.
Indeed, if this does not hold, we can find a sequence $L_k\uparrow \infty$
and, for every $k$, a point $x_k\in U$ that does not belong to $(\v_{L_k})_\ve$. Up 
to extracting a subsequence we can assume that $x_k\la x_\infty\in U$, and the fact 
that $x_\infty\in (\v_L)_{\ve/2}$ for $L$ sufficiently large forces $x_k\in (\v_{L_k})_\ve$
for $k$ large, which is a contradiction. 

To complete the proof of \eqref{conv-set}, we also need to verify that
$\v_L\subset U_\ve$ for $L$ large. If $v\in \v_L$, we can find a 
path from $v$ to $\partial\t_L$ that stays at distance at least $r_0^{(L)}-\sqrt{2/3}\,L^{-1/2}$
from $v_*^{(L)}$. Assuming that $L$ is large enough, Lemma \ref{easy-lemma} (ii) allows us to find $x\in \D$
with $\Delta(x,v)<\ve/2$, such that there exists a path from $x$ to $\partial \D$
that stays at distance at least $r_0-\delta$ from $\xx_*$, where
$\delta>0$ can be chosen arbitrarily small (note that the choice of $L$ ``large enough'' depends on $\ve$ and $\delta$, but
can be made uniformly in $v$). From the property stated in Proposition \ref{techni} (i), if
$\delta$ has been chosen small enough, this implies that
$x\in U_{\ve/2}$ and thus $v\in U_\ve$.

Let us turn to \eqref{conv-set2}. We need to verify that, for every $\ve >0$, we have for $L$ large 
enough,
\begin{equation}
\label{inclu1}
H\subset (\W_L)_\ve
\end{equation}
and
\begin{equation}
\label{inclu2}
\W_L \subset H_\ve.
\end{equation}
Consider the second inclusion \eqref{inclu2}. It is enough to verify that, for $L$ sufficiently large, if $v\in V(\t_L)$
and $v\notin H_\ve$, then necessarily $v\in V(\t_L)\backslash\W_L$. By the property stated in Proposition \ref{techni} (ii), we can find $\delta=\delta(\ve)>0$
such that, if $x\in\D\backslash H_{\ve/2}$, there is a path from $x$ to $\partial \D$ in $\D$ that stays at distance 
greater than $r_0+\delta(\ve)$ from $\xx_*$ (and thus at distance greater than $\delta(\ve)$ from $H$). Then, if $v\in V(\t_L)$
and $v\notin H_\ve$, we can find (for $L$ large, independently of the choice of $v$)
$x\in \D$ such that $\Delta(v,x)<\alpha$, where $\alpha=(\ve/2)\wedge (\delta(\ve)/4)$. Since 
$\alpha\leq \ve/2$, it is immediate that $x\notin H_{\ve/2}$, hence we can find a 
path from $x$ to $\partial \D$ that stays at distance greater than $r_0+\delta(\ve)$ from $\xx_*$.
Using Lemma \ref{easy-lemma} (i), we then obtain that, for $L$ large (independently of our choice of $v$),
there is a path from $v$ to $\partial\t_L$ in $V(\t_L)$ that stays at distance greater than $r_0^{(L)}+ \delta(\ve)/2$
from $v_*^{(L)}$. This implies that $v\in V(\t_L)\backslash\W_L$ as required. 

Let us finally discuss \eqref{inclu1}. Fix $\alpha\in(0,(\ve/3)\wedge r_0)$, and let
$x\in H$ such that $\Delta(x,U)\geq \alpha$. Since we know by the first part of
the proof that $\v_L$ converges to $U$ in the Hausdorff metric, we have 
$\v_L\subset U_{\alpha/2}$ for $L$ large, and this implies that $\Delta(x,\v_L)\geq \alpha/2$. 
Assuming that $L$ is even larger, we can find $v_{(L)}\in \t_L$ such that 
$\Delta(x,v_{(L)})<\alpha/4$, and then $\Delta(v_{(L)},\v_L)\geq \alpha/4$, which implies
that $v_{(L)}\in V(\t_L)\backslash \v_L$. We conclude that, for $L$ large,
$$\{x\in H:\Delta(x,U)\geq \alpha\} \subset (V(\t_L)\backslash \v_L)_{\alpha/4}\subset (\W_L)_{\alpha/4}.$$
On the other hand, if $x\in H$ is such that $\Delta(x,U)<\alpha$, then necessarily 
$\Delta(\xx_*,x)\leq r_0+\alpha$ (points of $\partial H=\partial U$ are at
distance $r_0$ from $\xx_*$) and we can find on a geodesic
from $x$ to $\xx_*$ a point $x'$ such that $\Delta(\xx_*,x')=r_0-\alpha$ and $\Delta(x,x')\leq 2\alpha$. Clearly $\Delta(x',U)\geq \alpha$,
and $x'$ belongs to the set in the left-hand side of the last display. We conclude that $x\in 
(\W_L)_{9\alpha/4}$. By combining the two cases, and using the property $\alpha<\ve/3$,
we obtain that $H\subset (\W_L)_\ve$, which completes the proof of \eqref{inclu1} and Lemma \ref{conv-set-lem}. \endproof

Recall from Section \ref{sec:identi} the definition of the metric space 
$(U',d_\infty)$ and of its ``boundary'' $\partial U'=U'\backslash \mathrm{Int}(U)$. For $\delta>0$, we set
$$
U_{(\delta)}:=\{x\in U': d_\infty(x,\partial U')\geq \delta\}
$$
(this set is not empty if  $\delta$ is small enough).
We state an analog of Lemma \ref{avoiding1}, which is  easier to establish.

\begin{lemma}
\label{avoiding2}
Let $\delta>0$ and $\eta>0$. Then 
there exists $\ve_0>0$
such that, for every $x,y\in U_{(\delta)}$,
there is a path from $x$ to $y$ in $\mathrm{Int}(U)$ that stays at $d_\infty$-distance at least $\ve_0$
from $\partial U'$, and whose $d_\infty$-length is bounded above
by $d_\infty(x,y)+\eta$. 
\end{lemma}

\proof 
If the conclusion of the lemma fails, then we can find a sequence $\ve_k\da 0$,
and, for every $k$, two points $x_k$ and $y_k$ in $U_{(\delta)}$ such that any path connecting $x_k$ to $y_k$
and  staying at $d_\infty$-distance at least $\ve_k$
from $\partial U'$ has length at least $d_\infty(x_k,y_k)+\eta$. By compactness, we may assume that 
$x_k\la x_\infty$ and $y_k\la y_\infty$, with $x_\infty,y_\infty\in U_{(\delta)}$. By the definition
of the intrinsic distance $d_\infty$, there is a curve $(\gamma(t))_{0\leq t\leq 1}$ connecting $x_\infty$
to $y_\infty$ and staying in $\mathrm{Int}(U)$ whose length is smaller than $d_\infty(x,y)+\eta/2$.
We can find $\alpha>0$ such that $d_\infty(\gamma(t),\partial U')>\alpha$ for every $t\in[0,1]$.
By concatenating the path $\gamma$ with $D$-geodesics from $x_k$ to $x_\infty$ and from
$y_\infty$ to $y_k$, we get, for all large $k$, a path from $x_k$ to $y_k$ which stays at
distance at least $\alpha/2$ from $\partial U'$ and
whose length is smaller than $d_\infty(x,y)+\eta$. This gives a contradiction when $\ve_k<\alpha/2$. \endproof

Recall that Proposition \ref{ident-lim-space} yields a canonical identification of the measure metric spaces $(U',d_\infty,\mathbf{V}_U)$ and $(\U,D_\star,\mathbf{V}_\star)$,
so that
we may view $U_{(\delta)}$ as a closed subset of $\U$. We will point the metric space $(U_{(\delta)},d_\infty)$ at 
the distinguished points
$$z_\delta:=\Pi_\star(\inf\{t\in[0,\xi]:\ell_{\Pi_\star(t)}\geq \delta\}),\;z'_\delta:=\Pi_\star(\sup\{t\in[0,\xi]:\ell_{\Pi_\star(t)}\geq \delta\}),$$
which are the ``first'' and ``last'' points of $\partial_0\U$ at distance $\delta$ from $\partial_1\U$.
The definition of $z_\delta$ and $z'_\delta$ makes sense for $\delta$ small enough, which we assume from now on. In what follows, we may also assume that $U_{(\delta)}\not=\varnothing$.

\begin{lemma}
\label{main1}
We have 
$$\dd_{GH\bullet\bullet}((U_{(\delta)},d_\infty,z_\delta,z'_\delta), (\U,D_\star,\Pi_\star(0),\Pi_\star(\xi)))\build{\la}_{\delta\to 0}^{} 0.$$
\end{lemma}

This is immediate since $U_{(\delta)}$ is identified to the closed subset
$\{x\in\U:D_\star(x,\partial_1\U)\geq \delta\}$ which converges to $\U$
in the sense of the Hausdorff distance, and we have also $z_\delta\la \Pi_\star(0)$
and $z'_\delta\la \Pi_\star(\xi)$
as $\delta\to 0$.

\begin{lemma}
\label{main3}
For all but countably many values of $\delta$, we have
$$\dd_{GH\bullet\bullet}((\v_{L,\delta},d_L,z_{L,\delta},z'_{L,\delta}),(U_{(\delta)},d_\infty,z_\delta,z'_\delta))\build{\la}_{L\to\infty}^{} 0.$$
\end{lemma}

Assume for the moment that we have proved Lemma \ref{main3}. Then, writing
\begin{align*}
&\dd_{GH\bullet\bullet}((V(\u_L),d_L,\Theta'_L(0),\Theta'_L(n_L)), (\U,D_\star,\Pi_\star(0),\Pi_\star(\xi)))\\
&\qquad\leq \dd_{GH\bullet\bullet}((V(\u_L),d_L,\Theta'_L(0),\Theta'_L(n_L)),(\v_{L,\delta},d_L,z_{L,\delta},z'_{L,\delta}))\\
&\qquad+\dd_{GH\bullet\bullet}((\v_{L,\delta},d_L,z_{L,\delta},z'_{L,\delta}),(U_{(\delta)},d_\infty,z_\delta,z'_\delta))\\
&\qquad+\dd_{GH\bullet\bullet}((U_{(\delta)},d_\infty,z_\delta,z'_\delta), (\U,D_\star,\Pi_\star(0),\Pi_\star(\xi))),
\end{align*}
we get, except possibly for countably many values of $\delta$,
\begin{align*}
&\limsup_{L\to\infty}\dd_{GH\bullet\bullet}((V(\u_L),d_L,\Theta'_L(0),\Theta'_L(n_L)), (\U,D_\star,\Pi_\star(0),\Pi_\star(\xi)))\\
&\qquad\leq \dd_{GH\bullet\bullet}((U_{(\delta)},d_\infty,z_\delta,z'_\delta), (\U,D_\star,\Pi_\star(0),\Pi_\star(\xi))) \\
&\qquad+ \limsup_{L\to\infty} \dd_{GH\bullet\bullet}((V(\u_L),d_L,\Theta'_L(0),\Theta'_L(n_L)),(\v_{L,\delta},d_L,z_{L,\delta},z'_{L,\delta})).
\end{align*}
By Lemmas \ref{main2} and \ref{main1}, the right-hand-side can be made arbitrarily small
by taking $\delta$ small. Hence, we have proved that
$$\dd_{GH\bullet\bullet}((V(\u_L),d_L,\Theta'_L(0),\Theta'_L(n_L)), (\U,D_\star,\Pi_\star(0),\Pi_\star(\xi))) \build{\la}_{L\to\infty}^{} 0.$$
Now notice that \eqref{convGHPU2} (which we have assumed to hold pointwise
for the fixed value of $\omega$) implies the $\dd_{GH\bullet\bullet}$-convergence of the pointed
spaces $(V(\u_L),d_L,\Theta'_L(0),\Theta'_L(n_L))$ to $(\wt\D,\wt D,\wt\Gamma(0),\wt\Gamma(\xi))$.
The statement of Proposition \ref{key-techni} follows by comparing with the last display. \endproof

\smallskip

\proof[Proof of Lemma \ref{main3}]
Although $\v_{L,\delta}$
was defined as a subset of $V(\u_L)$, we may and will identify $\v_{L,\delta}$ with the
corresponding subset of $\v_L$ (recall that points of $\v_L$ are in one-to-one 
correspondence with the points of $V(\u_L)$, with the exception of $\alpha_L$ that corresponds to
two points of $V(\u_L)$), and we have also $\v_{L,\delta}=\{x\in \v_L: d_L(x,\mathcal{W}_L)\geq \delta\}$. 
We first observe that
\begin{equation}
\label{main31}
\v_{L,\delta}\build{\la}_{L\to\infty}^{} U_{(\delta)},
\end{equation}
in the sense of the Hausdorff distance $\Delta_{\mathrm{Haus}}$, 
except possibly for countably many values of $\delta$. The proof 
of \eqref{main31} is essentially the same as the proof of \eqref{main21} above, using
now the convergence of $\W_L$ to $H$ in Lemma \ref{conv-set-lem}. The reader may be 
puzzled by the fact that the limit of $\v_{L,\delta}$ in \eqref{main31} is different than
the limit in \eqref{main21}: The point is that $\v_{L,\delta}$ was viewed as a subset of $V(\u_L)$
in \eqref{main21}, and we were using (in the proof of Lemma \ref{main2}) a specific
embedding of the metric spaces $V(\u_L)$ whereas in this section we are using another specific embedding
of the spaces $V(\t_L)$. 

We then also observe that, for all but countably many values of $\delta$, we must have
\begin{equation}
\label{main32}
z_{L,\delta}\build{\la}_{L\to\infty}^{} z_\delta\,, \quad z'_{L,\delta}\build{\la}_{L\to\infty}^{} z'_\delta.
\end{equation}
Indeed, in the cyclic exploration of $\partial\u_L$ induced by the path $\Theta_L$,
$z_{L,\delta}$ is the first visited point after $\alpha'_L$ such that $d_L(z_{L,\delta},\W_L)\geq \delta$
and similarly, in the cyclic exploration
of $\partial\D$ induced by $\Gamma$, $z_\delta$ is the first point after $\xx_0$ such that $D(z_\delta,H)\geq \delta$. 
Moreover, it is easy to verify that $\alpha_L$ converges as $L\to\infty$
to the point $\xx_0$ (which is the unique point of $\partial\D$
at minimal distance from the distinguished point of $\D$). 
Taking the preceding facts into account, one can use the convergence of $\W_L$ to $H$ and the uniform convergence of $\wh\Theta_L$
to $\Gamma$ to derive the convergence of $z_{L,\delta}$ toward $z_{\delta}$. A similar argument
applies to the convergence of $z'_{L,\delta}$ toward $z'_{\delta}$.

Fix $\delta>0$ such that \eqref{main31} and \eqref{main32} hold. We then define a correspondence between $(\v_{L,\delta},d_L)$ and $(U_{(\delta)},d_\infty)$
by setting
$$\mathcal{C}_L=\{(x,x')\in \v_{L,\delta}\times U_{(\delta)}: \Delta(x,x')\leq \kappa_L\}$$
where 
$$\kappa_L=\Delta_{\mathrm{Haus}}(\v_{L,\delta},U_{(\delta)})+ \Delta(z_{L,\delta},z_\delta)
+ \Delta(z'_{L,\delta},z'_\delta).$$
Note that $\kappa_L\la 0$ as $L\to\infty$.
By construction, the pairs $(z_{L,\delta},z_\delta)$ and $(z'_{L,\delta},z'_\delta)$ belongs to $\mathcal{C}_L$. 
From the characterization of $\dd_{GH\bullet\bullet}$ in terms 
of correspondences (Section \ref{sec:convmetric}), it therefore suffices to verify that the distortion of $\mathcal{C}_L$
tends to $0$ as $L\to\infty$. 

We first verify that
\begin{equation}
\label{main33}
\sup_{(x,x')\in\cc_L,(y,y')\in\cc_L} \Big(d_L(x,y)-d_\infty(x',y')\Big) \build{\la}_{L\to\infty}^{} 0.
\end{equation}
To this end, fix $\eta\in(0,1)$, and choose $\ve_0$ as in Lemma \ref{avoiding2}. We can assume that
$\ve_0\leq \eta\wedge \frac{\delta}{2}$. 
Then choose $L_0$ as in Lemma \ref{easy-lemma}, with $\ve=\frac{\ve_0}{2}\wedge \frac{\eta}{8}$
and $A=\mathrm{diam}(U',d_\infty)+1$, where $\mathrm{diam}(U',d_\infty)$ denotes the diameter of the
metric space $(U',d_\infty)$. Taking $L_0$
larger if necessary, we can assume that 
$\kappa_L\leq \frac{\ve_0}{4}$ and $\Delta_{\mathrm{Haus}}(\W_L,H)\leq \frac{\ve_0}{2}$, for every $L\geq L_0$. 
Then, let $L\geq L_0$ and let $(x,x')\in\cc_L$ and $(y,y')\in\cc_L$. By Lemma \ref{avoiding2}, we can find
a path $(\gamma(t))_{t\in[0,1]}$ from $x'$ to $y'$ in $\mathrm{Int}(U)$, which stays at distance at least 
$\ve_0$ from $H$ and whose length is smaller than $d_\infty(x',y')+\eta$. 
By Lemma \ref{easy-lemma} (i), we can then find a discrete path $\gamma_L=(u_0,u_1,\ldots,u_p)$
such that $\Delta(\gamma(0),u_0)<\ve$, $\Delta(\gamma(1),u_p)<\ve$
and $\Delta(u_i,\{\gamma(t):0\leq t\leq 1\})<\ve$ for every $i\in\{1,\ldots,p-1\}$, and 
the rescaled $\dg$-length of $\gamma_L$ is bounded above by $\mathrm{Length}(\gamma)+\ve$. Since 
$\Delta_{\mathrm{Haus}}(\W_L,H)\leq \frac{\ve_0}{2}$ and $\gamma$ stays at distance at least
$\ve_0$ from $H$, it follows that the path $\gamma_L$ does not visit $\W_L$, and 
therefore $d_L(u_0,u_p)$ is bounded above by the rescaled length of $\gamma_L$
which is bounded by $d_\infty(x',y')+\eta+\ve$. 
Since $\Delta(x,x')\leq \kappa_L\leq \frac{\ve_0}{4}$, we have $\Delta(u_0,x)\leq \ve+\frac{\ve_0}{4}
<\ve_0\leq \frac{\delta}{2}$, and (since $x\in \v_{L,\delta}$) it follows that a geodesic from $u_0$ to $x$ in $\t_L$
does not intersect $\W_L$, and therefore $d_L(x,u_0)=\Delta(u_0,x)<\ve_0$,
and similarly $d_L(y,u_p)<\ve_0$. Finally we get $d_L(x,y)\leq d_L(u_0,u_p)+2\ve_0
\leq d_\infty(x',y')+3\eta+\ve$. This shows that the supremum in 
\eqref{main33} is bounded above by $4\eta$ when $L\geq L_0$, which completes the proof of \eqref{main33}.

It remains to verify that
\begin{equation}
\label{main34}
\sup_{(x,x')\in\cc_L,(y,y')\in\cc_L} \Big(d_\infty(x',y')-d_L(x,y)\Big) \build{\la}_{L\to\infty}^{} 0.
\end{equation}
This is very similar to the proof of \eqref{main33}. Fix $\eta>0$. Then,
we can use Lemma \ref{avoiding1} to find $\ve_0>0$ such that the following holds 
for all $L\geq L_1$, for some integer $L_1\geq 1$. For any 
$x,y\in \v_{L,\delta}$, there is a discrete path $\gamma_L=(u_0,\ldots,u_p)$ from 
$x$ to $y$ 
in $\v_L$, which stays at rescaled graph distance at least $\ve_0$ from $\W_L$
and whose rescaled length is bounded above by $d_L(x,y)+\eta$. Let $L_0$
be as in Lemma \ref{easy-lemma}, with $A$ such that 
$\mathrm{diam}(V(\u_L),d_L)\leq A$ for every $L$ (the convergence
\eqref{convGHPU2} allows us to find $A$ with this property) and $\ve=\frac{\ve_0}{4}\wedge\frac{\eta}{4}\wedge \frac{\delta}{4}$. 
Taking $L_0$ larger if necessary, we can assume that $L_0\geq L_1$ and
$\Delta_{\mathrm{Haus}}(\W_L,H)<\ve$ for every $L\geq L_0$, and moreover
$\kappa_L<\ve$ for every $L\geq L_0$.
Then, if $L\geq L_0$ and  $(x,x')\in\cc_L$ and $(y,y')\in\cc_L$, and if the discrete path $\gamma_L=(u_0,\ldots,u_p)$
from $x$ to $y$
is chosen as explained above, Lemma \ref{easy-lemma} (ii) allows us to find a 
continuous path $(\gamma(t))_{t\in[0,1]}$ such that
$\Delta(\gamma(0),u_0)<\ve$, $\Delta(\gamma(1),u_p)<\ve$
and, for every $t\in [0,1]$, $\Delta(\gamma(t),\{u_0,u_1,\ldots,u_p\})<\ve$, and moreover
the length of $\gamma$ is bounded above by the rescaled length of $\gamma_L$
plus $\ve$. Recalling that $\Delta_{\mathrm{Haus}}(\W_L,H)<\ve$
and that $\gamma_L$  stays at rescaled graph distance at least $\ve_0$ from $\W_L$, we see that
$\gamma$ does not visit $H$ because this would contradict the property 
$\Delta(\gamma(t),\{u_0,u_1,\ldots,u_p\})<\ve$ for every $t\in[0,1]$. It follows
that $d_\infty(\gamma(0),\gamma(1))\leq d_L(x,y)+\eta+\ve$. 
However, we have $\Delta(\gamma(0),x')\leq \Delta(\gamma(0),u_0)+\Delta(x,x')\leq 2\ve$,
and similarly $\Delta(\gamma(1),y')\leq 2\ve$, and it follows that a geodesic (in $\D$)
from $x'$ to $\gamma(0)$, or from $y'$ to $\gamma(1)$, does not hit $H$. We finally
conclude that $d_\infty(x',y')\leq d_\infty(\gamma(0),\gamma(1)) + 4\ve\leq d_L(x,y)+3\eta$.
Thus, for $L\geq L_0$, the supremum in \eqref{main34} is bounded above by $3\eta$.
This completes the proof of Lemma \ref{main3}. \endproof

\subsection{The main theorem}
\label{main-Th}

Recall the definition of the exit measures $\z_y$
in Section \ref{sna-mea}, and the point measures $\n=\sum_{i\in I}\delta_{(t_i,\omega_i)}$ and $\wt{\n}=\sum_{i\in I}\delta_{(t_i,\tilde\omega_i)}$ 
that were used in
Section \ref{sec:limspace} to define $\D$ and $\U$ respectively. We set
$$\z_\U:=\sum_{i\in I} \z_0(\omega_i).$$
A first-moment argument shows that $\z_\U<\infty$ a.s. Specifically, using formula \eqref{Lap-exit} one gets that $\mathbb{N}_y(\z_0)=1$ for every $y>0$, and consequently:
\begin{equation}\label{computation:expectation:z_u}
\mathbb{E}[\z_\U]=2\,\mathbb{E}\big[ \int_0^\xi\dd t~\N_{\sqrt{3}\be_t}( \z_0)\big]=2 \,\xi. 
\end{equation}

Recall the notation
$\partial_0\U=\Pi_\star([0,\xi])$ and
$\partial_1\U=\Pi_\star(\partial\mathfrak{T}^\star)$
introduced at the end of Section \ref{conf-lim-space}, We next introduce a path $\Gamma^\star $ parametrizing $\partial_0\U\cup \partial_1 \U$. 
For every  $i\in I$, we let $(\wt L^{0}_t(\omega_i))_{t\geq 0}$
be the (time changed) exit local time of $\omega_i$ from $(0,\infty)$, as given by formula \eqref{formu-exit-bis} with $y=0$.   In the time scale of the clockwise exploration $(\ee^\star_s)_{0\leq s\leq \Sigma^\star}$ of $\mathfrak{T}^\star$, each snake trajectory $\wt \omega_i$
corresponds to the  interval $[r_i,r_i+\sigma(\wt\omega_i)]$,
where
$$r_i:=\sum_{j\in I:t_j<t_i}\sigma(\wt\omega_i).$$ 
We then set, for every $s\in[0,\Sigma^\star]$,
$$\wt L^{\star}_s:= \sum_{i\in I}\wt L^{0}_{(s-r_i)^+}(\omega_i),$$
so that $\wt L^{\star}_s$ represents the total exit local time accumulated at $0$
by the clockwise exploration of $\mathfrak{T}^\star$ up to time $s$. 
We set $\kappa_\star(t):=\inf\{s\geq 0: \wt L^{\star}_s> t\}$
for every $t\in[0,\mathcal{Z}_\U)$, and $\kappa_\star(\z_\U):= \Sigma^\star$. 

Then, for every $t\in[0,\xi]$, we set $\Gamma^\star(t):=\Pi_\star(t)$,
and we also define $\Gamma^\star(t)$ for $\xi\leq t\leq\xi+\z_\U$ by setting 
\begin{equation}
\label{gamma-star}
\Gamma^\star(\xi+\z_\U-s):=\Pi_\star(\ee^\star_{\kappa_\star(s)}),\quad\hbox{for every }s\in [0,\z_\U].
\end{equation}
The two definitions of $\Gamma^\star(\xi)$ are consistent since $\ee^\star_{\kappa_\star(\z_\U)}=\ee^\star_{\Sigma^\star}=\xi$.
We also note that $\Gamma^\star(0)=\Gamma^\star(\xi+\z_\U)=\Pi_\star(0)$. Furthermore, we have
$\{\Gamma^\star(t):0\leq t\leq\xi\}=\partial_0\U$, and 
$\{\Gamma^\star(t):\xi\leq t\leq \xi+\z_\U\}\}=\Pi_\star(\partial\mathfrak{T}^\star)=\partial_1\U$, by the support property of the exit local time.
So the range of $(\Gamma^\star(s))_{0\leq s\leq \z_\U+\xi}$ is $\partial_0\U\cup\partial_1\U$.

We next observe that $\Gamma^\star$ is continuous, and is injective on the interval $[0,\xi+\z_\U)$, The injectivity property is a direct consequence of Proposition \ref{def-lim-space},
and we omit the details. On the other hand,  we already know that the function $[0,\xi]\ni t\mapsto \Pi_\star(t)$ is continuous
(from the continuity properties of 
$(a,b)\mapsto D_\star(a,b)$, cf. Section \ref{conf-lim-space}). So it remains to show that  $(\Gamma^\star(\xi+\z_\U-s))_{0\leq s\leq  \z_\U}$ is  continuous.
 Let us briefly justify this point. Since $[0,\Sigma^\star]\ni t\mapsto \Pi_\star(\mathcal{E}^\star_t)$ is continuous and $t\mapsto \kappa_\star(t)$
 is right-continuous, we only need to show that  if $\kappa_\star(t-)<\kappa_\star(t)$ then $\Pi_\star(\ee^\star_{\kappa_\star(t-)})
=\Pi_\star(\ee^\star_{\kappa_\star(t)})$. To this end, notice that if $\kappa_\star(t-)<\kappa_\star(t)$, the support property of the exit local time implies that  
all points of the form $\ee_u^\star$ with $u\in(\kappa_\star(t-),\kappa_\star(t))$ have a positive label.
Proposition \ref{def-lim-space} then ensures  that $\Pi_\star(\ee^\star_{\kappa_\star(t-)})
=\Pi_\star(\ee^\star_{\kappa_\star(t)})$, as desired.

\begin{theorem}
\label{identif-hull-comple}
Let $(\U,D_\star,\mathbf{V}_\star)$ 
be defined as in Section \ref{sec:limspace}, and let $\Gamma^\star$ be as above.
Then 
$(\U,D_\star,\mathbf{V}_\star,\Gamma^\star)$
is a curve-decorated free Brownian disk with a random perimeter $\xi+\z_\U$ distributed according to the
measure $\frac{3}{2}\,\xi^{3/2}\,z^{-5/2}\,\mathbf{1}_{\{z>\xi\}}\,\dd z$.
\end{theorem}

\proof By  Proposition \ref{key-techni}, we may suppose that there exists a curve-decorated free Brownian disk
$(\wt \D, \wt D, \wt {\mathbf{V}},\wt\Gamma)$ with perimeter $\xi+\y$ such that the random bipointed metric spaces
$(\U,D_\star,\Pi_\star(0),\Pi_\star(\xi))$ and $(\wt \D, \wt D, \wt\Gamma(0),\wt\Gamma(\xi))$
are a.s. equal, so that, in particular, we have $\Gamma^\star(0)=\wt\Gamma(0)$ and $\Gamma^\star(\xi)=\wt\Gamma(\xi)$. We now want to argue that the preceding equality 
of bipointed metric spaces carries over to the volume measures
and to the decorating curves, meaning that we have also $\mathbf{V}_\star=\wt{\mathbf{V}}$ and
$\Gamma^\star=\wt \Gamma$ (the equality $\Gamma^\star=\wt \Gamma$ will imply in particular that
$\z_\U=\y$, which has density $\frac{3}{2}\,\xi^{3/2}\,(\xi+x)^{-5/2}$ by Proposition \ref{asymptotics-Z}).

The equality $(\U,D_\star)=(\wt\D, \wt D)$ allows us to define the boundary $\partial\U$
(as the set of all points of $\U$ that have no neighborhood
homeomorphic to the unit disk), and we also know that $\partial \U$ is the range of a simple loop. It is easy to verify
that any point $x$ of $\U\backslash (\partial_0\U\cup\partial_1\U)$ has a neighborhood homeomorphic to
the unit disk (just note that the ``interior'' $\U\backslash (\partial_0\U\cup\partial_1\U)$ is identified 
to a subset of $\D$ in such a way that, in a sufficiently small neighborhood of $x$, the distance $D_\star$ coincides
with the distance $D$). It follows that $\partial \U$ is contained in $\partial_0\U\cup\partial_1\U$.
Now recall that the range of  the simple loop $(\Gamma^\star(s))_{0\leq s\leq \z_\U+\xi}$ is precisely $\partial_0\U\cup\partial_1\U$. It follows that the boundary of $\U$ is contained in a simple loop, which is
only possible if $\partial\U$ is the whole loop. We have thus $\partial \U= \partial_0\U\cup\partial_1\U$.

Let us then verify that $\mathbf{V}_\star=\wt{\mathbf{V}}$.
It follows from the main result of \cite{Hausdorff} that $\wt{\mathbf{V}}$ may
be defined by
$$\wt{\mathbf{V}}(A)=\mathbf{c}\,m_h(A)$$
for every Borel subset $A$ of $\wt\D$, where $\mathbf{c}>0$
is a constant, and $m_h$ stands for the Hausdorff 
measure with gauge function $h(r)=r^4\log\log(1/r)$.
To be precise, the results of \cite{Hausdorff} apply to
the Brownian sphere and not to the Brownian disk.
However, we may use the connections between
the Brownian sphere and the Brownian disk, and 
in particular Theorem 8 in \cite{Stars} showing that the 
complement of a hull in the Brownian sphere is a Brownian 
disk (this complement needs to be equipped with the
intrinsic distance, but this makes no difference for Hausdorff
measures as long as
we consider sets that do not intersect the boundary, and
on the other hand we know that the boundary has zero
volume measure, and also zero $h$-Hausdorff measure since its Hausdorff dimension is $2$, by \cite{Bet}). 

So to prove that $\mathbf{V}_\star=\wt{\mathbf{V}}$, it is enough to
verify that we have also $\mathbf{V}_\star(A)=\mathbf{c}\,m_h(A)$
for every Borel subset $A$ of $\U$. We may restrict our attention
to subsets of $\U\backslash\partial_1\U$ since we
know that $m_h(\partial\U)=0$ and $\mathbf{V}_\star(\partial_1\U)=0$. 
Then we can use the fact that $\U\backslash\partial_1\U$
is identified to the open subset $\mathrm{Int}(U)$,
and the previously mentioned result for the Brownian disk $\D$
to obtain that the equality $\mathbf{V}_\star(A)=\mathbf{c}\,m_h(A)$ 
holds for every Borel subset of $\U\backslash\partial_1\U$ --- here 
again the fact that we deal with the intrinsic distance on $\mathrm{Int}(U)$
instead of the distance of $\D$ makes no difference for Hausdorff measures.
This completes the proof of the equality $\mathbf{V}_\star=\wt{\mathbf{V}}$.

To complete the proof of the theorem, we rely on the next lemma.

\begin{lemma}
\label{ident-boundarysize}
Almost surely,  $\Gamma^\star$ is a standard boundary curve of  $(\U,D_\star,\mathbf{V}_\star)$. \end{lemma}

Assuming the result of the lemma, the proof of the theorem is easily completed.
Indeed, since $\Gamma^\star$ and $\wt\Gamma$ 
are both standard boundary curves, the property
$(\Gamma^\star(0),\Gamma^\star(\xi))=(\wt\Gamma(0),\wt\Gamma(\xi))$
can only hold if $\Gamma^\star=\wt\Gamma$, which was the desired result. \endproof

\noindent{\it Remark.} The reason for dealing with bipointed spaces throughout this section
is the fact that we need the equality $(\Gamma^\star(0),\Gamma^\star(\xi))=(\wt\Gamma(0),\wt\Gamma(\xi))$ to
identify $\Gamma^\star$. The equality $\Gamma^\star(0)=\wt\Gamma(0)$ alone would 
not be sufficient for this identification.

\smallskip

\proof[Proof of Lemma \ref{ident-boundarysize}]
Let $\mu_{\partial \U}$ be the boundary measure of the Brownian disk $\U$. The total mass $\mu_{\partial \U}(\partial\U)$ 
is the boundary size of $\U$, or equivalently of $\wt\D$, and is therefore equal to $\xi+\mathcal{Y}$, which has density 
$\frac{3}{2}\,\xi^{3/2}\,z^{-5/2}\,\mathbf{1}_{\{z>\xi\}}$ by Proposition \ref{asymptotics-Z} and the fact that
$\mathcal{Y}$ is distributed as $\xi\,\Lambda$.  

We already know that $(\Gamma^\star(t))_{0\leq t\leq\xi+ \z_\U}$ 
is a simple loop, and
we will verify that, almost surely for any continuous function $\Phi$ on $\partial \U$, we have
\begin{equation}
\label{identbr1}
\int_0^{\xi+\z_\U} \dd t \,\Phi(\Gamma^\star(t))=\int \mu_{\partial\U}(\dd x)\,\Phi(x).
\end{equation}
This will imply that the perimeter of  $\U$ is $\mathcal{Z}_\U+\xi$ (take $\Phi=1$), and 
then, by the very definition, that  $\Gamma^\star$ is a standard boundary curve of $U$.

Let us prove \eqref{identbr1}. We first observe that
the restriction of $\mu_{\partial\U}$ to $\partial_0\U$ is the pushforward of Lebesgue measure on $[0,\xi]$ under $\Pi_\star$. This is
easy from the identification of Proposition \ref{ident-lim-space}, the fact that the boundary measure $\mu_{\partial\D}$
of $\partial\D$ is the pushforward of
Lebesgue measure on $[0,\xi]$ under $\Pi$, and the approximations of the boundary measure by the volume measure restricted to a tubular
neighborhood of the boundary. It follows that
\begin{equation}
\label{identbr2}
\int_{0}^{\xi} \dd t \,\Phi(\Gamma^\star(t))=\int_0^\xi \dd s\,\Phi(\Pi_\star(s))=\int_{\partial_0\U} \mu_{\partial\U}(\dd x)\,\Phi(x),
\end{equation}
and in particular $\mu_{\partial\U}(\partial_0\U)=\xi$.
We then claim that we have also
\begin{equation}
\label{identbr3}
\int_{\xi}^{\xi+\z_\U} \dd t \,\Phi(\Gamma^\star(t))=\int_{\partial_1\U} \mu_{\partial\U}(\dd x)\,\Phi(x),
\end{equation}
which will complete the proof of \eqref{identbr1}. To derive \eqref{identbr3}, introduce, for every $\ve>0$, the measure $\nu_\ve$ on $\U$
defined by
$$\langle\nu_\ve,\Phi\rangle=\ve^{-2} \int \mathbf{V}_\star(\dd x)\,\mathbf{1}_{\{D_\star(x,\partial_1 \U)<\ve\}} \Phi(x).$$
It follows from the approximations of the boundary measure of $\U$ that
$$\lim_{\ve\to 0} \langle\nu_\ve,\Phi\rangle = \int_{\partial_1\U} \mu_{\partial\U}(\dd x)\,\Phi(x).$$
On the other hand, we can also verify that $\langle\nu_\ve,\Phi\rangle$ converges to the left-hand side of \eqref{identbr3}.
Consider first the case $\Phi=1$. Then, we know that $\langle\nu_\ve,1\rangle$ converges to 
$\mu_{\partial \U}(\partial_1\U)= \mu_{\partial \U}(\partial \U)-\xi$. From the first observation of the proof, we get that
 $\mu_{\partial \U}(\partial_1\U)$ has density $\frac{3}{2}\,\xi^{3/2}\,(\xi+z)^{-5/2}$, and in particular, $\E[\mu_{\partial \U}(\partial_1\U)]=2\xi$. On the other hand, recall the notation 
introduced at the beginning of Section \ref{main-Th}, and in particular the fact that each $\wt \omega^i$, $i\in I$,
corresponds to the interval $[r_i,r_i+\sigma(\wt\omega_i)]$ in the time scale of the clockwise exploration $\mathcal{E}^\star$.
Also note that we have $D_\star(x,\partial_1 \U)=\ell_x=\wh W_{s-r_i}(\wt\omega_i)$ when $x=\Pi_\star(\ee^\star_s)$ with $s\in[r_i,r_i+\sigma(\wt\omega_i)]$.
Using Fatou's lemma, we have
$$\mu_{\partial \U}(\partial_1\U)=\lim \limits_{\ve\to 0} \langle\nu_\ve, 1 \rangle
\geq   \sum \limits_{i\in I}\Big(\liminf \limits_{\ve \to 0}\ve^{-2}\int_0^{\sigma(\tilde\omega_i)} \dd t\,\mathbf{1}_{\{\wh{W}_t(\tilde{\omega}_i)<\ve\}} \Big)= \sum \limits_{i\in I}\z_{0}(\omega_i)=\z_\U,$$
and we know from \eqref{computation:expectation:z_u} that $\E[\z_\U]=2\xi= \E[\mu_{\partial \U}(\partial_1\U)]$. It
follows from the last display that $\mu_{\partial \U}(\partial_1\U)=\z_\U$ a.s. 

Let us consider now a general continuous function $\Phi$ on $\U$. Without loss of generality we can assume that $0\leq\Phi\leq 1$. 
By Fatou's lemma and the definition of $V_\star$, we have
 $$\int_{\partial_1\U} \mu_{\partial\U}(\dd x)\,\Phi(x)=\lim \limits_{\ve\to 0} \langle\nu_\ve, \Phi\rangle \geq  \sum \limits_{i\in I}\Big(\liminf \limits_{\ve \to 0}\ve^{-2}\int_0^{\sigma(\tilde\omega_i)} \dd t\,\mathbf{1}_{\{\wh{W}_t(\tilde{\omega}_i)<\ve\}} \Phi\big(\Pi_\star(\mathcal{E}^\star_{r_i+t})\big)\Big).$$
 Moreover, \eqref{formu-exit-bis} entails that, for every $i\in I$, the measures $\ve^{-2} \mathbf{1}_{\{\wh{W}_t(\tilde{\omega}_i)<\ve\}} \mathbf{1}_{[0,\sigma(\tilde\omega_i)]}(t) \dd t$ converge  weakly to $\dd \wt L_t^0(\omega^i)$ as $\ve \to 0$. Using the definition of $L^\star_t$, it follows that
 $$\int_{\partial_1\U} \mu_{\partial\U}(\dd x)\,\Phi(x)=\lim \limits_{\ve\to 0} \langle\nu_\ve, \Phi\rangle\geq 
 \int_0^{\Sigma^\star} \dd \wt L^\star_t~ \Phi\big( \Pi_\star(\mathcal{E}^\star_t)\big)= \int_\xi^{\xi+\z_\U}\dd s ~\Phi( \Gamma^\star(s)),$$
 where the last equality holds
 by the definition of $\Gamma^\star(s)$ when $\xi\leq s\leq \xi+\z_\U$.
If we combine the bound of the last display with the same bound when $\Phi$ is replaced by $1-\Phi$ (using $\mu_{\partial \U}(\partial_1\U)=\z_\U$),
we arrive at the desired equality \eqref{identbr3}. This completes the proof of Lemma \ref{ident-boundarysize} and Theorem \ref{identif-hull-comple}. \endproof

\noindent{\it Remarks.} (i) Another way of verifying that the boundary size of $\U$ is $\xi+\z_\U$ would have been to prove the convergence
in distribution of $L^{-1}Z_L$ to 
$\z_\U$ (independently of Proposition \ref{asymptotics-Z})  and to check that this convergence holds jointly with the convergence
of $(V(\u_L),d_L)$ to $(U',d_\infty)=(\U,D_\star)$.

\smallskip
\noindent (ii) It follows from Theorem \ref{identif-hull-comple} 
that the density of the distribution of $\z_\U$ is 
 $\frac{3}{2}\,\xi^{3/2}\,(\xi+z)^{-5/2}$. This can be verified by the following direct 
 calculation. 
 For $\lambda> 0$, we have
$$
\E[e^{-\lambda \z_\U}\mid \be]= \exp\Big(-2\int_0^\xi \dd t\, \N_{\sqrt{3}\be_t}(1-\exp(-\lambda \z_0))\Big)=\exp\Big(-\int_0^\xi \dd t\,
(\be_t+ (2\lambda)^{-1/2})^{-2}\Big),
$$
using formula \eqref{Lap-exit}. 
Then, for every $\alpha>0$,
$$\E\Big[\exp\Big(-\int_0^\xi \dd t\,
(\be_t+\alpha)^{-2}\Big)\Big]=1-\sqrt{\frac{\pi}{2}}\,\Big(\frac{\xi}{\alpha^{2}}\Big)^{3/2}\,\chi_1(\frac{\xi}{2\alpha^2}),$$
where $\chi_1(x)=\frac{1}{\sqrt{\pi}}x^{-1/2} -e^x\,\mathrm{erfc}(\sqrt{x})$. 
This formula is the special case $F=1$ of Lemma \ref{key-form} below.
It follows that
$$\E[e^{-\lambda \z_\U}]= 1-2\sqrt{\pi}\,(\lambda\xi)^{3/2}\,\chi_1(\lambda\xi).$$
To invert the Laplace transform, start from
the classical formula
$$\int_0^\infty e^{-\lambda x}\,\frac{\dd x}{\sqrt{\pi(x+1)}} = \frac{1}{\sqrt{\lambda}}\,e^\lambda \mathrm{erfc}(\sqrt{\lambda}).$$
Two integrations by parts then give
$$\frac{3}{2}\int_0^\infty (x+1)^{-5/2}\,e^{-\lambda x}\,\dd x= 1- 2\lambda +2\sqrt{\pi} \lambda^{3/2}\,e^{\lambda}\,\mathrm{erfc}(\sqrt{\lambda})=1-2\sqrt{\pi}\,\lambda^{3/2}\,\chi_1(\lambda)=\E[e^{-\lambda \z_\U/\xi}],$$
and we conclude that the density of $\z_\U/\xi$ is $\frac{3}{2}(x+1)^{-5/2}$.

\section{Peeling the Brownian disk}
\label{section:comple}
Our goal in this section is to discuss a peeling exploration for Brownian disks. In particular, we will see that the complement of a hull centered at a boundary point in a Brownian disk
is again a Brownian disk. Our study relies on the representation derived in Theorem \ref{identif-hull-comple}. 
\subsection{Preliminary distributional identities}
\label{preli-dist}

We first need to introduce some notation. We fix $\xi>0$ and $r>0$. As previously, we write $(\be_t)_{0\leq t\leq \xi}$ for a positive Brownian excursion 
of duration $\xi$. We then consider a five-dimensional Bessel process $(X_t)_{t\geq 0}$ started from $r$. Since $r$ is a regular point
for this Markov process, we can define an (infinite) excursion measure away from $r$, and  we can also make sense of the law of the excursion of duration $\xi$ above level $r$ for the process $X$. We let 
$(\ov\be^{(r)}_t)_{0\leq t\leq \xi}$ be distributed according to this law, and set $\be^{(r)}_t:=\ov\be^{(r)}_t-r$ for $0\leq t\leq \xi$,
so that $\be^{(r)}$ starts and ends at $0$. We let $C([0,\xi],\R)$ stand for the set of 
all continuous functions from $[0,\xi]$ into $\R$, which is equipped with the sup norm. 

\begin{lemma}
\label{key-form}
For every $u\geq 0$, set
$$\Phi(u):= 1-2u + 2\sqrt{\pi}u^{3/2}e^u\mathrm{erfc}(\sqrt{u})=1-2\sqrt{\pi}u^{3/2}\chi_1(u).$$
Then, for every bounded continuous function $F:C([0,\xi],\R)\la \R$, we have
$$\E\Big[F\Big((\be^{(r)}_t)_{0\leq t\leq \xi}\Big)\Big] 
=\Phi(\frac{\xi}{2r^2})^{-1}\,\E\Big[ F\Big((\be_t)_{0\leq t\leq \xi}\Big)\,\exp\Big(-\int_0^\xi \frac{\dd s}{(\be_s+r)^2}\Big)\Big] .$$
\end{lemma}

\proof
For every $\ve>0$, let $(\bbf^{(\ve,r)}_t)_{0\leq t\leq \xi}$ be distributed as a five-dimensional Bessel process 
started from $r+\ve$ and conditioned to hit $r$ exactly at time $\xi$. See Proposition 3 in \cite{Repre} for a 
precise definition, noting that this proposition deals with a Bessel process of dimension $-1$
instead of a five-dimensional Bessel process, but this replacement gives the same conditioned process
because of the $h$-process relation linking the Bessel processes of dimension $5$
and of dimension $-1$ (cf. formula (4) in \cite{Repre}). By standard arguments, $(\bbf^{(\ve,r)}_t)_{0\leq t\leq \xi}$
converges in distribution to $(\ov\be^{(r)}_t)_{0\leq t\leq \xi}$ as $\ve\to 0$, and therefore
\begin{equation}
\label{key-form10}
\lim_{\ve \to 0} \E\Big[F\Big((\bbf^{(\ve,r)}_t-r)_{0\leq t\leq \xi}\Big)\Big]
=\E\Big[F\Big((\be^{(r)}_t)_{0\leq t\leq \xi}\Big)\Big].
\end{equation}
On the other hand, Lemma 5 in \cite{Repre} shows that
\begin{equation}
\label{key-form11}
\E\Big[F\Big((\bbf^{(\ve,r)}_t-r)_{0\leq t\leq \xi}\Big)\Big]=\frac{r+\ve}{r}\,\frac{q_\xi(\ve)}{\rho_\xi(\ve,r)}\,
\E\Big[F\Big((\bg^{(\ve)}_t)_{0\leq t\leq \xi}\Big)\,\exp\Big(-\int_0^\xi \frac{\dd s}{(\bg^{(\ve)}_s+r)^2}\Big)\Big],
\end{equation}
where $(\bg^{(\ve)}_s)_{0\leq s\leq \xi}$ is distributed as a linear Brownian motion started from 
$\ve$ and conditioned to hit $0$ exactly at time $\xi$, and the functions $q_\xi(\ve)$
and $\rho_\xi(\ve,r)$ are given by
$$q_\xi(\ve):=\frac{\ve}{\sqrt{2\pi \xi^3}}\,\exp\Big(-\frac{\ve^2}{2\xi}\Big)$$
and
$$\rho_\xi(\ve,r):=\ve\,e^{-\ve^2/(2\xi)}
\Bigg( \frac{1}{2r^3} \mathrm{erfc}\Big(\frac{\sqrt{\xi}}{r\sqrt{2}} + \frac{\ve}{\sqrt{2\xi}}\Big)
\exp\Big( \Big(\frac{\sqrt{\xi}}{r\sqrt{2}} + \frac{\ve}{\sqrt{2\xi}}\Big)^2\Big) 
-\frac{1}{r^2\sqrt{2\pi \xi}}+ \frac{r+\ve}{r\sqrt{2\pi \xi^3}}\Bigg).$$
Since 
$$\lim_{\ve\to 0} \E\Big[F\Big((\bg^{(\ve)}_t)_{0\leq t\leq \xi}\Big)\,\exp\Big(-\int_0^\xi \frac{\dd s}{(\bg^{(\ve)}_s+r)^2}\Big)\Big]
= \E\Big[ F\Big((\be_t)_{0\leq t\leq \xi}\Big)\,\exp\Big(-\int_0^\xi \frac{\dd s}{(\be_s+r)^2}\Big)\Big],$$
the desired result will follow from \eqref{key-form10} and \eqref{key-form11} if we can verify that
$$\lim_{\ve\to 0} \frac{\rho_\xi(\ve,r)}{q_\xi(\ve)}= \Phi(\frac{\xi}{2r^2}).$$
This is immediately checked from the formulas for 
$q_\xi(\ve)$
and $\rho_\xi(\ve,r)$. \endproof

In what follows, we will consider $\be^{(r/\sqrt{3})}$ rather than $\be^{(r)}$, and in order to
simplify notation we set $r'=r/\sqrt{3}$.  As in Section \ref{sec:limspace}, we assume that, conditionally on $\be$, $\n=\sum_{i\in I}\delta_{(t_i,\omega_i)}$
is a   measure on $[0,\xi]\times \S$ with intensity 
$$2\,\dd t\,\N_{\sqrt{3}\be_t}(\dd\omega).$$
We write $\z=\int \z_0(\omega)\,\n(\dd t\,\dd\omega)=\sum_{i\in I} \z_0(\omega_i)$ for the total exit measure at $0$ of the atoms $\omega_i$
($\z=\z_\U$ in the notation of Section \ref{pass-lim}) and $\mathrm{tr}_0(\n)=\sum_{i\in I}\delta_{(t_i,\mathrm{tr}_0(\omega_i))}$.
Furthermore, conditionally given $\be^{(r')}$, we let $\wt \n$ be distributed as
a Poisson measure on $[0,\xi]\times \S$ with intensity
\begin{equation}
\label{inten-Pois}2\,\mathbf{1}_{\{W_*(\omega)>-r\}}\,\dd t\,\N_{\sqrt{3}\be^{(r')}_t}(\dd\omega).
\end{equation}
 We introduce the same notation $\wt\z=\int \z_0(\omega)\,\wt\n(\dd t\, \dd\omega)$ and
$\mathrm{tr}_0(\wt\n)$ for the analogs of $\z$ and $\mathrm{tr}_0(\n)$ when $\n$ is replaced by $\wt \n$.

\begin{lemma}
\label{key-form2}
For any nonnegative measurable functions $F$ and $G$,
$$\E\Big[F(\be^{(r')})\,G(\mathrm{tr}_0(\wt\n))\Big]=\Phi(\frac{3\xi}{2r^2})^{-1}\,
\E\Big[F(\be)\,G(\mathrm{tr}_0(\n))\,\exp\Big(-\frac{3\z}{2r^2}\Big)\Big].$$
\end{lemma}

\proof Suppose that, conditionally on $\be^{(r')}$, 
$\n'$ is distributed as
a Poisson measure on $[0,\xi]\times \S$ with intensity
$$2\,\dd t\,\N_{\sqrt{3}\be^{(r')}_t}(\dd\omega).
$$
Let $\min \n'$ stand for the minimal value of $W_*(\omega)$
for all atoms $(t,\omega)$ of $\n'$. Then, 
conditionally on $\be^{(r')}$, $\wt\n$ is distributed as $\n'$
conditioned to have $\min \n'>-r$. By formula \eqref{hittingpro},
$$\P\Big(\min \n' >-r \,\Big| \,\be^{(r')}\Big)= 
\exp\Big(-2\int_0^\xi \frac{3}{2\times(\sqrt{3}\be^{(r')}_t+r)^2}\,\dd t\Big)
= \exp\Big(-\int_0^\xi \frac{\dd t}{(\be^{(r')}_t+r')^2}\Big).$$
Hence, using the same notation as explained before the lemma to define $\mathrm{tr_0}(\n')$,
$$\E\Big[F(\be^{(r')})\,G(\mathrm{tr}_0(\wt\n))\Big]
=\E\Big[F(\be^{(r')})\,G(\mathrm{tr}_0(\n'))\,\mathbf{1}_{\{\min\n'>-r\}}\,\exp\Big(\int_0^\xi \frac{\dd t}{(\be^{(r')}_t+r')^2}\Big)\Big].$$
By the special Markov property and \eqref{hittingpro}, we have
$$\P\Big(\min \n' >-r \,\Big| \,\be^{(r')}, \mathrm{tr_0}(\n')\Big)=\exp\Big(-\frac{3\z'}{2r^2}\Big),$$
where $\z'$ is the total exit measure at $0$ of the atoms of $\n'$.
We thus arrive at the formula
$$\E\Big[F(\be^{(r')})\,G(\mathrm{tr}_0(\wt\n))\Big]
= \E\Big[F(\be^{(r')})\,G(\mathrm{tr}_0(\n'))\,\exp\Big(-\frac{3\z'}{2r^2}+\int_0^\xi \frac{\dd t}{(\be^{(r')}_t+r')^2}\Big)\Big].$$
At this stage, we use Lemma \ref{key-form} to observe that the law of the pair $(\be^{(r')},\n')$ has a density 
with respect to the law of $(\be,\n)$ which is given by 
$$\Phi(\frac{3\xi}{2r^2})^{-1}\,\exp\Big(-\int_0^\xi \frac{\dd t}{(e(t)+r')^2}\Big),
$$
where $(e(t))_{0\leq s\leq \xi}$ stands for the generic element of $C([0,\xi],\R)$. Thanks to this observation, we get
$$
\E\Big[F(\be^{(r')})\,G(\mathrm{tr}_0(\n'))\,\exp\Big(-\frac{3\z'}{2r^2}+\int_0^\xi \frac{\dd t}{(\be^{(r')}_t+r')^2}\Big)\Big]
=\Phi(\frac{3\xi}{2r^2})^{-1}\,
\E\Big[F(\be)\,G(\mathrm{tr}_0(\n))\,\exp\Big(-\frac{3\z}{2r^2}\Big)\Big],
$$
and this completes the proof. \endproof

Lemma \ref{key-form2} applied with $F=1$ and $G(\mathrm{tr}_0(\wt\n))=\varphi(\wt\z)$, for any test function $\varphi$, implies that
\begin{equation}
\label{relation-laws}
\E[\varphi(\wt\z)]= \Phi(\frac{3\xi}{2r^2})^{-1}\,\E\Big[\varphi(\z)\,\exp\Big(-\frac{3\z}{2 r^2}\Big)\Big].
\end{equation}
Since we know that the density of $\z$ is the function $z\mapsto \frac{3}{2} \xi^{3/2}(\xi+z)^{-5/2}$,
we get that the density of $\wt\z$
is
\begin{equation}\label{density-tilde-Z}
z\mapsto \frac{3}{2}\,\Phi(\frac{3\xi}{2r^2})^{-1}\,\xi^{3/2}(\xi+z)^{-5/2}\,\exp(-\frac{3z}{2r^2}).
\end{equation}

\begin{proposition}
\label{ident-cond-law}
For any nonnegative measurable functions $F$ and $G$, we have for Lebesgue almost every $z>0$,
$$\E\Big[F(\be^{(r')})\,G(\mathrm{tr}_0(\wt\n))\,\Big|\,\wt \z=z\Big]= \E\Big[F(\be)\,G(\mathrm{tr}_0(\n))\,\Big|\,\z=z\Big].$$
In other words, the conditional distributions of $(\be^{(r')},\wt\n)$ given $\wt \z$ and of $(\be,\n)$ given $\z$ are the same.
\end{proposition}

\proof This easily follows from Lemma \ref{key-form2}. Set
$$\gamma_1(z):=\E\Big[F(\be^{(r')})\,G(\mathrm{tr}_0(\wt\n))\,\Big|\,\wt \z=z\Big]\,,\quad\gamma_2(z):=\E\Big[F(\be)\,G(\mathrm{tr}_0(\n))\,\Big|\,\z=z\Big],$$
noting that both $\gamma_1$ and $\gamma_2$ are defined up to a set of zero Lebesgue measure of values of $z>0$. 
By Lemma \ref{key-form2} and formula \eqref{relation-laws},
for any nonnegative measurable function $\varphi$ on $\R_+$,
\begin{align*}
\E\Big[\gamma_1(\wt\z)\,\varphi(\wt\z)\Big]
&=\E\Big[F(\be^{(r')})\,G(\mathrm{tr}_0(\wt\n))\,\varphi(\wt\z)\Big]\\
&=\Phi(\frac{3\xi}{2r^2})^{-1}\,\E\Big[F(\be)\,G(\mathrm{tr}_0(\n))\,\varphi(\z)\,\exp\Big(-\frac{3\z}{2r^2}\Big)\Big]\\
&=\Phi(\frac{3\xi}{2r^2})^{-1}\,\E\Big[\gamma_2(\z)\,\varphi(\z)\,\exp\Big(-\frac{3\z}{2 r^2}\Big)\Big]\\
&=\E[\gamma_2(\wt\z)\,\varphi(\wt\z)].
\end{align*}
Since $\varphi$ was arbitrary it follows that $\gamma_1(z)=\gamma_2(z)$, $\dd z$ a.e. \endproof

\subsection{The complement of a hull in a Brownian disk}\label{sec:compl:hull}
We now recall the construction of the Brownian disk ``viewed from a boundary point'' which is
given in \cite{Repre}. As in \cite{Repre}, we deal with a Brownian disk of perimeter $1$, but the
construction and the results of this section can easily be extended to an arbitrary perimeter $\xi>0$
via scaling arguments.
We start from a pair $(\mathbf{b},\m)$, where 
$\mathbf{b}=(\mathbf{b}_t)_{0\leq t\leq 1}$ is a five-dimensional Bessel bridge from $0$ 
to $0$ over the time interval $[0,1]$ and, conditionally on $\mathbf{b}$, $\m(\dd t\dd\omega)$ is a Poisson measure
on $[0,1]\times \S$ with intensity 
$$2\,\mathbf{1}_{\{W_*(\omega)>0\}}\,\dd t\,\N_{\sqrt{3}\,\mathbf{b}_t}(\dd \omega).$$
 We write
$$\m=\sum_{j\in J} \delta_{(t'_j,\omega'_j)},$$
and $\Sigma':=\sum_{j\in J}\sigma(\omega'_j)$. 
From $\m$, we can define a compact measure metric space $\mathfrak{T}'$ exactly 
in the same way as $\mathfrak{T}$ was defined in \eqref{tree-disk}. We also introduce an
associated clockwise exploration $(\ee'_s)_{s\in[0,\Sigma']}$, and
intervals $[|a,b|]'$ in $\mathfrak{T}'$ are defined as previously from the clockwise exploration.
We specify labels $(\ell'_u)_{u\in\mathfrak{T}'}$ by setting
$\ell'_t:=\sqrt{3}\,\mathbf{b}_t$ for $t\in[0,1]$, and $\ell'_u:=\ell_u(\omega'_j)$ for $u\in\t_{(\omega'_j)}$, $j\in J$.
A fundamental difference is the fact that $\ell'_u\geq 0$ for every $u\in\mathfrak{T}'$
(because by construction $W_*(\omega'_j)>0$ for every $j\in J$).
Furthermore $0$ and $1$ are the only elements of $\mathfrak{T}'$ with zero label.

For every $a,b\in\mathfrak{T}'$ we set
\begin{equation}
\label{formula-Dprime0}
D'^\circ(a,b):=\ell'_a+\ell'_b - 2\max\Big(\min_{c\in[|a,b|]'}\ell'_c,\min_{c\in[|b,a|]'}\ell'_c\Big),
\end{equation}
and then
\begin{equation}\label{def:D^prime}
D'(a,b):=\inf_{a_0=a,a_1,\ldots,a_{p-1},a_p=b}\sum_{i=1}^p D'^\circ(a_{i-1},a_i)
\end{equation}
where the infimum is over all choices of the integer $p\geq 1$
and of the points $a_1,\ldots,a_{p-1}$ in $\mathfrak{T}'$. We notice that $D'(0,1)=D'^\circ(0,1)=0$ and that, for every $a,b\in \mathfrak{T}^\prime$,
\begin{equation}\label{eq:D^prime=0}
D^{\prime}(a,b)=0\quad \text{if and only if }D^{\prime\circ}(a,b)=0.
\end{equation}
This follows from the analogous result \cite[Theorem 13]{Bet} in the Bettinelli-Miermont
construction of Section \ref{Bet-Mie} via an absolute continuity argument.

 Finally, we set $\D'=\mathfrak{T}'/\{D'=0\}$. We observe that $D'(a,b)=0$ implies $\ell'_a=\ell'_b$ 
so that we can make sense of labels on $\D'$, for which we keep the
same notation $\ell'_x$, $x\in\D'$. We write $\Pi'$ for the canonical projection from $\mathfrak{T}'$ onto $\D'$,
and $\mathbf{V}^\prime$ for the pushforward of the volume measure on $\mathfrak{T}'$ under $\Pi'$.

\begin{theorem}{\bf\cite[Theorem 15]{Repre}}
\label{cons-disk}
The quotient space $(\D^\prime, D^\prime, \mathbf{V}^\prime)$ equipped with the distinguished
point $\Pi'(0)$, is a free Brownian disk with perimeter $1$ pointed at a uniform boundary point. 
\end{theorem}

See the end of Section \ref{Bet-Mie}, or \cite[Section 6]{Repre} for the definition of the free Brownian disk pointed at
a uniform boundary point.
The boundary $\partial\D'$ coincides with $\Pi'([0,1])$ (as in formula \eqref{tree-disk}, $[0,1]$ is viewed as a subset of
$\mathfrak{T}'$). In a way similar to the formula $D(\xx_*,u)=\ell_x-\ell_{\xx_*}$
in the Bettinelli-Miermont construction,
labels $\ell'_x$ exactly correspond to distances from the distinguished point $\Pi'(0)$, which lies
on $\partial\D'$. 

We fix $\alpha\in(0,1)$ and we set $\xx_1:=\Pi'(\alpha)$, which is a point of $\partial \D'$ distinct from $\Pi'(0)$. 
Note that $D'(\Pi'(0),\xx_1)=\sqrt{3}\,\bb_\alpha$. We also fix $r>0$ and write $B_r$ for the closed ball of radius $r$
centered at $\Pi'(0)$ in $\D'$.
On the event where $D'(\Pi'(0),\xx_1)>r$, we let $\wh B^{\circ,\xx_1}_r$ denote 
the connected component of $\D^\prime\backslash B_r$ that contains $\xx_1$, 
and define the hull $B^{\bullet,\xx_1}_r:=\D'\backslash \wh B^{\circ,\xx_1}_r$. Notice that
$D'(\Pi'(0),x)=r$ for every $x$ belonging to the topological boundary of $B^{\bullet,\xx_1}_r$. 
We also let $\wh B^{\bullet,\xx_1}_r$ be the closure of $\wh B^{\circ,\xx_1}_r$.

Let us argue on the event where $D'(\Pi'(0),\xx_1)=\sqrt{3}\,\bb_\alpha>r$. On this event, we set $T_-:=\sup\{t\in[0,\alpha]:\bb_t=r/\sqrt{3}\}$ and $T_+:=\inf\{t\in[\alpha,1]:\bb_t=r/\sqrt{3}\}$, so that $(T_-,T_+)$ 
is the excursion interval of $\bb$ above level $r/\sqrt{3}$ that straddles $\alpha$. We also 
set $P_0=T_+-T_-$, and 
$$P_1=\sum_{j\in J:T_-<t'_j<T_+}\z_r(\omega'_j).$$
 \begin{figure}[!h]
 \begin{center}
    \includegraphics[height=8cm,width=8cm]{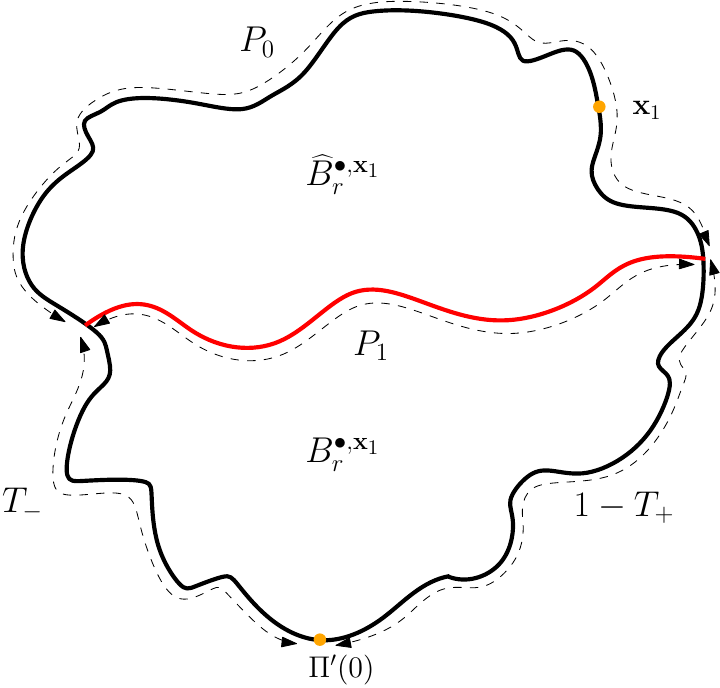} 
 \end{center}
    \caption{Illustration of the Brownian disk $\mathbb{D}^\prime$ and the subsets $B^{\bullet,\xx_1}_r$ and $\wh B^{\bullet,\xx_1}_r$. The intersection $B^{\bullet,\xx_1}_r\cap \wh B^{\bullet,\xx_1}_r$ is represented in red and we interpret $P_1$ as its length. The variables $T_-$, $P_0$ and $1-T_+$ can be thought of  as the lengths of the corresponding subsets of the boundary of $\mathbb{D}^\prime$. }
  \end{figure}
\begin{theorem}
\label{spatial-Markov}
Almost surely under the conditional probability $\P(\cdot\mid D'(\Pi'(0),\xx_1)>r)$, the intrinsic metric 
on $\wh B^{\circ,\xx_1}_r$ has a unique continuous extension to $\wh B^{\bullet,\xx_1}_r$,
which is a metric on $\wh B^{\bullet,\xx_1}_r$. Furthermore, conditionally on the pair $(P_0,P_1)$,  the resulting metric 
space equipped with the volume measure $\wh{\mathbf{V}}'_r$ which is the restriction of $\mathbf{V}'$
to $\wh B^{\bullet,\xx_1}_r$, and with the distinguished point $\Pi^\prime(T_-)$,  is a free Brownian disk with perimeter $P_0+P_1$ pointed at
a uniform boundary point.
\end{theorem}

\proof 
Throughout the proof, we argue under the conditional probability $\P(\cdot\mid D'(\Pi'(0),\xx_1)>r)$. For every $t\in[0,P_0]$,
we set $\bb^\diamond_t=\bb_{T_-+t}-r/\sqrt{3}$. Also define
$$\m^\diamond:=\sum_{j\in J:T_-<t'_j<T_+} \delta_{(t'_j-T_-,\,\omega'_j-r)},$$
with the abuse of notation consisting in writing $\omega-r$  for the snake trajectory $\omega$
shifted by $-r$. We have then
$P_1=\int \z_0(\omega)\,\m^\diamond(\dd t\dd \omega)$.
Let $\mathfrak{T}^\diamond$ be the metric space constructed from $\mathrm{tr}_0(\m^\diamond)$ in the same 
way as $\mathfrak{T}^\star$ was constructed from $\mathrm{tr}_0(\n)$ at the beginning of Section \ref{conf-lim-space}. More precisely, 
$\mathfrak{T}^\diamond$ is obtained from the disjoint union
$$[0,P_0] \cup \Bigg(\bigcup_{j\in J:T_-<t'_j<T_+} \t_{(\mathrm{tr}_0(\omega'_j-r))}\Bigg)$$
by identifying the root of $\t_{(\mathrm{tr}_0(\omega'_j-r))}$ with the point $t'_j-T_-$ of $[0,P_0]$, for every $j\in J$ 
such that $T_-<t'_j<T_+$. 
We also assign labels $\ell^\diamond_a$ to the points of $\mathfrak{T}^\diamond$ in the same 
manner as in Section \ref{conf-lim-space}, so that $\ell^\diamond_a=\sqrt{3}\,\bb^\diamond_a$ if $a\in[0,P_0]$, and 
points of $\t_{(\mathrm{tr}_0(\omega'_j-r))}$ keep their labels. We set $\partial\mathfrak{T}^\diamond:=\{a\in \mathfrak{T}^\diamond : \ell^\diamond_a=0\}$ and
$$\Sigma^\diamond:=\int \sigma(\tr_0(\omega))\, \m^\diamond(\dd t\dd \omega)=\sum_{j\in J:T_-<t'_j<T_+}\sigma\big(\mathrm{tr}_r(\omega'_j)\big). $$

As in Section \ref{conf-lim-space}, we can introduce the clockwise exploration $(\mathcal{E}^\diamond_t)_{0\leq t\leq \Sigma^\diamond}$ of $\mathfrak{T}^\diamond$,
which allows us to define
intervals in $\mathfrak{T}^\diamond$. Then,  
for every $a,b\in\mathfrak{T}^\diamond\backslash \partial\mathfrak{T}^\diamond$, we define the functions $D^\circ_\diamond(a,b)$
and $D_\diamond(a,b)$ by the analogs of formulas \eqref{Dbar1} and \eqref{Dbar2}.  

Let $\xi\in(0,1)$. We observe that, conditionally on $P_0=\xi$, $(\bb_{T_-+t})_{0\leq t\leq \xi}$ is distributed as a five-dimensional Bessel process excursion of
length $\xi$ above level $r/\sqrt{3}$ and thus has the same distribution as the process $(\ov\be^{(r')}_t)_{0\leq t\leq \xi}$
introduced at the beginning of Section \ref{preli-dist}. Hence, conditionally on $P_0=\xi$, $(\bb^\diamond_t)_{0\leq t\leq \xi}$
has the same distribution as $(\be^{(r')}_t)_{0\leq t\leq \xi}$. By construction, conditionally on $P_0=\xi$ and on $\bb^\diamond$,
$\m^\diamond$ is Poisson on $[0,\xi]\times \S$ with intensity
$$2\,\mathbf{1}_{\{W_*(\omega)>-r\}}\,\dd t\,\N_{\sqrt{3}\,\bb^\diamond_t}(\dd\omega).
$$
Comparing with formula \eqref{inten-Pois} for the intensity of the Poisson measure $\wt\n$, and using both identities
$\wt\z=\int\z_0(\omega)\wt\n(\dd t\dd\omega)$ and $P_1=\int \z_0(\omega)\,\m^\diamond(\dd t\dd \omega)$, we get
that the conditional distribution of the pair $(\bb^\diamond,\tr_0(\m^\diamond))$
knowing $(P_0,P_1)=(\xi,z)$ coincides with the conditional distribution of $(\be^{(r')},\mathrm{tr}_0(\wt\n))$ given $\wt\z=z$.
By Proposition \ref{ident-cond-law}, this is also the conditional distribution of $(\be,\mathrm{tr}_0(\n))$ given $\z=z$.

As a consequence of the latter identity in distribution, we can use Proposition \ref{def-lim-space} to get that  the function $(a,b)\mapsto D_\diamond(a,b)$ has 
almost surely a continuous extension to $\mathfrak{T}^\diamond\times \mathfrak{T}^\diamond$. We then consider the quotient metric space $\D_\diamond:=\mathfrak{T}^\diamond/\{D_\diamond=0\}$, and the canonical projection $\Pi_\diamond:\mathfrak{T}^\diamond\la \D_\diamond$.  As usual,  the metric space 
$(\D_\diamond,D_\diamond)$ is equipped with the pushforward of the 
volume measure on $\mathfrak{T}^\diamond$ under $\Pi_\diamond$, which is denoted by  $\mathbf{V}_\diamond$. Moreover, Theorem  \ref{identif-hull-comple} 
implies that, conditionally on the pair $(P_0,P_1)$, the space $(\D_\diamond,D_\diamond,\mathbf{V}_\diamond)$ with the distinguished point $\Pi_\diamond(0)$ is a free Brownian disk with perimeter $P_0+P_1$
pointed at a uniform boundary point. 

In order to complete the proof of  Theorem \ref{spatial-Markov}, we now need to explain that the space $\wh B^{\bullet,x_1}_r$ equipped with
its (extended) intrinsic metric and with the restriction of the volume measure $\mathbf{V}'$ is identified isometrically
to $(\D_\diamond,D_\diamond,\mathbf{V}_\diamond)$, and this identification maps $\Pi'(T_-)$ to $\Pi_\diamond(0)$. 
To this end, we first  introduce the subset of $\mathfrak{T}'$ defined by
$$\mathcal{G}^{\circ,\xx_1}_r:=(T_-,T_+) \cup \Bigg(\bigcup_{j\in J:T_-<t'_j<T_+} \{a\in \t_{(\omega'_j)}:m'_a>r\}\Bigg),$$
where, if $a\in \t_{(\omega'_j)}$, $m'_a$ denotes the minimal label along the ancestral line $\llbracket \rho_{(\omega'_j)},a\rrbracket$, and we 
(of course) make the same identifications as in the definition of $\mathfrak{T}'$. 
Let us verify that $\wh  B^{\circ,\xx_1}_r=\Pi'(\mathcal{G}^{0,\xx_1}_r)$. We know that a point $x$ of $\D'$ belongs to $\wh B^{\circ,\xx_1}_r$
if and only if there is a continuous path from $x$ to $\xx_1$ that does not intersect the ball $B_r$.
From this observation and the cactus bound already used in the proof
of Proposition \ref{ident-lim-space}, it is not hard to verify that 
$x$ belongs to $\wh B^{\circ,\xx_1}_r$ if and only if $x=\Pi'(a)$ where either
$a\in(T_-,T_+)$ or $a$ belongs to one of the trees $\t_{(\omega'_j)}$, with
$T_-<t'_j<T_+$, and labels along the ancestral line of $a$ stay greater than $r$. This leads to the desired identity $\wh  B^{\circ,\xx_1}_r=\Pi'(\mathcal{G}^{0,\xx_1}_r)$.

Next,  set
$$\partial \mathcal{G}^{\circ,\xx_1}_r:=\{T_-,T_+\} \cup \Bigg(\bigcup_{j\in J:T_-<t'_j<T_+} \big\{a\in \t_{(\omega'_j)}: \ell^\prime_a=r\text{ and }\ell'_b>r \text{ for every } b\in \llbracket \rho_{(\omega'_j)},a\rrbracket\backslash\{a\}\big\}\Bigg).$$
One easily verifies that  the topological boundary of $\wh  B^{\circ,\xx_1}_r$ is $\Pi'(\partial \mathcal{G}^{\circ,\xx_1}_r)$. Consequently, $\wh B^{\bullet,\xx_1}_r=\Pi'(\mathcal{G}^{\bullet,\xx_1}_r)$ where $\mathcal{G}^{\bullet,\xx_1}_r=\mathcal{G}^{\circ,\xx_1}_r\cup \partial \mathcal{G}^{\circ,\xx_1}_r$.  Note that we have 
$$ \mathcal{G}^{\bullet,\xx_1}_r=[T_-,T_+] \cup \Bigg(\bigcup_{j\in J:T_-<t'_j<T_+}  \mathcal{T}_{(\mathrm{tr}_r(\omega'_j))}\Bigg),$$
where as usual we identify $\mathcal{T}_{(\mathrm{tr}_r(\omega'_j))}$ with a subset of $\mathcal{T}_{(\omega'_j)}$.

We can then
identify $\mathcal{G}^{\bullet,\xx_1}_r$ with $\mathfrak{T}^\diamond$ in the 
following manner: $a\in[T_-,T_+]$
is identified with $a-T_-$, and a point $a\in \mathcal{T}_{(\mathrm{tr}_r(\omega'_j))}$ (where $j$ is such that $T_-< t_j^\prime<T_+$) is identified to the corresponding point of $\t_{(\mathrm{tr}_0(\omega'_j-r))}$. Moreover, using \eqref{eq:D^prime=0} and Proposition \ref{def-lim-space}, one checks that, for every $a,b\in \mathcal{G}^{\bullet,\xx_1}_r$, the 
property $D^\prime(a,b)=0$ holds if and only if the points $\wt a$ and $\wt b$ of $\mathfrak{T}^\diamond$ corresponding to $a$ and $b$ 
satisfy $D_\diamond(\wt a,\wt b)=0$. This leads to the desired identification of the sets $\wh B^{\bullet,\xx_1}_r$ and $\D_\diamond$. 
Then one verifies that the intrinsic metric on $\wh  B^{\circ,\xx_1}_r$ coincides (modulo the preceding identification)
with the restriction of the metric $D_\diamond$ to $\D_\diamond\backslash \partial\D_\diamond$. This relies on arguments very similar 
to the proof of Proposition \ref{ident-lim-space}, and we omit the details. Finally, it is immediate that the restriction of 
$\mathbf{V}'$ to $\wh B^{\bullet,\xx_1}_r$ corresponds to the volume measure on $\D_\diamond$. This completes the proof. \endproof

We will write $\wh D'_r$ for the metric on $\wh B^{\bullet,\xx_1}_r$ constructed in Theorem \ref{spatial-Markov} as the extension of the 
intrinsic metric on $\wh B^{\circ,\xx_1}_r$. In view of future applications, it will also be convenient to
introduce the standard boundary curve of the Brownian disk $\wh B^{\bullet,\xx_1}_r$ that is defined as follows. 
Recall the construction of the boundary curve $(\Gamma^\star(t))_{0\leq t\leq \xi+\z_\U}$
from the pair $(\be,\tr_0(\n))$ as explained before Theorem \ref{identif-hull-comple}. Since the conditional distribution of the pair
$(\bb^\diamond,\tr_0(\m^\diamond))$
knowing $(P_0,P_1)=(\xi,z)$ coincides with the conditional distribution of $(\be,\mathrm{tr}_0(\n))$ given $\z=z$, we
can use the same construction to get a standard boundary curve of the Brownian disk $\D_\diamond$, hence (via the
identification of the preceding proof) a standard boundary curve $(\wh\Gamma'_r(t))_{0\leq t\leq P_0+P_1}$
of the Brownian disk $\wh B^{\bullet,\xx_1}_r$. 
More precisely, $\wh \Gamma'_r(t)=\Pi'(T_-+t)$ for $0\leq t\leq P_0$, and the values of 
$\wh\Gamma'_r(t)$ for $P_0\leq t\leq P_0+P_1$ are defined by the analog of formula \eqref{gamma-star}
(cf.~formula \eqref{ident-bdry-curve} below).
We have then $\{\wh\Gamma'_r(t):0\leq t\leq P_0\}=
\wh B^{\bullet,\xx_1}_r\cap \partial\D'$, and the set $\{\wh\Gamma'_r(t):P_0\leq t\leq P_0+P_1\}$ is the 
topological boundary of $\wh B^{\bullet,\xx_1}_r$ or, equivalently, of the hull $B^{\bullet,\xx_1}_r$.

\subsection{A spatial Markov property}

Our goal in this section is to prove that the free Brownian disk $(\wh B^{\bullet,\xx_1}_r,\wh D'_r,\wh{\mathbf{V}}'_r,\Pi^\prime(T_{-}))$ in Theorem \ref{spatial-Markov} is independent
of the hull $B^{\bullet,\xx_1}_r$ conditionally on the pair $(P_0,P_1)$. To make this assertion precise, we need to explain 
how the hull $B^{\bullet,\xx_1}_r$ is viewed as a random measure metric space. We argue on the event $\{D'(0,\xx_1)>r\}$
and we keep the notation of the previous section. We introduce the subset $\k_r$ of $\mathfrak{T}'$
defined by
$$\k_r= [0,T_-]\cup[T_+,1]\cup  \Big(\bigcup_{j\in J:\,t'_j\in[0,T_-]\cup[T_+,1]}  \t_{(\omega'_j)}\Big)\cup\Big(\bigcup_{j\in J:\,T_-<t'_j<T_+} \{a\in \t_{(\omega'_j)}:m'_a\leq r\}\Big),$$
where we recall that $m_a^\prime$ stands for the minimal label along the ancestral line of $a\in \t_{(\omega'_j)}$. Note that labels on $\mathfrak{T}'\backslash \k_r$
are greater than $r$. We have $\k_r=\mathfrak{T}'\backslash \mathcal{G}^{\circ,\xx_1}_r$,
and therefore $\Pi'(\k_r)=B^{\bullet,\xx_1}_r$ as a consequence of the equality 
$\wh  B^{\circ,\xx_1}_r=\Pi'(\mathcal{G}^{0,\xx_1}_r)$.   In view of forthcoming applications, we also mention the following simple fact. 
Let $a,b\in \k_r$. Then, in formula \eqref{formula-Dprime0} defining $D'^\circ(a,b)$, we may replace the intervals $[|a,b|]'$
and $[|b,a|]'$ by $[|a,b|]'\cap\k_r$
and $[|b,a|]'\cap\k_r$ respectively: the point is that, if the interval $[|a,b|]'$ contains a point $c\notin \k_r$,
then, necessarily, it contains another point $c'$ (belonging to $\k_r$) whose label is $r$ and is thus smaller than the label of $c$. 
Informally, the definition of $D'^\circ(a,b)$, when $a,b\in \k_r$ only depends on the labels on $\k_r$, despite
the fact that the interval $[|a,b|]'$ may not be contained in $\k_r$. 

For every $a,b\in\k_r$, we set
\begin{equation}
\label{pseudo-hull2}
D^{\prime}_r(a,b) := \mathop{\inf \limits_{a_0,a_1,\ldots,a_p\in \k_r}}_{a_0=a,\,a_p=b} \sum_{i=1}^p D^{\prime\circ}(a_{i-1},a_i), 
\end{equation}
where the infimum is over all choices of the integer $p\geq 1$ and of the
finite sequence $a_0,a_1,\ldots,a_p$ in $\k_r$ such that $a_0=a$ and
$a_p=b$. This is similar to the definition \eqref{def:D^prime} of $D'(a,b)$, but we 
restrict the infimum to ``intermediate'' points $a_1,\ldots,a_{p-1}$ that belong to $\k_r$. 
Clearly, we have $D^\prime(a,b)\leq D^\prime_r(a,b)\leq D^{\prime\circ}(a,b)$  for every $a,b\in \k_r$. 
Since  the condition $D^\prime(a,b)=0$ can only hold if $D^{\prime\circ}(a,b)=0$, we get that, for every $a,b\in \k_r$, we have $D^\prime_r(a,b)=0$ if and only if $D'(a,b)=0$. Hence $D^\prime_r$
induces a metric on $\Pi^\prime(\k_r)=B^{\bullet,\xx_1}_r$ and we keep the notation $D^\prime_r$ for this metric. Using simple geodesics (as defined at the end 
of Section \ref{Bet-Mie}) and  the definition 
\eqref{def:D^prime} of $D^\prime$ as an infimum, one verifies that the restriction of $D^\prime_r$ 
to the interior of $B^{\bullet,\xx_1}_r$ coincides with the intrinsic distance induced by $D^\prime$. This follows by an adaptation of the proof of Proposition \ref{ident-lim-space}
and we omit the details since this is not really needed in what follows. Additionally, we have $D'_r(\Pi'(0),x)=D'(\Pi'(0),x)$ for every $x\in B^{\bullet,\xx_1}_r$ (note that a $D'$-geodesic from $x$ to $\Pi'(0)$ 
cannot exit the hull $B^{\bullet,\xx_1}_r$). In particular, $D'_r(\Pi'(0),\wh\Gamma'_r(t))=D'(\Pi'(0),\wh\Gamma'_r(t))=r$ for every 
$t\in[P_0,P_0+P_1]$.

We equip the metric space $(B^{\bullet,\xx_1}_r,D^\prime_r)$ with the restriction of the volume measure $\mathbf{V}^{\prime}$, which we  denote by $\mathbf{V}^{\prime}_r$.  It is also convenient to introduce a boundary curve of $B^{\bullet,\xx_1}_r$, which we define as follows. We set,
for every $t\in[0,P_1+1-P_0]$,
$$\Gamma_r^\prime(t):=\left\{\begin{array}{ll}
\Pi'(t)\quad&\hbox{if } t\in [0,T_-],\\
\wh\Gamma'_r(P_0+P_1-(t-T_-))\quad&\hbox{if } t\in (T_-,T_-+P_1),\\
\Pi'(t-(T_-+P_1)+T_+)\quad&\hbox{if } t\in [T_-+P_1,P_1+1-P_0].
\end{array}
\right.
$$

Note that $\Gamma'_r$ is a simple loop taking values in $B^{\bullet,\xx_1}_r$, and $\Gamma'_r(0)=\Gamma'_r(P_1+1-P_0)=\Pi'(0)$. In fact, using Jordan's theorem, it is
not hard to to verify that $B^{\bullet,\xx_1}_r$ has the topology of the closed unit disk, which makes it possible to consider the ``boundary'' of $B^{\bullet,\xx_1}_r$, and this
boundary (which is not the topological boundary) is precisely the range of $\Gamma'_r$. We will consider the hull $(B^{\bullet,\xx_1}_r,D^\prime_r,\mathbf{V}'_r)$ as equipped with the curve $\Gamma'_r$: we view
the 4-tuple
$$(B^{\bullet,\xx_1}_r,D^\prime_r,\mathbf{V}'_r,\Gamma'_r)$$
as a random variable taking values in the space $\M^{GHPU}$. Then it is not hard to verify that the quantities $T_-,T_+,P_0,P_1$
are measurable functions of the latter random space. In fact, recalling that $D'_r(\Pi'(0), \Pi'(t))=D'(\Pi'(0), \Pi'(t))=\sqrt{3}\,\bb_t$ 
for $t\in[0,T_-]\cup[T_+,1]$, one sees that $T_-$ is the first time $t\geq 0$ 
such that there exists $\ve>0$ verifying  
$$D^{\prime}_r\big(\Pi'(0),\Gamma^{\prime}_r(t+s)\big)=r,\quad \forall s\in[0,\ve],$$
and an analogous representation holds for $T_+$. Furthermore $P_1=\inf\{t\geq 0:D^{\prime}_r(0,\Gamma^{\prime}_r(T_-+t))\neq r\}$.

\begin{theorem}\label{indep-hull}
Under the conditional probability $\P(\cdot\mid D'(0,\xx_1)>r)$, the space $(\wh B^{\bullet,\xx_1}_r,\wh D^\prime_r,\wh{\mathbf{V}}'_r,\wh\Gamma'_r)$
is independent of $(B^{\bullet,\xx_1}_r,D^\prime_r,\mathbf{V}'_r,\Gamma'_r)$ conditionally on the pair $(P_0,P_1)$.
\end{theorem}

\proof
The general strategy of the proof is to show that the space $(B^{\bullet,\xx_1}_r,D^\prime_r,\mathbf{V}'_r,\Gamma'_r)$ can be constructed  from random quantities that are independent of 
$(\wh B^{\bullet,\xx_1}_r,\wh D^\prime_r,\wh{\mathbf{V}}'_r,\wh\Gamma'_r)$ conditionally on $(P_0,P_1)$. To this end, we will describe the hull $B^{\bullet,\xx_1}_r$ in terms of the labeled tree $\mathfrak{T}^\prime$. 
Recall the notation
$(\ee'_s)_{s\in[0,\Sigma']}$ for the clockwise exploration of $\mathfrak{T}'$. For every $j\in J$ such that $T_-<t'_j<T_+$, write
$[u_j,u_j+\sigma(\omega'_j)]$ for the interval corresponding to $\omega'_j$ in the time scale of the clockwise 
exploration $\ee'$ (meaning that $\ee'_s\in\t_{(\omega'_j)}$ if and only if $s\in[u_j,u_j+\sigma(\omega'_j)]$), and recall the notation $(L^r_t(\omega'_j))_{t\geq 0}$
for the exit local time of $\omega'_j$ from $(r,\infty)$. We set, for every $t\geq 0$,
\begin{equation}
\label{def-loc-ti}
L'_t:=\sum_{j\in J: T_-<t'_j<T_+} L^r_{(t-u_j)^+}(\omega'_j),
\end{equation}
Then, for every index $j$ such that $T_-<t'_j<T_+$, we
denote the excursions of $\omega'_j$ outside $(r,\infty)$ by $(\omega^\#_{j,k})_{k\in I_{j}}$
where $I_j$ is an appropriate indexing set (if the index $j$ is such that $W_*(\omega'_j)\geq r$,
$I_j$ is the empty set). In the time scale of 
$\mathcal{E}^\prime$, the excursion $\omega^\#_{j,k}$ corresponds to an interval
$(\alpha_{j,k},\beta_{j,k})$ and we set $t^\#_{j,k}=L^\prime_{\alpha_{j,k}}=L^\prime_{\beta_{j,k}}$.
Finally,  we introduce the point measure 
$$\n_\#(\dd t\dd\omega):=\sum_{j\in J:T_-<t'_j<T_+}\,\sum_{k\in I_j} \delta_{(t^\#_{j,k},\omega^\#_{j,k})}(\dd t\dd\omega).$$
We will show that the hull $B^{\bullet,\xx_1}_r$ is a function of the triple
\begin{equation}
\label{triple}
\mathbb{T}_\#:=\Bigg(\n_\#(\dd t\dd\omega);\Big((\bb_t)_{0\leq t\leq T_-},\sum_{j\in J:0\leq t'_j\leq T_-}\delta_{(t'_j,\omega'_j)}\Big); \Big((\bb_t)_{T_+\leq t\leq 1},\sum_{j\in J:T_+\leq t'_j\leq 1}\delta_{(t'_j,\omega'_j)}\Big)\Bigg).
\end{equation}
To this end, we consider the interval
$[0,1-P_0+P_1]$, viewed as the union of the three intervals $[0,T_-]$, $[T_-,T_-+P_1]$
and $[T_-+P_1,T_-+P_1+1-T_+]$. We define 
$\mathfrak{T}^\#$ as the union 
$$[0,1-P_0+P_1] \cup \Bigg(\bigcup_{j\in J: t'_j\in[0,T_-]\cup[T_+,1]} \t_{(\omega'_j)}\Bigg)
\cup\Bigg(\bigcup_{j\in J:T_-<t'_j<T_+}\bigcup_{k\in I_j} \t_{(\omega^\#_{j,k})}\Bigg),$$
 where, for $j\in J$,
 \begin{itemize}
 \item if $0\leq t'_j\leq T_-$, the root 
 of $\t_{(\omega'_j)}$ is identified to $t'_j\in[0,T_-]$;
 \item if $T_-<t'_j<T_+$, then, for every $k\in I_j$, the root of $\t_{(\omega^\#_{j,k})}$ is identified to 
 $T_- + t^\#_{j,k}\in [T_-,T_-+P_1]$;
 \item if $T_+\leq  t'_j\leq 1$, the root 
 of $\t_{(\omega'_j)}$ is identified to $T_-+P_1+(t'_j-T_+)\in [T_-+P_1,1-P_0+P_1]$.
 \end{itemize}

As in the Bettinelli-Miermont construction of Section \ref{Bet-Mie}, we view $\mathfrak{T}^\#$ as a measure metric space
and we write $\Sigma^\#$ for the total mass of its volume measure. We can also assign labels to the points
of $\mathfrak{T}^\#$: points of the trees $\t_{(\omega'_j)}$ and $\t_{(\omega^\#_{j,k})}$ obviously keep their labels,
the label of $s\in [0,T_-]$, resp. of $s\in [T_-+P_1,T_-+P_1+1-T_+]$, is $\sqrt{3}\,\bb_s$,
resp. $\sqrt{3}\,\bb_{T_++s-(T_-+P_1)}$, and finally the label  of each $s\in [T_-,T_-+P_1]$
is $r$. 
We define the exploration
function $(\ee^\#_t)_{0\leq t\leq \Sigma^\#}$ in a way similar to Section \ref{Bet-Mie}: informally, we concatenate the 
exploration functions of the trees $\t_{(\omega'_j)}$ and $\t_{(\omega^\#_{j,k})}$ in the order prescribed 
by their roots viewed as elements of $[0,1-P_0+P_1]$. This exploration function allows us to define intervals on $\mathfrak{T}^\#$,
and then to introduce the functions $D_\#^\circ(a,b)$ and $D_\#(a,b)$ for $a,b\in\mathfrak{T}^\#$,
exactly as we did in Section~\ref{Bet-Mie} to define $D^\circ(a,b)$ and $D(a,b)$. Similarly, we
consider the quotient space $\D_\#:=\mathfrak{T}^\#/\{D_\#=0\}$ and we write $\Pi_\#$ for the canonical projection. 
 As usual, we write $\mathbf{V}_\#$ for the pushforward of the 
volume measure on $\mathfrak{T}^\#$ under $\Pi_\#$. Finally, we set $\Gamma_\#(t):=\Pi_\#(t)$, for $t\in[0,1-P_0+P_1]$.

We then claim that we have the almost sure equality
\begin{equation}
\label{identif-MMS}
(B^{\bullet,\xx_1}_r,D^{\prime}_r,\mathbf{V}^{\prime}_r,\Gamma'_r)= (\D_\#,D_\#,\mathbf{V}_\#,\Gamma_\#).
\end{equation}
Let us explain why Theorem \ref{indep-hull} follows from \eqref{identif-MMS}. On one hand, 
we know that the space $(\D_\#,D_\#,\mathbf{V}_\#,\Gamma_\#)$ is obtained as a function of the triple $\mathbb{T}_\#$
in \eqref{triple}. On the other hand, the proof of Theorem \ref{spatial-Markov} shows that the space 
$(\wh B^{\bullet,\xx_1}_r,\wh D^\prime_r,\wh{\mathbf{V}}'_r,\wh\Gamma'_r)$ is a function of the pair
$(\bb^\diamond,\tr_0(\m^\diamond))$ introduced in this proof. So the statement of Theorem \ref{indep-hull} reduces to
checking that
 the pair  $(\bb^\diamond,\tr_0(\mathcal{M}^\diamond))$  is independent of  $\mathbb{T}_\#$ conditionally on $(P_0,P_1)$. 
To this end,  notice that  the excursion $\bb^\diamond$
is independent of $\big((\bb_t)_{0\leq t\leq T_-},(\bb_t)_{T_+\leq t\leq 1}\big)$
conditionally on $P_0$. It follows that,  conditionally on $P_0$, the pair $(\bb^\diamond,\mathcal{M}^\diamond)$ is independent of the pair
$$\Bigg((\bb_t)_{0\leq t\leq T_-},\sum_{j\in J:0\leq t'_j\leq T_-}\delta_{(t'_j,\omega'_j)}\Bigg)
\;,\quad \Bigg((\bb_t)_{T_+\leq t\leq 1},\sum_{j\in J:T_+\leq t'_j\leq 1}\delta_{(t'_j,\omega'_j)}\Bigg).$$
Furthermore, an application of the special Markov property  \eqref{measure:special:mark} entails that, 
conditionally on $(P_0,P_1)$, the point measure $\n_\#(\dd t\dd\omega)$
is Poisson with intensity $\mathbf{1}_{[0,P_1]}(t)\,\dd t\,\N_r(\dd \omega\cap\{W_*(\omega)>0\})$,
and is independent of the pair $(\bb^\diamond,\tr_0(\m^\diamond))$. It follows that conditionally on $(P_0,P_1)$, the variable  $(\bb^\diamond,\tr_0(\mathcal{M}^\diamond))$  is independent of   $\mathbb{T}_\#$, which was the desired result.

It only remains to justify our claim \eqref{identif-MMS}. This relies on arguments similar
to the proof of Theorem 31 in \cite{spine}, which is a statement analogous to Theorem \ref{indep-hull} 
for the hull centered at the distinguished point of the Brownian plane. For this reason, we will skip some details.
Recall the definition \eqref{def-loc-ti} of the process $(L'_t)_{t\geq 0}$, and, for every $s\in[0,P_1)$, set 
$$\kappa'(t):=\inf\{s\geq 0:L'_s>t\}.$$
By convention we let $\kappa'(P_1)=\kappa'(P_1-)$ be the left limit of $t\mapsto \kappa'(t)$ at $P_1$.
We have then $\ee'_{\kappa'(0)}=T_-$ and $\ee'_{\kappa'(P_1)}=T_+$. From the construction of
the boundary curve $\wh\Gamma'_r$, we have 
\begin{equation}
\label{ident-bdry-curve}
\Gamma'_r(T_-+s)=\wh\Gamma'_r(P_0+P_1-s)=\Pi'(\ee'_{\kappa'(s)})\quad\hbox{for every }s\in [0,P_1].
\end{equation}

We then define a mapping $\Phi:\k_r\la\mathfrak{T}^\#$ by the following prescriptions. Let $a\in\k_r$:
\begin{itemize}
\item[(i)] If $a\in[0,T_-]$, $\Phi(a)$ is the ``same'' point of $\mathfrak{T}^\#$.
\item[(ii)] If $a\in[T_+,1]$, $\Phi(a)$  is the point $a+P_1-P_0$ of $\mathfrak{T}^\#$.
\item[(iii)] If $a\in \mathcal{T}_{(\omega_j^\prime)}$, for some $j\in J$ with $t_j^\prime \in [0,T_-]\cup [T_+,1]$, $\Phi(a)$ is the ``same'' point in $\mathfrak{T}^\#$.
\item[(iv)] If $a\in \t_{(\omega^\#_{j,k})}$, for some $j\in J$ such that $t'_j\in(T_-,T_+)$ and some $k\in I_j$, $\Phi(a)$ is the ``same'' point in $\mathfrak{T}^\#$ ---
we use the fact that $\t_{(\omega^\#_{j,k})}$ can be viewed as a subset of the tree $\t_{(\omega'_j)}$.
\item[(v)]  If $a$ is of the form $\ee'_{\kappa'(s)}$ with $s\in[0,P_1]$, or of the form $\ee'_{\kappa'(s-)}$ with $s\in (0,P_1]$, $\Phi(a)$ is the point $T_-+s\in \mathfrak{T}^\#$.
\end{itemize} 
 One verifies that these
prescriptions are consistent with the identifications made when defining $\k_r$ and $\mathfrak{T}^\#$. In particular, for $j\in J$ such that $t'_j\in(T_-,T_+)$, and $k\in I_j$,
the root $\rho_{(\omega^\#_{j,k})}$ of $\t_{(\omega^\#_{j,k})}$ (viewed as an element of $\k_r$) is easily seen 
to coincide with $\ee'_{\kappa'(t^\#_{j,k})}$ and thus (by property (v)) is mapped to $T_-+t^\#_{j,k}$, which is identified to $\rho_{(\omega^\#_{j,k})}$
in $\mathfrak{T}^\#$. Moreover, (i) --- (v) define $\Phi(a)$ for {\it every} $a\in\k_r$. The point is that, if $a\in\k_r$
belongs to a tree $\t_{(\omega'_j)}$, for some $j$ such that $T_-<t'_j<T_+$, and if $a$ does not belong to any of the subtrees $\t^\#_{j,k}$ 
with $k\in I_j$, then necessarily $\ell'_a=m'_a=r$, and the support property of the exit
local time ensures that we have  $a=\ee'_{\kappa'(s)}$ or $a=\ee'_{\kappa'(s-)}$ for some $s\in[0,P_1]$.

We note that $\Phi$ preserves labels
and is surjective. However, $\Phi$ is not injective because $\Phi(\ee'_{\kappa'(s)})=\Phi(\ee'_{\kappa'(s-)})$ for $s\in(0,P_1]$. Nonetheless,
it follows from our definitions and the support property of the exit local time that $\Pi'(\ee'_{\kappa'(s)})=\Pi'(\ee'_{\kappa'(s-)})$ for every $s\in(0,P_1]$.
We also note that, if $a,b\in\k_r$, the image under $\Phi$
of $[|a,b|]'\cap\k_r$ is the ``interval'' from $\Phi(a)$ to $\Phi(b)$ in $\mathfrak{T}^\#$. 
Using this observation and the simple fact stated before formula \eqref{pseudo-hull2}, it is 
straightforward to verify that,  for every $a,b\in \k_r$, we have $D'^\circ(a,b)=D_\#^\circ(\Phi(a),\Phi(b))$, and then $D^\prime_r(a,b)=D_\#(\Phi(a),\Phi(b))$. It follows that $\Phi$ induces an isometry from 
$B^{\bullet,\xx_1}_r=\Pi^\prime(\k_r)$
onto $\D_\#$, and this isometry, for which we keep the same notation $\Phi$, is easily seen to preserve the volume measures. 
Finally, from \eqref{ident-bdry-curve} and properties (i),(ii),(v), we immediately get that $\Phi(\Gamma'_r(s))=\Pi_\#(s)=\Gamma_\#(s)$
for every $s\in[0,1-P_0+P_1]$. This completes the proof of \eqref{identif-MMS} and of Theorem \ref{indep-hull}.
\endproof

Theorems  \ref{spatial-Markov} and \ref{indep-hull} should be interpreted as giving a way to define a peeling exploration of the
Brownian disk $\D'$. Starting from the point $\Pi'(0)$ uniformly distributed on the boundary, one may start by ``peeling'' the hull
$B^{\bullet,\xx_1}_r$
of small radius $r>0$
centered at this point (relative to another fixed point $\xx_1$ of the boundary). Then the remaining part $\wh B^{\bullet,x_1}_r$ of the initial Brownian disk
is a Brownian disk with a different perimeter (Theorem \ref{spatial-Markov}) and, if we choose a point on its boundary as a function of the
part that has been removed, this point will again be distributed uniformly on the boundary of the new Brownian disk (as a
consequence of the independence property in Theorem \ref{indep-hull}). In the notation of this section, we may choose
the next point to be ``peeled'' as $\wh\Gamma'_r(U)$ where $U$ is a measurable function of 
$(B^{\bullet,\xx_1}_r,D^\prime_r,\mathbf{V}'_r,\Gamma'_r)$.

\section*{Appendix}
 In this appendix, we briefly explain how formula \eqref{triangu21} is derived from \cite{Kri}. We consider the
 generating function
 $$G(y,z)= \sum_{k=0}^\infty\sum_{L=1}^\infty \sum_{p=1}^\infty \#\T^2(L,p,k)\,(12\sqrt{3})^{-k}\,y^L\,z^p.$$
 A direct application of formula (27) in \cite{Kri} gives
 $$G(y,z)=\int_0^3  \frac{2sy}{(1-4sy)^{3/2}}\frac{2sz}{(1-4sz)^{3/2}}\,\frac{\dd s}{2s}.$$
Note that the parameter $k$ corresponding to the number of inner vertices is replaced in \cite{Kri} by a 
parameter counting the number of edges of the triangulation, but of course this makes 
no difference thanks to Euler's formula. 

Let us also consider the generating function
$$\wt G(y,z)=\sum_{L= 1}^\infty\sum_{p=1}^\infty \frac{1}{2}\,\frac{3^{L+p}}{L+p}\,\,L {2L\choose L}\,p{2p\choose p}\,y^L\,z^p.$$
If we set $F_{y,z}(t)=\wt G(\frac{ty}{3},\frac{tz}{3})$, we have
$$F_{y,z}(t)= \frac{1}{2} \sum_{L= 1}^\infty\sum_{p=1}^\infty \frac{1}{L+p}\,\,L {2L\choose L}\,p{2p\choose p}\,t^{L+p}\,y^L\,z^p$$
and 
$$tF'_{y,z}(t)= \frac{1}{2} \sum_{L= 1}^\infty\sum_{p=1}^\infty L {2L\choose L}\,p{2p\choose p}\,t^{L+p}\,y^L\,z^p = \frac{1}{2}\,\varphi(ty)\,\varphi(tz),$$
where
$$\varphi (x)=\sum_{n=1}^\infty n{2n\choose n} x^n= \frac{2x}{(1-4x)^{3/2}}.$$
We conclude that
$$\wt G(y,z)=F_{y,z}(3)= \int_0^3\varphi(ty)\,\varphi(tz)\,\frac{\dd t}{2t} = G(y,z)$$
as desired.

\begin{funding}
This work was supported by the ERC Advanced Grant 740943 {\sc GeoBrown}.
\end{funding}


\begin{thebibliography}{99}


\bibitem{ALG}
{\sc C. Abraham, J.-F. Le Gall}, Excursion theory for Brownian motion
indexed by the Brownian tree. {\it J. Eur. Math. Soc. (JEMS)} 20, 2951--3016 (2018)


\bibitem{AHS}
{\sc M. Albenque, N. Holden, X. Sun},
Scaling limit of large triangulations of polygons. 
{\it Electron. J. Probab.} 25, Paper No. 135, 43 pp. (2020)

\bibitem{AC}
{\sc O. Angel, N. Curien},
Percolations on random maps I: Half-plane models.
{\it Ann. Inst. H. Poincar\'e Probab. Statist.} 51, 405--431 (2015)

\bibitem{BMR}
{\sc E. Baur, G. Miermont, G. Ray},
Classification of scaling limits of uniform quadrangulations with a boundary.
{\it Ann. Probab.}
47, 3397--3477 (2019)

\bibitem{BF}
{\sc O. Bernardi, \'E. Fusy}, Bijections for planar maps
with boundaries. {\it J. Combin. Theory Ser. A} 158, 176--227 (2018)

\bibitem{Bet}
{\sc J. Bettinelli}, Scaling limit of random planar 
quadrangulations with a boundary. {\it Ann. Inst. H. Poincar\'e Probab. Stat.} 51, 432--477 (2015)

\bibitem{BM}
{\sc J. Bettinelli, G. Miermont},
Compact Brownian surfaces I. Brownian disks. 
{\it Probab. Theory Related Fields} 167, 555-614 (2017)

\bibitem{BBI}
{\sc D. Burago, Y. Burago, S. Ivanov}, {\it A Course in Metric Geometry}.
Graduate Studies in Mathematics, vol. 33. Amer. Math. Soc., Boston, 2001.


\bibitem{CC}
{\sc A. Caraceni, N. Curien},
Geometry of the Uniform Infinite Half-Planar Quadrangulation.
{\it Random Struct. Alg.} 52, 454--494 (2018)


\bibitem{Cur}
{\sc N. Curien},
{\it Peeling Random Planar Maps}. Lecture notes from the 2019 Saint-Flour Probability
Summer School. {Lecture Notes in Mathematics} 2335. Springer, Berlin, 2023.

\bibitem{CLG}
{\sc N. Curien, J.-F. Le Gall},
Scaling limits for the peeling process on random maps.
{\it Ann. Inst. H. Poincar\'e Probab. Stat.} 53, 322--357 (2017)


\bibitem{GM0}
{\sc E. Gwynne, J. Miller}, Scaling limit of the uniform infinite half-plane quadrangulation in
the Gromov-Hausdorff-Prokhorov-uniform topology.
{\it Electron. J. Probab.} 22, Paper No. 84, 47 pp. (2017)

\bibitem{GM1}
{\sc E. Gwynne, J. Miller}, 
Convergence of the free Boltzmann quadrangulation with simple boundary to the Brownian disk.
{\it Ann. Inst. Henri Poincar\'e Probab. Stat.} 55, 551--589 (2019)

\bibitem{Kri}
{\sc M. Krikun},
Explicit enumeration of triangulations with multiple boundaries.
{\it Electron. J. Combin.} 14, Research Paper 61, 14 pp. (2007)

\bibitem{Zurich} {\sc J.-F. Le Gall}, {\it Spatial Branching Processes, Random Snakes and 
Partial Differential Equations}. {Lectures in Mathematics ETH Z\"urich}. Birkh\"auser, Boston, 1999.

\bibitem{Geodesics}
{\sc J.-F. Le Gall}, Geodesics in large planar maps and
in the Brownian map. {\it Acta Mathematica} 205, 287--360 (2010)


\bibitem{subor} 
{\sc J.-F. Le Gall}, Subordination of trees and the Brownian map.
{\it Probab. Theory Related Fields} {171}, 819--864 (2018)


\bibitem{Disks}
{\sc J.-F. Le Gall}, Brownian disks and the Brownian snake.
{\it Ann. Inst. H. Poincar\'e Probab. Stat.} 55, 237--313 (2019)

\bibitem{Repre}
{\sc J.-F. Le Gall}, The Brownian disk viewed from a boundary point.
{\it Ann. Inst. H. Poincar\'e Probab. Stat.} 58, 1091-1119 (2022)

\bibitem{Stars}
{\sc J.-F. Le Gall}, Geodesic stars in random geometry. {\it Ann. Probab.} 50, 1013--1058 (2022)

\bibitem{Hausdorff}
{\sc J.-F. Le Gall}, The volume measure of the Brownian sphere is
a Hausdorff measure. {\it Electron. J. Probab.} 27, article no.113, 1--28 (2022)

\bibitem{spine}
{\sc J.-F. Le Gall, A. Riera}, Spine representations for non-compact models of 
random geometry. {\it Probab. Theory Related Fields} 181, 571--645 (2021)

\bibitem{ALG3}
{\sc J.-F. Le Gall, A. Riera}, Peeling the Brownian half-plane. In preparation.

\bibitem{MS}
{\sc J. Miller, S. Sheffield}, An axiomatic characterization of the Brownian map.
J. \'Ec. polytech. Math. 8, 609--731 (2021)

\bibitem{MS2}
{\sc J. Miller, S. Sheffield},
Liouville quantum gravity and the Brownian map II: Geodesics and continuity of the embedding. 
Ann. Probab. 49, 273--2829 (2021)
\end{thebibliography}
\end{document}